\documentclass[11pt]{article}
\pdfoutput=1
\usepackage{jheppub}

\usepackage{blindtext}

\usepackage{mathrsfs}
\usepackage{bm}
\usepackage{paralist}
\usepackage{url}
\usepackage[]{microtype}

\usepackage{tikz-cd}
\usepackage{tikz}
\usepackage{pict2e}
\usepackage{float}
\usepackage{stmaryrd}
\usepackage{geometry}
\usepackage{tensor}
\usepackage{mathtools}

\usepackage{amsmath,amssymb,amsthm,amscd,graphicx}
\usepackage{psfrag}

\usepackage[utf8]{inputenc}
\usepackage{graphicx}
\usepackage{xcolor}
\usepackage{hyperref}
\hypersetup{colorlinks,allcolors=blue}

\definecolor{frangreen}{rgb}{0.040, 0.475, 0.435}

\usepackage[T1]{fontenc}

\usepackage[utf8]{inputenc}
\usepackage{graphicx}
\usepackage{amsmath,amsthm,amssymb,physics,mathtools}
\usepackage{tikz,tikz-cd, dynkin-diagrams}
\usepackage{epigraph}
\usepackage{hhline}
\usepackage{ytableau}
\usepackage{natbib}
\usetikzlibrary{positioning,arrows}
\usetikzlibrary{decorations.pathmorphing}
\usetikzlibrary{decorations.markings}
\usetikzlibrary{matrix}
\usepackage{bbold}

\usepackage{environ}         
\usepackage{etoolbox}        
\usepackage{graphicx}        
\usepackage{nomencl}
\makenomenclature

 \usepackage{etoolbox}
\renewcommand\nomgroup[1]{%
  \item[\bfseries
  \ifstrequal{#1}{C}{CFT symbols}{%
  \ifstrequal{#1}{H}{Heun symbols}{%
  \ifstrequal{#1}{F}{CFT symbols - semiclassics}{}}}%
]}
	\newlength{\myl}
\let\origequation=\equation
\let\origendequation=\endequation

\RenewEnviron{equation}{
  \settowidth{\myl}{$\BODY$}                       
  \origequation
  \ifdimcomp{\the\linewidth}{>}{\the\myl}
  {\ensuremath{\BODY}}                             
  {\resizebox{\linewidth}{!}{\ensuremath{\BODY}}}  
  \origendequation
}

\theoremstyle{definition}

\newcommand{\eqlb}[2]{\begin{equation} \label{#1} #2 \end{equation}}
\newcommand{\eq}[1]{\begin{equation} #1 \end{equation}}
\newcommand{\eqn}[1]{\begin{equation*} #1 \end{equation*}}

\newcommand{\brc}[1]{\left(#1\right)}
\newcommand{\bsq}[1]{\left[#1\right]}
\newcommand{\bfi}[1]{\left\{ #1\right\}}

\newcommand{\rme}{\textrm{e}}

\newcommand{\rmd}{{\rm d}}
\newcommand{\ops}{{\rm op}}

\newcommand{\be}{\begin{equation}}
\newcommand{\ee}{\end{equation}}
\newcommand{\ba}{\begin{aligned}}
\newcommand{\ea}{\end{aligned}}
\newcommand{\ben}{\begin{eqnarray}\displaystyle}
\newcommand{\een}{\end{eqnarray}}

\makeatletter
\gdef\@fpheader{}
\makeatother

\title{Basics of Multiple Polyexponential Integrals}

\author[a,b]{Gleb Aminov}
\author[c,d,e]{Paolo Arnaudo}

\affiliation[a]{ C.~N. Yang Institute for Theoretical Physics, State University of New York, Stony Brook, NY 11794-3840, USA}
\affiliation[b]{ Simons Center for Geometry and Physics, State University of New York, Stony Brook, NY 11794-3636, USA}
\affiliation[c]{International School of Advanced Studies (SISSA), via Bonomea 265, 34136 Trieste, Italy}
\affiliation[d]{INFN Sezione di Trieste, via Valerio 2, 34127 Trieste, Italy}
\affiliation[e]{Institute for Geometry and Physics, IGAP, via Beirut 2, 34151 Trieste, Italy}

\emailAdd{gleb.aminov@stonybrook.edu}
\emailAdd{parnaudo@sissa.it}
		
\date{\today}

\abstract{We introduce a set of special functions called \emph{multiple polyexponential integrals}, defined as iterated integrals of the exponential integral $\text{Ei}(z)$. 
These functions arise in certain perturbative expansions of the local solutions of second-order ODEs around an irregular singularity. In particular, their recursive definition describes the asymptotic behavior of these local solutions.
To complement the study of the multiple polyexponential integrals on the entire complex plane, we relate them with two other sets of special functions -- the \emph{undressed} and \emph{dressed multiple polyexponential functions} -- which are characterized by their Taylor series expansions around the origin. }

\begin{document}

\maketitle

\section{Introduction}

The physical motivation for studying the \emph{multiple polyexponential integrals} and \emph{multiple polyexponential functions} arises in the context of the linear perturbations of asymptotically flat black holes. After the separation of variables, the radial part of these perturbations is described by the second-order ordinary differential equation (ODE) with two regular and one irregular singularities. In the small-frequency expansion, the local solutions in the near-horizon region can be written using multiple polylogarithm functions. At the same time, in the near-spatial infinity region, one needs to introduce generalizations of the exponential integral - \emph{multiple polyexponential integrals}.

Among the gravitational quantities that can be computed via numerical or analytical methods, a special role is played by the quasinormal mode frequencies (QNMs). These realize a set of discrete and complex frequencies that are responsible for the damped oscillations appearing in the last phase of the merger of two black holes (see \cite{Berti:2009kk} for a review). In particular, since the recent experimental verification of gravitational waves \cite{PhysRevLett.116.061102}, there has been a growing interest in finding new methods to compute with high precision these quantities and compare the results with the observations \cite{PhysRevD.30.295,PhysRevD.35.3632,Mano:1996mf,Mano:1996vt,Mano:1996gn,Cardoso:2003vt,Konoplya:2004ip,novaes, novaes2014,CarneirodaCunha:2015hzd,Novaes:2018fry,   CarneirodaCunha:2019tia,Amado:2020zsr,Aminov:2020yma,Hatsuda:2020sbn,Hatsuda:2020egs,BarraganAmado:2021uyw,Bonelli:2021uvf,Bianchi:2021xpr,Bianchi:2021mft,Amado:2021erf,Hatsuda:2021gtn,Fioravanti:2021dce,Bonelli:2022ten,Dodelson:2022yvn,Consoli:2022eey,Imaizumi:2022qbi,Ivanov:2022qqt,daCunha:2022ewy,Imaizumi:2022dgj,Bianchi:2022qph, Gregori:2022xks,Dodelson:2023vrw,Bianchi:2023rlt,Bianchi:2023sfs,Giusto:2023awo,Hatsuda:2023geo, Fioravanti:2023zgi,Saketh:2023bul,Lei:2023mqx,Bautista:2023sdf,Ivanov:2024sds,Arnaudo:2024rhv}.

In the problems we are interested in, the perturbation field satisfies a second-order ODE and is constrained to two suitable boundary conditions. The quantization of the QNMs can be obtained by solving the differential equation locally around the boundary points and gluing the resulting solutions. In the case of asymptotically flat black holes, one of the boundary conditions is imposed at spatial infinity, which is an irregular singularity of rank one of the differential equation.
When working in the small-frequency regime around an irregular singularity \cite{Aminov:2024mul}, we encounter leading order solutions of the ODE which are given by products of rational functions and exponential functions. To compute the higher-order corrections, we apply the perturbative method described in \cite{Aminov:2023jve}. 
The first correction involves the exponential integral
\begin{equation}
\mathrm{Ei}(z)=\int_{-\infty}^z\frac{\mathrm{e}^t}{t}\mathrm{d}t.
\end{equation}
This function can also be written as a Taylor series expansion around $z=0$ plus a logarithm:
\begin{equation}\label{Eiandseries}
\mathrm{Ei}(z)=\gamma+\log(-z)+\sum_{k=1}^{\infty}\frac{z^k}{k!\,k},\quad |\mathrm{Arg}(-z)|<\pi,
\end{equation}
where $\gamma$ is the Euler-Mascheroni constant. 
When computing the next orders of perturbation, we need to take further integrations, and an iterative structure similar to the one defining the multiple polylogarithms appears (for the latter, see \cite{goncharov2001multiple,Wald}). Therefore, we will define functions that generalize the exponential integral as multiple polylogarithms do for the logarithm. We are also interested in the properties of the multiple polyexponential integrals for both $z\to \infty$ and $z\to 0$.

Taking inspiration from the recursive integral structure of the multiple polylogarithms, we define the set of \emph{multiple polyexponential integrals} as
\begin{equation}
\begin{aligned}
\text{ELi}_{1,s_2,\dots,s_n}\brc{z} &= - \int_{-\infty}^{z} \frac{\rme^{t}}{t}\, \text{ELi}_{s_2,\dots,s_n}\brc{-t} \rmd t,\\
s_1>1:\quad
\text{ELi}_{s_1,\dots,s_n}\brc{z}&=\int_{-\infty}^{z} \frac{1}{t}\,\text{ELi}_{s_1-1,s_2,\dots,s_n}\brc{t} \rmd t,
\end{aligned}
\end{equation}
where $n, s_1,\dots,s_n\in\mathbb{Z}_{>0}$ and the starting point is
\begin{equation}
\mathrm{ELi}_1(z)=\mathrm{Ei}(z).
\end{equation}
The integer $n$ is called \emph{level} and the sum $s_1+\dots +s_n$ is called \emph{weight}. 
These functions have an explicit asymptotic behavior when $z\to\infty$ and satisfy the following simple recursive derivative relations:
\begin{equation}\label{intro:recurrenceELi}
z\frac{\mathrm{d}}{\mathrm{d}z}\mathrm{ELi}_{s_1,s_2,\dots,s_n}(z)=\begin{cases}
-\mathrm{e}^z\,\mathrm{ELi}_{s_2,\dots,s_n}(-z)\quad &s_1=1,\\
\mathrm{ELi}_{s_1-1,s_2,\dots,s_n}(z)\quad &s_1>1.
\end{cases}
\end{equation}

To make contact with the behavior at $z=0$, we will generalize the Taylor series in \eqref{Eiandseries}
and define the set of \emph{undressed multiple polyexponential functions}:
\begin{equation}
el_{s_1,s_2\dots,s_n}(z)=\sum_{k_1>k_2>\dots>k_n\ge 1}\frac{1}{k_1^{s_1}k_2^{s_2}\dots k_n^{s_n}}\frac{z^{k_1}}{k_1!},
\end{equation}
where again $n, s_1,\dots,s_n\in\mathbb{Z}_{>0}$. Functions with similar Taylor expansions have been discussed in \cite{boyadzhiev2007polyexponentials,kim2020degenerate,Kim2019ANO,KIM2020124017,komatsu1,komatsu2,lacpao2019hurwitz}, but their asymptotic behavior at $z\to\infty$ was not studied systematically.
As follows from the definition, undressed multiple polyexponential functions do not satisfy simple derivative relations like (\ref{intro:recurrenceELi}). Instead, their derivatives involve sums over all ordered partitions of their indices.
Thus, we also define \emph{dressed multiple polyexponential functions} $\mathrm{EL}_{s_1,\dots,s_n}(z)$ that satisfy derivative relations analogous to \eqref{intro:recurrenceELi}:
\begin{equation}\label{recurrenceEL}
z\frac{\mathrm{d}}{\mathrm{d}z}\mathrm{EL}_{s_1,s_2,\dots,s_n}(z)=\begin{cases}
-\mathrm{e}^z\,\mathrm{EL}_{s_2,\dots,s_n}(-z)\quad &s_1=1,\\
\mathrm{EL}_{s_1-1,s_2,\dots,s_n}(z)\quad &s_1>1.
\end{cases}
\end{equation}
The corresponding Taylor expansions around $z=0$ are obtained using the relations between dressed and undressed multiple polyexponential functions.
Finally, the behavior of the multiple polyexponential integrals at $z=0$ is given by the relations similar to \eqref{Eiandseries}, which reads in the new notations as
\begin{equation}
\mathrm{ELi}_1(z)=\gamma+\log(-z)+\mathrm{EL}_1(z).
\end{equation}
In general, these relations involve powers of the logarithm function and different constants such as $\gamma$ and multiple zeta values (MZVs). For useful properties of the MZVs, see \cite{Borwein1996,Blumlein:2009cf,em/1062621000,Flajolet1998EulerSA,xu2020explicit,XU2017443,CHEN2015107,kuba2019multisets,hoffman2021logarithmic,hoffman2015quasi,seki2020ohno,Si+2021+1612+1619,kuba2019note,batir2017some,muneta2007some}.

This paper is structured as follows.
In Sec.~\ref{sec:mpf_def}, we introduce the undressed and dressed multiple polyexponential functions. In Sec.~\ref{sec:taylor}, we study the relations between the said functions, which allow us to derive the corresponding Taylor series expansions around $z=0$. In Sec.~\ref{sec:polyexp}, we define the polyexponential integrals and introduce the relations with the corresponding polyexponential functions. We generalize these relations to multiple polyexponential integrals in Sec.~\ref{sec:relationsELandELi}. The behavior of the multiple polyexponential integrals in the asymptotic region $z\to\infty$ is studied in Sec.~\ref{sec:asymptotics}.
Finally, in Sec.~\ref{sec:quadid}, we look into quadratic identities satisfied by the dressed multiple polyexponential functions. These identities are inspired by the ones for the multiple polylogarithms, such as the shuffle and stuffle relations \cite{Wald, Lewin'81, Goncharov:1998kja, MINH2000217,MINH2000273}.
We postpone to the appendices the proofs of more involved relations and identities.
Appendix \ref{sec:App_diff} proves the derivative relation satisfied by the undressed multiple polyexponential functions. Appendix \ref{app:rel_Taylor} proves the general relation between undressed and dressed multiple polyexponential functions. Appendix \ref{app:rel23_proof} proves the relations between multiple polyexponential integrals of levels 2 and 3 with the corresponding dressed multiple polyexponential functions. Appendix \ref{app:reln_proof} proves the general case of level $n\geq 2$.

\vskip 0.5cm
\noindent {\large {\bf Acknowledgments}}
\vskip 0.5cm
We thank G. Bonelli, A. Grassi, A. Grekov, and A. Tanzini for their valuable remarks and discussions.   

The research of P.A.\ is partly supported by the INFN Iniziativa Specifica GAST and
by the MIUR PRIN Grant 2020KR4KN2 
String Theory as a bridge between Gauge Theories and Quantum Gravity.
P.A. acknowledges funding from the EU project  
Caligola HORIZON-MSCA-2021-SE-01), Project ID: 101086123.

\section{Multiple polyexponential functions: definitions}\label{sec:mpf_def}

First, we define a set of functions we call undressed multiple polyexponential functions:
\eqlb{eq:polyExp_def}{el_n\brc{z}=\sum_{k=1}^{\infty} \frac{z^k}{k^n\, k!},\quad
el_{s_1,\dots,s_n}\brc{z}=\sum_{k_1>k_2>\dots>k_n\geq 1} 
\frac{1}{k_1^{s_1} \dots k_n^{s_n}}\, \frac{z^{k_1}}{ k_1 !}.}
The latter functions are straightforward in terms of their series expansion around $z=0$, but taking the derivative proves more challenging. If the first index $s_1$ is greater than one, we have a simple polylogarithm-like derivative
\eqlb{eq:el2_diff}{s_1>1:\quad z\,\frac{\rmd }{\rmd z} el_{s_1,\dots,s_n}\brc{z}= el_{s_1-1,s_2,\dots,s_n}\brc{z}.}
However, when $s_1=1$, the derivative rule becomes harder:
\begin{equation}
\label{eq:el1_diff}
z\,\frac{\rmd }{\rmd z} el_{1,s_2,\dots,s_n}\brc{z}=-el_{s_2,\dots,s_n}\brc{z} - \brc{-1}^n \rme^{z}
\sum_{\ops\brc{s_2}}\dots \sum_{\ops\brc{s_n}} el_{\ops\brc{s_2},\dots,\ops\brc{s_n}} \brc{-z},
\end{equation}
where the sum is over all ordered partitions of  $s_i\in \mathbb{N}$, $i\geq 2$:
\eq{\ops\brc{1}=\bfi{1},\quad \ops\brc{2}=\bfi{2,\brc{1,1}},\quad
\ops\brc{3}=\bfi{3,\brc{2,1},\brc{1,2},\brc{1,1,1}}}
and so on. Some particular cases of (\ref{eq:el1_diff}) are
\eqlb{eq:diff_ex11}{z\,\frac{\rmd}{\rmd z} el_{1,1}\brc{z} = - el_1\brc{z} - \rme^{z} el_1\brc{-z},}
\eq{z\,\frac{\rmd}{\rmd z} el_{1,2}\brc{z} = - el_2\brc{z} - \rme^{z} el_2\brc{-z}- \rme^{z} el_{1,1}\brc{-z},}
\eq{z\,\frac{\rmd}{\rmd z} el_{\bfi{1}_n}\brc{z} = - el_{\bfi{1}_{n-1}}\brc{z} - \brc{-1}^n \rme^{z} el_{\bfi{1}_{n-1}}\brc{-z}.}
We prove the second derivative rule (\ref{eq:el1_diff}) in Appendix \ref{sec:App_diff}.
Similarly to multiple polyexponential functions, we can introduce the notions of level and weight. For a given function $el_{s_1,\dots,s_n}\brc{z}$, the level is $n$, and the weight $w_n$ is the sum of all indices:
\eq{w_n=\sum_{i=1}^{n} s_i.}
Although undressed multiple polyexponential functions are linearly independent, there are identities that involve non-linear terms. Organizing these identities according to the weight, we have for weights $2$ and $3$:
\begin{equation}
\label{eq:el_id1}
el_{1,1}\brc{z}+el_{1,1}\brc{-z}+ 2\, el_2\brc{z}+2\, el_2\brc{-z} +el_1\brc{z}el_1\brc{-z}=0,
\end{equation}
\begin{equation}
\begin{aligned}
&el_{1,1,1}\brc{z}+el_{1,1,1}\brc{-z}+el_{1,2}\brc{z}+el_{1,2}\brc{-z}+
3\,el_{2,1}\brc{z}+3\,el_{2,1}\brc{-z}\\
&+6\, el_3\brc{z}+6\, el_3\brc{-z}
+el_1\brc{z}el_2\brc{-z}+el_1\brc{-z}el_2\brc{z}=0.
\end{aligned}
\end{equation}
Later, we will prove more general forms of such identities.

The next important observation is
\eqlb{eq:obs_el1n}{z\,\frac{\rmd}{\rmd z} \sum_{\ops\brc{n}}el_{1,\ops\brc{n}}\brc{z} =
- \sum_{\ops\brc{n}}el_{\ops\brc{n}}\brc{z} -\rme^z \,el_n\brc{-z},}
which is a particular case of a more general identity given in Appendix \ref{app:rel_Taylor}.
This leads us to a definition of a first dressed multiple polyexponential function:
\eqlb{eq:EL1n_def}{\text{EL}_{1,n}\brc{z}\equiv\sum_{\ops\brc{n+1}}el_{\ops\brc{n+1}}\brc{z}}
with the following derivative rule:
\eqlb{eq:EL1n_der}{z\,\frac{\rmd}{\rmd z}\, \text{EL}_{1,n}\brc{z} = - \rme^z \,\text{EL}_{n}\brc{-z},}
where
\eq{\text{EL}_{n}\brc{z}\equiv el_{n}\brc{z}.}

To prove (\ref{eq:EL1n_der}), we will use the fact that the sum over all ordered partitions of $n+1$ can be written in terms of the sum over all ordered partitions of $n$ in the following way:
\eqlb{eq:op_dec}{\sum_{\ops\brc{n+1}}el_{\ops\brc{n+1}}\brc{z}=\sum_{\ops\brc{n}}
\bsq{el_{1,\ops\brc{n}}\brc{z}+ el_{1\oplus\ops\brc{n}}\brc{z}},}
where the operator $\oplus$ takes the number to the left of it and adds it to the first entry of the vector to the right of it, e.g., $1\oplus(1,1)=(2,1)$. In a more general case, when this operator is applied between two vectors $\textbf{v}=\brc{v_1,\dots,v_i}$ and $\textbf{u}=\brc{u_1,\dots,u_j}$, we define
\eq{\textbf{v} \oplus \textbf{u}=\brc{v_1,\dots, v_{i-1},v_i+u_1,u_2,\dots,u_j}.}
If the second vector is a scalar, we have:
\eq{\textbf{v} \oplus u= \brc{v_1,\dots, v_{i-1},v_i+u}.}
Using (\ref{eq:EL1n_def}) and (\ref{eq:op_dec}), we get
\eq{z\,\frac{\rmd}{\rmd z}\, \text{EL}_{1,n}\brc{z} = - \sum_{\ops\brc{n}}el_{\ops\brc{n}}\brc{z}-\rme^z \,el_n\brc{-z} +
z\,\frac{\rmd}{\rmd z} \sum_{\ops\brc{n}}el_{1\oplus\ops\brc{n}}\brc{z}.}
Since the following is true for any ordered partition of $n$:
\eq{z\,\frac{\rmd}{\rmd z}\,el_{1\oplus\ops\brc{n}}\brc{z}=el_{\ops\brc{n}}\brc{z},}
we arrive at (\ref{eq:EL1n_der}).

Motivated by the definition of $\text{EL}_{1,n}\brc{z}$, we define the complete set of dressed multiple polyexponential functions via the following recursive derivative rules:
\eq{z\,\frac{\rmd }{\rmd z}\, \text{EL}_{1,s_2,\dots,s_n}\brc{z}= - \rme^z \, \text{EL}_{s_2,\dots,s_n}\brc{-z},}
\eqlb{eq:EL_gen_der}{s_1>1:\quad z\,\frac{\rmd }{\rmd z} \, \text{EL}_{s_1,\dots,s_n}\brc{z}= \text{EL}_{s_1-1,s_2,\dots,s_n}\brc{z},}
where the integration constants are fixed by
\eq{\text{EL}_{s_1,\dots,s_n}\brc{0}=0.}
Compared with undressed functions $el_{s_1,\dots,s_n}\brc{z}$, the new set of functions has simple derivative rules. However, their series expansions around $z=0$ still need to be derived. In our first example (\ref{eq:EL1n_def}), the Taylor expansion is given by
\eq{\text{EL}_{1,n}\brc{z}=\sum_{k_1\geq k_2\geq\dots \geq k_{n+1}\geq 1}
\frac{1}{k_1 k_2 \dots k_{n+1}}\, \frac{z^{k_1}}{ k_1 !},}
where we used the following simple fact:
\eqlb{eq:op_Taylor}{\sum_{k_1\geq k_2\geq\dots \geq k_{n+1}\geq 1} \frac{1}{k_1 k_2 \dots k_{n+1}}=
\sum_{\textbf{v}\in \ops\brc{n+1}} \sum_{k_1>\dots > k_{j}\geq 1} \frac{1}{k_1^{v_1} \dots k_{j}^{v_j}}.}
In the latter equation, vector $\textbf{v}=\brc{v_1,\dots,v_j}$ represents elements of the set of all ordered partitions of $n+1$ (including the scalar contribution $n+1$).

\section{Multiple polyexponential functions: Taylor expansions}\label{sec:taylor}

In this section, we will derive the relations between the dressed and undressed multiple polyexponential functions, which will lead to Taylor expansions of the dressed functions $\text{EL}_{s_1,\dots,s_n}\brc{z}$. The first example of such a relation was given previously in (\ref{eq:EL1n_def}). We derive the general relation for level $2$ dressed functions $\text{EL}_{m,n}\brc{z}$ by integrating the derivative rule (\ref{eq:EL_gen_der}) recursively:
\eq{\text{EL}_{m,n}\brc{z}=\sum_{\ops\brc{n+1}}el_{\brc{m-1}\oplus\ops\brc{n+1}}\brc{z}.}
The corresponding Taylor expansion follows from (\ref{eq:op_Taylor}) and the derivative rule (\ref{eq:el2_diff}):
\eqlb{eq:Taylor_ELmn}{\text{EL}_{m,n}\brc{z}=\sum_{k_1\geq k_2\geq\dots \geq k_{n+1}\geq 1}
\frac{1}{k_1^m k_2 \dots k_{n+1}}\, \frac{z^{k_1}}{ k_1 !}.}

To figure out the general relation at level $3$, we start with the function $\text{EL}_{1,1,n}\brc{z}$ and its derivative:
\eq{z\,\frac{\rmd }{\rmd z}\, \text{EL}_{1,1,n}\brc{z}= - \rme^z \, \text{EL}_{1,n}\brc{-z},}
where the right-hand side can be rewritten using the definition (\ref{eq:EL1n_def}):
\eq{z\,\frac{\rmd }{\rmd z}\, \text{EL}_{1,1,n}\brc{z}= - \rme^z \sum_{\ops\brc{n+1}}el_{\ops\brc{n+1}}\brc{-z}.}
Now, we can use the derivative rule (\ref{eq:el1_diff}) for the undressed function
\eq{z\,\frac{\rmd }{\rmd z}\, el_{1,n+1}\brc{z}=-el_{n+1}\brc{z}-\rme^z \sum_{\ops\brc{n+1}}el_{\ops\brc{n+1}}\brc{-z}}
to derive
\eq{z\,\frac{\rmd }{\rmd z}\, \text{EL}_{1,1,n}\brc{z}=z\,\frac{\rmd }{\rmd z}\bsq{el_{n+2}\brc{z} +el_{1,n+1}\brc{z}}.}
Since $\text{EL}_{1,1,n}\brc{0}=0$, one gets
\eqlb{eq:EL_rel11n}{\text{EL}_{1,1,n}\brc{z} = el_{n+2}\brc{z} +el_{1,n+1}\brc{z}.}
Repeating this procedure for increasing values of the second index, we arrive at the following relation:
\eqlb{eq:EL_rel1mn}{\text{EL}_{1,m,n}\brc{z} =\sum_{\ops\brc{m+1}}el_{\ops\brc{m+1}\oplus n}\brc{z},}
which can be proved by induction. The base case of this induction is (\ref{eq:EL_rel11n}). Assuming that (\ref{eq:EL_rel1mn}) is true for some fixed $m=k-1$, we take the derivative of the sum in (\ref{eq:EL_rel1mn}) with $m=k\geq 2$:
\eq{z\,\frac{\rmd }{\rmd z}\sum_{\ops\brc{k+1}}el_{\ops\brc{k+1}\oplus n}\brc{z}=z\,\frac{\rmd }{\rmd z}\sum_{\ops\brc{k}}
\bsq{el_{\ops\brc{k},\brc{n+1}}\brc{z}+el_{\ops\brc{k}\oplus \brc{n+1}}\brc{z}},}
where we used the decomposition similar to (\ref{eq:op_dec}):
\eqlb{eq:op_dec2}{\sum_{\ops\brc{n+1}}el_{\ops\brc{n+1}}\brc{z}=\sum_{\ops\brc{n}}
\bsq{el_{\ops\brc{n},1}\brc{z}+ el_{\ops\brc{n}\oplus 1}\brc{z}}.}
Since we assume (\ref{eq:EL_rel1mn}) is valid for $m=k-1$ and any $n$, we have:
\eqlb{derivativegeneratesexpEL}{z\,\frac{\rmd }{\rmd z}\sum_{\ops\brc{k}} el_{\ops\brc{k}\oplus \brc{n+1}}\brc{z}= - \rme^z \, \text{EL}_{k-1,n+1}\brc{-z} =
- \rme^z \sum_{\ops\brc{n+2}} el_{\brc{k-2}\oplus \ops\brc{n+2}}\brc{-z}.}
For the remaining sum we use decomposition (\ref{eq:op_dec}) to get
\eq{z\,\frac{\rmd }{\rmd z}\sum_{\ops\brc{k}}el_{\ops\brc{k},\brc{n+1}}\brc{z} = \sum_{\ops\brc{k-1}}el_{\ops\brc{k-1},\brc{n+1}}\brc{z} +
z\,\frac{\rmd }{\rmd z}\sum_{\ops\brc{k-1}}el_{1,\ops\brc{k-1},\brc{n+1}}\brc{z}.}
Analogously to (\ref{eq:obs_el1n}), the following relation holds (see Appendix \ref{app:rel_Taylor} for the proof):
\eqlb{eq:threeblockderivativeel}{z\,\frac{\rmd }{\rmd z}\sum_{\ops\brc{k-1}}el_{1,\ops\brc{k-1},\brc{n+1}}\brc{z} =-\sum_{\ops\brc{k-1}}el_{\ops\brc{k-1},\brc{n+1}}\brc{z} +
\rme^z \sum_{\ops\brc{n+1}} el_{\brc{k-1}, \ops\brc{n+1}}\brc{-z}.}
Summing up all the contributions leads to
\begin{equation}
\begin{aligned}
&z\,\frac{\rmd }{\rmd z}\sum_{\ops\brc{k+1}}el_{\ops\brc{k+1}\oplus n}\brc{z}=- \rme^z \sum_{\ops\brc{n+2}} el_{\brc{k-2}\oplus \ops\brc{n+2}}\brc{-z} + \rme^z \sum_{\ops\brc{n+1}} el_{\brc{k-1}, \ops\brc{n+1}}\brc{-z}=\\
&=- \rme^z \sum_{\ops\brc{n+1}} el_{\brc{k-1}\oplus \ops\brc{n+1}}\brc{-z}=- \rme^z \, \text{EL}_{k,n}\brc{-z},
\end{aligned}
\end{equation}
which finishes the induction step and proves (\ref{eq:EL_rel1mn}). Finally, the general relation at level $3$ is
\eq{\text{EL}_{l,m,n}\brc{z} =\sum_{\ops\brc{m+1}}el_{\brc{l-1}\oplus\ops\brc{m+1}\oplus n}\brc{z},}
and the corresponding Taylor expansion is
\eq{\text{EL}_{l,m,n}\brc{z} = \sum_{k_1\geq k_2\geq\dots \geq k_{m+1}\geq 1}
\frac{1}{k_1^l k_2 \dots k_m k_{m+1}^{n+1}}\, \frac{z^{k_1}}{ k_1 !}.}

Repeating the same steps at level $4$, we get
\eqlb{eq:EL_reljlmn}{\text{EL}_{j,l,m,n}\brc{z} = \sum_{\ops\brc{l+1}}\sum_{\ops\brc{n+1}} el_{\brc{j-1}\oplus\ops\brc{l+1}\oplus\brc{m-1}\oplus \ops\brc{n+1}}\brc{z}.}
To construct the corresponding Taylor expansion, we split the indices of the undressed functions into two parts: $\brc{j-1}\oplus\ops\brc{l+1}$ and $\brc{m-1}\oplus \ops\brc{n+1}$. Then, we consider two vectors $\textbf{v}\in \brc{j-1}\oplus\ops\brc{l+1}$ and $\textbf{u}\in \brc{m-1}\oplus \ops\brc{n+1}$. Assuming $\textbf{v}=\brc{v_1,\dots, v_i}$ and $\textbf{u}=\brc{u_1,\dots, u_j}$, the corresponding contribution to the sum on the right-hand side of (\ref{eq:EL_reljlmn}) is:
\eq{el_{\textbf{v}\oplus \textbf{u}}\brc{z}=\sum_{k_1>\dots >k_{i+j-1}\geq 1}\frac{1}{k_1^{v_1} \dots k_{i-1}^{v_{i-1}} k_i^{v_i+u_1} k_{i+1}^{u_2} \dots k_{i+j-1}^{u_j}}\, \frac{z^{k_1}}{ k_1 !}.}
Summing over all possible vectors $\textbf{v}$, $\textbf{u}$ and using the previous result (\ref{eq:Taylor_ELmn}), we arrive at
\eq{\text{EL}_{j,l,m,n}\brc{z} = \sum_{k_1\geq\dots \geq k_{l+n+1}\geq 1}
\frac{1}{k_1^j \, k_2 \dots k_l\, k_{l+1}^{m+1} k_{l+2}\dots k_{l+n+1}}\, \frac{z^{k_1}}{ k_1 !}.}

The generalization to an arbitrary number of indices is straightforward now (see Appendix \ref{app:rel_Taylor} for proof). For an even number of indices $n=2\,k$, the relation between the dressed and undressed functions is
\eqlb{eq:EL_rel_even}{\boxed{\text{EL}_{r_1,s_1,\dots,r_k,s_k}\brc{z} = \sum_{\ops\brc{s_1+1}}\dots\sum_{\ops\brc{s_k+1}} el_{\brc{r_1-1}\oplus\ops\brc{s_1+1}\oplus\dots\oplus\brc{r_k-1}\oplus \ops\brc{s_k+1}}\brc{z}.}}
If the number of indices $n$ is odd, $n=2\,k+1$ with $k\geq 1$, we have
\eqlb{eq:EL_rel_odd}{\boxed{\text{EL}_{r_1,s_1,\dots,r_k,s_k,r_{k+1}}\brc{z} = \sum_{\ops\brc{s_1+1}}\dots\sum_{\ops\brc{s_k+1}} el_{\brc{r_1-1}\oplus\ops\brc{s_1+1}\oplus\dots\oplus\brc{r_k-1}\oplus \ops\brc{s_k+1}\oplus r_{k+1}}\brc{z}.}}
We introduce the following notation
\eq{1\leq j \leq k:\quad w_j=\sum_{i=1}^{j} s_i}
to write down the corresponding Taylor expansions:
\eqlb{taylor_el}{\text{EL}_{r_1,s_1,\dots,r_k,s_k}\brc{z} = \sum_{k_1\geq\dots \geq k_{w_k+1}\geq 1}
\brc{\prod_{j=2}^{w_k+1}\frac{1}{k_j}} \brc{\prod_{j=1}^{k-1}\frac{1}{k_{w_j+1}^{r_{j+1}}}}
\frac{z^{k_1}}{k_1^{r_1} k_1 !},}
\eq{\text{EL}_{r_1,s_1,\dots,r_k,s_k,r_{k+1}}\brc{z} = \sum_{k_1\geq\dots \geq k_{w_k+1}\geq 1}
\brc{\prod_{j=2}^{w_k+1}\frac{1}{k_j}} \brc{\prod_{j=1}^{k}\frac{1}{k_{w_j+1}^{r_{j+1}}}}
\frac{z^{k_1}}{k_1^{r_1} k_1 !}.}

\section{Polyexponential integrals}\label{sec:polyexp}

The asymptotic behavior of polyexponential functions $\text{EL}_n\brc{z}$ at $z \rightarrow \pm\infty$ is determined by their relations with a set of functions called polyexponential integrals $\text{ELi}_{n}\brc{z}$, defined as
\eqlb{eq:polyexpI_def}{n\geq2:\quad \text{ELi}_{n}\brc{z} = \int_{-\infty}^{z} \frac{\text{ELi}_{n-1}\brc{t}}{t} \rmd t,\quad
\text{ELi}_{1}\brc{z}\equiv \text{Ei}\brc{z}.}
Using L'H\^{o}pital's rule, the leading asymptotic behavior of $\text{ELi}_{n}\brc{z}$ at $z \rightarrow \pm \infty$ can be derived:
\eq{\text{ELi}_{n}\brc{z}= \frac{\rme^{z}}{z^n}\brc{1+O\brc{\abs{z}^{-1}}}.}
For the full asymptotic series of the polyexponential integrals see Sec.~\ref{sec:asymptotics}.

Starting with the known relation for the exponential integral
\eqlb{eq:ELi1}{\text{Ei}\brc{z}=\gamma +\log\brc{-z}+\text{EL}_{1}\brc{z}}
and applying recursively the definition (\ref{eq:polyexpI_def}), we get the following general result for $n\geq 1$:
\eqlb{eq:ELn_gen}{\boxed{
\text{ELi}_{n}\brc{z}=\text{EL}_{n}\brc{z}+\sum_{k=0}^{n} \frac{\brc{-1}^{n-k}}{k! \brc{n-k}!}\, \Gamma^{\brc{n-k}}\brc{1} \log\brc{-z}^k .}}
The latter can be proven by induction with a base case given by (\ref{eq:ELi1}). Assuming that the claim is true for some $n>1$, we use the definition to get for $n+1$:
\begin{equation}
\text{ELi}_{n+1}\brc{z} = \int_{-\infty}^{z} \frac{\text{ELi}_{n}\brc{t}}{t} \rmd t,
\end{equation}
where the integral can be split into two parts:
\begin{itemize}
\item[1.] $\displaystyle \int_{-\infty}^{-1} \frac{\text{ELi}_{n}\brc{t}}{t} \rmd t$,
\item[2.] $\displaystyle \int_{-1}^{z} \frac{\text{ELi}_{n}\brc{t}}{t} \rmd t=\int_{-1}^{0} \frac{\text{EL}_{n}\brc{t}}{t} \rmd t+\int_{0}^{z} \frac{\text{EL}_{n}\brc{t}}{t} \rmd t+\sum_{k=0}^n\frac{(-1)^{n-k}}{k!\,(n-k)!}\Gamma^{(n-k)}(1)\int_{-1}^z\frac{\log(-t)^k}{t}\rmd t$.
\end{itemize}  
The last two integrals in the second part can be taken straightforwardly:
\begin{equation}
\begin{aligned}
&\int_{0}^{z} \frac{\text{EL}_{n}\brc{t}}{t} \rmd t=\text{EL}_{n+1}\brc{z},\\
&\sum_{k=0}^n\frac{(-1)^{n-k}}{k!\,(n-k)!}\Gamma^{(n-k)}(1)\int_{-1}^z\frac{\log(-t)^k}{t}\rmd t=\sum_{k=1}^{n+1}\frac{(-1)^{n+1-k}}{k!\,(n+1-k)!}\Gamma^{(n+1-k)}(1)\log(-z)^{k}.
\end{aligned}
\end{equation}
For the remaining integrals, we repeatedly integrate by parts until the weights of the polyexponential integral and the polyexponential function are reduced to 0:
\begin{equation}\label{eq:poly_int_parts}
\begin{aligned}
&\int_{-\infty}^{-1} \frac{\text{ELi}_{n}\brc{t}}{t} \rmd t=-\int_1^{\infty}\frac{\text{ELi}_{n}\brc{-t}}{t} \rmd t=-\int_1^{\infty}\text{ELi}_{n}\brc{-t} \rmd \log(t)=\dots=\frac{(-1)^{n+1}}{(n+1)!}\int_1^{\infty}\log^{n+1}(t)\,\mathrm{e}^{-t}\, \rmd t,\\
&\int_{-1}^{0} \frac{\text{EL}_{n}\brc{t}}{t} \rmd t=-\int_0^{1}\frac{\text{EL}_{n}\brc{-t}}{t} \rmd t=-\int_0^{1}\text{EL}_{n}\brc{-t} \rmd \log(t)=\dots=\frac{(-1)^{n+1}}{(n+1)!}\int_0^{1}\log^{n+1}(t)\,\mathrm{e}^{-t}\, \rmd t,
\end{aligned}
\end{equation}
where the following properties were used:
\begin{equation}
\text{EL}_{n}\brc{0}=0,\quad\lim_{z\to-\infty}\text{ELi}_{n}\brc{z}=0,\quad \forall\,n\in\mathbb{N}.
\end{equation}
Combining the two results in (\ref{eq:poly_int_parts}) gives
\eq{\frac{(-1)^{n+1}}{(n+1)!}\int_0^{\infty}\log^{n+1}(t)\,\mathrm{e}^{-t}\, \rmd t=\frac{(-1)^{n+1}}{(n+1)!}\Gamma^{(n+1)}(1).}
Putting everything together, we end up with
\begin{equation}
\text{ELi}_{n+1}\brc{z}=\text{EL}_{n+1}\brc{z}+\sum_{k=0}^{n+1}\frac{(-1)^{n+1-k}}{k!\,(n+1-k)!}\Gamma^{(n+1-k)}(1)\log(-z)^{k},
\end{equation}
which concludes the inductive proof.

\section{Multiple polyexponential integrals}\label{sec:relationsELandELi}

The previously defined dressed multiple polyexponential functions have the following property:
\eq{\text{EL}_{s_1,\dots,s_n}\brc{0}=0.}
Thus, they are naturally described in terms of their Taylor expansions.
The corresponding multiple polyexponential integrals tend to zero as $z\rightarrow - \infty$ and can be defined using definite integrals:
\eqlb{eq:ELi_def1}{\text{ELi}_{1,s_2,\dots,s_n}\brc{z} = - \int_{-\infty}^{z} \frac{\rme^{t}}{t}\, \text{ELi}_{s_2,\dots,s_n}\brc{-t} \rmd t,}
\eqlb{eq:ELi_def}{s_1>1:\quad
\text{ELi}_{s_1,\dots,s_n}\brc{z}=\int_{-\infty}^{z} \frac{1}{t}\,\text{ELi}_{s_1-1,s_2,\dots,s_n}\brc{t} \rmd t.}
These integrals have natural asymptotic expansions around infinity, which do not include logarithm function $\log \brc{z}$. We postpone the discussion of the asymptotic expansions to the next section and focus on the relations between multiple polyexponential functions and multiple polyexponential integrals. In what follows, we use the fact that all multiple polyexponential integrals $\text{ELi}_{s_1,\dots,s_n}\brc{z}$ behave at most like $1/z$ when $z\rightarrow -\infty$ (see Sec.~\ref{sec:asymptotics}).

\subsection{Level \texorpdfstring{$2$}{} integrals}
Level $2$ multiple polyexponential integrals are the functions of the general form $\text{ELi}_{m,n}\brc{z}$.
The first example is:
\eqlb{eq:def_ELi11}{\text{ELi}_{1,1}\brc{z}=-\int_{-\infty}^{z}\frac{\rme^{t}}{t}\, \text{ELi}_{1}\brc{-t} \rmd t.}
One can derive the expansion around $z=0$ using relations (\ref{eq:ELn_gen}):
\eq{\text{ELi}_{1,1}\brc{z}=-\int_{-\infty}^{z}\frac{\rme^{t}}{t}\bsq{\gamma+\log t + \text{EL}_1\brc{-t}} \rmd t.}
First, we can take the following integrals:
\eq{\int_{-\infty}^{z}\gamma\, \frac{\rme^{t}}{t} \,\rmd t=\gamma\, \text{ELi}_{1}\brc{z},}
\eq{\int_{-\infty}^{z} \frac{\rme^{t}}{t} \,\log t\,\rmd t = \int_{-\infty}^{z}\log t\, \rmd\, \text{ELi}_{1}\brc{t}=
\text{ELi}_{1}\brc{z} \log z - \text{ELi}_{2}\brc{z}.}
The remaining integral can be split into two parts:
\eq{\int_{-\infty}^{z} \frac{\rme^{t}}{t} \, \text{EL}_{1}\brc{-t}\rmd t = \int_{-\infty}^{0} \frac{\rme^{t}}{t} \, \text{EL}_{1}\brc{-t}\rmd t+\int_{0}^{z} \frac{\rme^{t}}{t} \, \text{EL}_{1}\brc{-t}\rmd t.}
The first part can be computed using the series expansion of $\text{EL}_{1}$ function:
\eq{\int_{-\infty}^{0} \frac{\rme^{t}}{t} \, \text{EL}_{1}\brc{-t}\rmd t =-\sum_{k=1}^{\infty} \int_{0}^{\infty} \rme^{-t} \frac{t^{k-1}}{k\,k!} \rmd t= -\sum_{k=1}^{\infty} \frac{\Gamma\brc{k}}{k\,k!}=-\frac{\pi^2}{6}.}
And the second part is given by the derivative rule for $\text{EL}_{1,1}$ function:
\eq{\int_{0}^{z} \frac{\rme^{t}}{t} \, \text{EL}_{1}\brc{-t}\rmd t=-\text{EL}_{1,1}\brc{z}}
Summing up all the contributions, we get
\eqlb{eq:ELi11_sub}{\text{ELi}_{1,1}\brc{z}=\text{EL}_{1,1}\brc{z}+\text{ELi}_{2}\brc{z}-\brc{\gamma +\log z }\text{ELi}_{1}\brc{z} + \frac{\pi^2}{6}.}
Then, for $\text{ELi}_{2,1}\brc{z}$, we have:
\eqlb{eq:ELi21_sub}{\text{ELi}_{2,1}\brc{z} = \text{EL}_{2,1}\brc{z}+2\,\text{ELi}_{3}\brc{z}-\brc{\gamma +\log z }\text{ELi}_{2}\brc{z} + \frac{\pi^2}{6}\bsq{\gamma+ \log\brc{-z}} -\zeta\brc{3}.}
To derive the latter expression, we use the definition
\eq{\text{ELi}_{2,1}\brc{z} \equiv \int_{-\infty}^{z} \text{ELi}_{1,1}\brc{t} \frac{\rmd t}{t} =
\int_{-\infty}^{z}\brc{\text{EL}_{1,1}\brc{t}+\text{ELi}_{2}\brc{t}-\brc{\gamma +\log t }\text{ELi}_{1}\brc{t} + \frac{\pi^2}{6}}
\frac{\rmd t}{t}}
and split the integral into two parts:
\begin{itemize}
\item[1.] $\displaystyle \int_{-\infty}^{z}\brc{\text{ELi}_{2}\brc{t}-\brc{\gamma +\log t }\text{ELi}_{1}\brc{t}}\frac{\rmd t}{t}$,
\item[2.] $\displaystyle \int_{-\infty}^{z}\brc{\text{EL}_{1,1}\brc{t} + \frac{\pi^2}{6}}
\frac{\rmd t}{t}$,
\end{itemize}
where the integrands in both cases tend to zero as $t \rightarrow - \infty$. In the first case, we integrate by parts:
\eq{\text{ELi}_{3}\brc{t}-\int_{-\infty}^{z}\brc{\gamma +\log t }\rmd\, \text{ELi}_{2}\brc{t} =\text{ELi}_{3}\brc{t} -\brc{\gamma +\log z} \text{ELi}_{2}\brc{z} + \int_{-\infty}^{z} \text{ELi}_{2}\brc{t}\frac{\rmd t}{t}=}
\eqn{ =2\,\text{ELi}_{3}\brc{z} -\brc{\gamma +\log z} \text{ELi}_{2}\brc{z}.}
In the second case, we split the integration domain and then integrate by parts:
\eq{\int_{-\infty}^{-1}\brc{\text{EL}_{1,1}\brc{t} + \frac{\pi^2}{6}} \rmd \, \log\brc{-t} =
\int_{-\infty}^{-1} \rme^t \, \log\brc{-t} \text{EL}_{1}\brc{-t} \frac{\rmd t}{t} ,}
\eq{\int_{-1}^{z} \text{EL}_{1,1}\brc{t} \frac{\rmd t}{t} = \int_{-1}^{0}\text{EL}_{1,1}\brc{t} \rmd \, \log\brc{-t} +
\int_{0}^{z}\text{EL}_{1,1}\brc{t} \frac{\rmd t}{t}=}
\eqn{= \int_{-1}^{0} \rme^t \, \log\brc{-t} \text{EL}_{1}\brc{-t} \frac{\rmd t}{t} +
\text{EL}_{2,1}\brc{z}, }
\eq{\int_{-1}^{z} \frac{\pi^2}{6} \, \frac{\rmd t}{t} = \frac{\pi^2}{6} \log\brc{-z}.}
Putting everything back together, we arrive at
\eqlb{eq:ELi21_drv}{\text{ELi}_{2,1}\brc{z} =\text{EL}_{2,1}\brc{z} + 2\, \text{ELi}_{3}\brc{z}
-\brc{\gamma +\log z} \text{ELi}_{2}\brc{z} + \frac{\pi^2}{6} \log\brc{-z} -
\int_{0}^{\infty} \rme^{-t} \, \log\brc{t} \text{EL}_{1}\brc{t} \frac{\rmd t}{t},}
where the change of integration variable $t\rightarrow -t$ was also performed. The last step is to compute the constant represented by the integral in (\ref{eq:ELi21_drv}). To do so, we use power series expansion of $\text{EL}_{1}\brc{t}$ and integrate each term of the series:
\eq{\frac{1}{k\,k!} \int_{0}^{\infty} \rme^{-t} \, \log\brc{t} t^{k-1} \rmd \, t =\frac{ \Gamma\brc{k} \psi\brc{k}}{k\,k!} =
\frac{1}{k^2} \brc{-\gamma +\sum_{m=1}^{k-1} \frac{1}{m}}.}
Taking the sum over $k$ gives
\eq{\int_{0}^{\infty} \rme^{-t} \, \log\brc{t} \text{EL}_{1}\brc{t} \frac{\rmd t}{t} = - \gamma\, \frac{\pi^2}{6} +\zeta\brc{3},}
which completes the derivation of (\ref{eq:ELi21_sub}).

Repeating the same steps, we derive a more general expression:
\eqlb{eq:ELin1_sub}{\text{ELi}_{n,1}\brc{z} = \text{EL}_{n,1}\brc{z}+n\,\text{ELi}_{n+1}\brc{z}-\brc{\gamma +\log z }\text{ELi}_{n}\brc{z} + \sum_{k=1}^{n}\frac{\brc{-1}^{k-1}}{\brc{k-1}! \brc{n-k}!}\, \text{cLi}_{k,1} \log\brc{-z}^{n-k},}
where the constants $\text{cLi}_{k,1}$ are the following integrals
\eq{k\geq 1:\quad \text{cLi}_{k,1}=\int_{0}^{\infty} \rme^{-t} \, \log\brc{t}^{k-1} \text{EL}_{1}\brc{t} \frac{\rmd t}{t},}
which we will compute in Sec.~\ref{sec:leveln_int}.
The proof of (\ref{eq:ELin1_sub}) is a straightforward induction with a base case given by (\ref{eq:ELi11_sub}).

Finally, we are ready to discuss the rest of the level $2$ functions, starting with the $\text{ELi}_{1,n}\brc{z}$. By definition, we have
\eqlb{eq:def_ELi1n}{\text{ELi}_{1,n}\brc{z}=-\int_{-\infty}^{z}\frac{\rme^{t}}{t}\, \text{ELi}_{n}\brc{-t} \rmd t.}
Using (\ref{eq:ELn_gen}), we split the integral into two parts:
\begin{itemize}
\item[1.] $\displaystyle -\int_{-\infty}^{z} \frac{\rme^{t}}{t}\, \text{EL}_{n}\brc{-t} \rmd t$,
\item[2.] $\displaystyle -\int_{-\infty}^{z} \frac{\rme^{t}}{t} \sum_{k=0}^{n} \frac{\brc{-1}^k}{k! \brc{n-k}!}\,
\Gamma^{\brc{k}}\brc{1} \log\brc{t}^{n-k} \rmd t$.
\end{itemize}
The first part is straightforward:
\eq{-\int_{-\infty}^{z} \frac{\rme^{t}}{t}\, \text{EL}_{n}\brc{-t} \rmd t = -\int_{-\infty}^{0} \frac{\rme^{t}}{t}\, \text{EL}_{n}\brc{-t} \rmd t -\int_{0}^{z} \frac{\rme^{t}}{t}\, \text{EL}_{n}\brc{-t} \rmd t = \text{EL}_{1,n}\brc{z} + \text{cLi}_{1,n},}
where
\eq{\text{cLi}_{1,n}=\int_{0}^{\infty} \rme^{-t}\, \text{EL}_{n}\brc{t} \frac{\rmd t}{t} =
\sum_{k=1}^{\infty} \frac{1}{k^{n+1}}=\zeta\brc{n+1}.}
To compute the second part,  we integrate each term in the sum by parts:
\eq{-\int_{-\infty}^{z} \frac{\rme^{t}}{t} \log\brc{t}^{n-k} \rmd t = -\int_{-\infty}^{z} \log\brc{t}^{n-k} \rmd \,\text{ELi}_{1}\brc{t}=
-\log\brc{z}^{n-k} \text{ELi}_{1}\brc{t} }
\eq{+\brc{n-k}\int_{-\infty}^{z} \log\brc{t}^{n-k-1} \rmd \,\text{ELi}_{2}\brc{t}.}
We repeat this until the power of $\log\brc{t}$ is reduced to zero, which gives
\eqlb{eq:exp_log_int}{-\int_{-\infty}^{z} \frac{\rme^{t}}{t} \log\brc{t}^{n-k} \rmd t = \sum_{m=0}^{n-k} \brc{-1}^{m+1} \frac{\brc{n-k}!}{\brc{n-k-m}!} \log\brc{z}^{n-k-m} \text{ELi}_{m+1}\brc{z}.}
Putting everything back together, we arrive at
\eqlb{eq:ELi1n_sub}{\text{ELi}_{1,n}\brc{z}= \text{EL}_{1,n}\brc{z} + \zeta\brc{n+1} + \sum_{k=0}^{n}\sum_{m=k}^{n}
\frac{\brc{-1}^{m+1}}{k!\brc{n-m}!} \, \Gamma^{\brc{k}}\brc{1} \log\brc{z}^{n-m} \text{ELi}_{m-k+1}\brc{z}.}
We increase the first index of $\text{ELi}_{1,n}\brc{z}$ by repeatedly integrating and following the same steps as in the case of $\text{ELi}_{2,1}\brc{z}$, which gives the following general relation for level $2$ multiple polyexponential integrals:
\begin{equation}
\label{eq:ELimn_sub}
\boxed{
\begin{aligned}
\text{ELi}_{m,n}\brc{z}= \,&\text{EL}_{m,n}\brc{z} + \sum_{k=1}^{m} \frac{\brc{-1}^{k-1}}{\brc{k-1}!\brc{m-k}!} \,
\text{cLi}_{k,n} \log\brc{-z}^{m-k}\\
& +\sum_{k=0}^{n}\sum_{j=k}^{n}
\frac{\brc{-1}^{j+1}}{k!\brc{n-j}!} \, \tbinom{j+m-k-1}{m-1} \Gamma^{\brc{k}}\brc{1} \log\brc{z}^{n-j} \text{ELi}_{j+m-k}\brc{z}.
\end{aligned}}
\end{equation}
The constants $\text{cLi}_{k,n}$ are in agreement with the previous definitions and are given by the following integrals:
\eqlb{eq:cLimn_def}{k\geq1:\quad \text{cLi}_{k,n}= \int_{0}^{\infty} \rme^{-t} \log\brc{t}^{k-1}\, \text{EL}_{n}\brc{t} \frac{\rmd t}{t}.}
Again, we postpone the computation of these constants to Sec.~\ref{sec:leveln_int}.
To prove (\ref{eq:ELimn_sub}), we apply an induction with a base case (\ref{eq:ELi1n_sub}) (see Appendix \ref{app:rel23_proof}).

\subsection{Level \texorpdfstring{$3$}{} integrals}
Before we move on to an arbitrary level, let's discuss level $3$ multiple polyexponential integrals of the form $\text{ELi}_{l,m,n}\brc{z}$.
As usual, we start with the definition of $\text{ELi}_{1,m,n}\brc{z}$:
\eq{\text{ELi}_{1,m,n}\brc{z} = -\int_{-\infty}^{z}\frac{\rme^{t}}{t}\, \text{ELi}_{m,n}\brc{-t} \rmd t.}
Substituting the previous result (\ref{eq:ELimn_sub}), we split the integral into three parts:
\begin{itemize}
\item[1.] $\displaystyle -\int_{-\infty}^{z} \frac{\rme^{t}}{t}\, \text{EL}_{m,n}\brc{-t} \rmd t$,
\item[2.] $\displaystyle -\int_{-\infty}^{z} \frac{\rme^{t}}{t} \sum_{k=1}^{m} \frac{\brc{-1}^{k-1}}{\brc{k-1}!\brc{m-k}!} \,
\text{cLi}_{k,n} \log\brc{t}^{m-k} \rmd t$,
\item[3.] $\displaystyle -\int_{-\infty}^{z} \frac{\rme^{t}}{t} \sum_{k=0}^{n}\sum_{j=k}^{n}
\frac{\brc{-1}^{j+1}}{k!\brc{n-j}!} \, \tbinom{j+m-k-1}{m-1} \Gamma^{\brc{k}}\brc{1} \log\brc{-t}^{n-j} \text{ELi}_{j+m-k}\brc{-t} \rmd t$.
\end{itemize}
The first part is the simplest one to compute:
\eq{-\int_{-\infty}^{z} \frac{\rme^{t}}{t}\, \text{EL}_{m,n}\brc{-t} \rmd t = \text{EL}_{1,m,n}\brc{z}+\text{cLi}_{1,m,n},}
where
\eq{\text{cLi}_{1,m,n}= \int_{0}^{\infty} \rme^{-t}\, \text{EL}_{m,n}\brc{t} \frac{\rmd t}{t}.}
We use the result obtained previously (\ref{eq:exp_log_int})
\eq{-\int_{-\infty}^{z} \rme^{t} \log\brc{t}^{m-k} \frac{\rmd t}{t} = \sum_{j=k}^{m} \brc{-1}^{j-k+1} \frac{\brc{m-k}!}{\brc{m-j}!} \log\brc{z}^{m-j} \text{ELi}_{j-k+1}\brc{z}}
to compute the second part:
\begin{equation}
\begin{aligned}
&-\int_{-\infty}^{z} \frac{\rme^{t}}{t} \sum_{k=1}^{m} \frac{\brc{-1}^{k-1}}{\brc{k-1}!\brc{m-k}!} \,
\text{cLi}_{k,n} \log\brc{t}^{m-k} \rmd t =\\ &\sum_{k=1}^{m} \sum_{j=k}^{m} \frac{\brc{-1}^{j}}{\brc{k-1}!\brc{m-j}!} \, \text{cLi}_{k,n} \log\brc{z}^{m-j} \text{ELi}_{j-k+1}\brc{z}.
\end{aligned}
\end{equation}
The new type of integral appears in the third part:
\begin{equation}
\begin{aligned}
-&\int_{-\infty}^{z} \rme^{t} \log\brc{-t}^{n-j} \text{ELi}_{j+m-k}\brc{-t} \frac{\rmd t}{t} =
\int_{-\infty}^{z} \log\brc{-t}^{n-j} \rmd \,\text{ELi}_{1,j+m-k}\brc{t} =\\
&\log\brc{-z}^{n-j} \text{ELi}_{1,j+m-k}\brc{z}-\brc{n-j} \int_{-\infty}^{z} \log\brc{-t}^{n-j-1} \rmd \,\text{ELi}_{2,j+m-k}\brc{t}.
\end{aligned}
\end{equation}
We continue integrating by parts until the power of the logarithm is reduced to zero:
\eq{-\int_{-\infty}^{z} \rme^{t} \log\brc{-t}^{n-j} \text{ELi}_{j+m-k}\brc{-t} \frac{\rmd t}{t} =
\sum_{i=j}^{n} \brc{-1}^{i-j} \frac{\brc{n-j}!}{\brc{n-i}!} \log\brc{-z}^{n-i} \text{ELi}_{i-j+1,j+m-k}\brc{z}.}
The latter gives the following answer for the third integral:
\begin{equation}
\begin{aligned}
\sum_{k=0}^{n}\sum_{j=k}^{n}\sum_{i=j}^{n}
\frac{\brc{-1}^{i+1}}{k!\brc{n-i}!} \, \binom{j+m-k-1}{m-1} \Gamma^{\brc{k}}\brc{1} \log\brc{-z}^{n-i} \text{ELi}_{i-j+1,j+m-k}\brc{z}.
\end{aligned}
\end{equation}
Putting all three results together, we get
\eq{\text{ELi}_{1,m,n}\brc{z} = \text{EL}_{1,m,n}\brc{z}+\text{cLi}_{1,m,n} +
\sum_{k=1}^{m} \sum_{j=k}^{m} \frac{\brc{-1}^{j}}{\brc{k-1}!\brc{m-j}!} \, \text{cLi}_{k,n} \log\brc{z}^{m-j} \text{ELi}_{j-k+1}\brc{z}}
\eqn{+\sum_{k=0}^{n}\sum_{j=k}^{n}\sum_{i=j}^{n}
\frac{\brc{-1}^{i+1}}{k!\brc{n-i}!} \, \tbinom{j+m-k-1}{m-1} \Gamma^{\brc{k}}\brc{1} \log\brc{-z}^{n-i} \text{ELi}_{i-j+1,j+m-k}\brc{z}.}
Now, we recursively apply the definition of $\text{ELi}_{l,m,n}\brc{z}$ for $l\geq 2$
\eq{\text{ELi}_{l,m,n}\brc{z}= \int_{-\infty}^{z} \text{ELi}_{l-1,m,n}\brc{t} \frac{\rmd t}{t}}
to derive the general relation for level $3$ multiple polyexponential integrals:
\eqn{\text{ELi}_{l,m,n}\brc{z} = \text{EL}_{l,m,n}\brc{z}+ \sum_{k=1}^{l} \frac{\brc{-1}^{k-1}}{\brc{k-1}!\brc{l-k}!} \,
\text{cLi}_{k,m,n} \log\brc{-z}^{l-k} }
\eqlb{eq:ELi_lmn_sub}{+\sum_{k=1}^{m} \sum_{j=k}^{m} \frac{\brc{-1}^{j}}{\brc{k-1}!\brc{m-j}!} \, \tbinom{l+j-k-1}{l-1} \text{cLi}_{k,n} \log\brc{z}^{m-j} \text{ELi}_{l+j-k}\brc{z}}
\eqn{+\sum_{k=0}^{n}\sum_{j=k}^{n}\sum_{i=j}^{n}
\frac{\brc{-1}^{i+1}}{k!\brc{n-i}!} \, \tbinom{l+i-j-1}{l-1} \tbinom{m+j-k-1}{m-1} \Gamma^{\brc{k}}\brc{1} \log\brc{-z}^{n-i} \text{ELi}_{l+i-j,m+j-k}\brc{z},}
where the new constants $\text{cLi}_{k,m,n}$ are defined as
\eqlb{eq:cLi_kmn_def}{k\geq1:\quad \text{cLi}_{k,m,n}= \int_{0}^{\infty} \rme^{-t} \log\brc{t}^{k-1}\, \text{EL}_{m,n}\brc{t} \frac{\rmd t}{t}.}
For the proof of (\ref{eq:ELi_lmn_sub}), see Appendix \ref{app:rel23_proof}.

\subsection{Level \texorpdfstring{$n\geq 2$ }{} integrals}
\label{sec:leveln_int}
First, we rewrite (\ref{eq:ELi_lmn_sub}) by replacing indices $l$, $m$, $n$ with $s_1$, $s_2$, $s_3$ and rearranging the order of summations:
\eqn{\text{ELi}_{s_1,s_2,s_3}\brc{z} = \text{EL}_{s_1,s_2,s_3}\brc{z}+ \sum_{k_1=1}^{s_1} \frac{\brc{-1}^{k_1-1}}{\brc{k_1-1}!\brc{s_1-k_1}!} \, \text{cLi}_{k_1,s_2,s_3} \log\brc{-z}^{s_1-k_1} }
\eq{+\sum_{s_2\geq k_1\geq k_2\geq 1} \frac{\brc{-1}^{k_1}}{\brc{k_2-1}!\brc{s_2-k_1}!} \, \tbinom{s_1-1+k_1-k_2}{s_1-1}\, \text{cLi}_{k_2,s_3} \log\brc{z}^{s_2-k_1} \text{ELi}_{s_1+k_1-k_2}\brc{z}}
\eqn{+\sum_{s_3\geq k_1\geq k_2\geq k_3\geq 0} \frac{\brc{-1}^{k_1+1}}{k_3!\brc{s_3-k_1}!}\prod_{j=1}^{2} \tbinom{s_j-1+k_j-k_{j+1}}{s_j-1} \, \Gamma^{\brc{k_3}}\brc{1} \log\brc{-z}^{s_3-k_1} \text{ELi}_{s_1+k_1-k_2,s_2+k_2-k_3}\brc{z},}
where the following notations were used:
\eq{\sum_{s_2\geq k_1\geq k_2\geq 1} = \sum_{k_1=1}^{s_2} \sum_{k_2=1}^{k_1},\quad
\sum_{s_3\geq k_1\geq k_2\geq k_3\geq 0}= \sum_{k_1=0}^{s_3} \sum_{k_2=0}^{k_1} \sum_{k_3=0}^{k_2}.}
The generalization to an arbitrary level $n\geq 2$ is
\begin{equation}
\label{eq:ELi_sub_gen}
\begin{aligned}
\text{ELi}_{s_1,\dots,s_n}\brc{z} = \,&\text{EL}_{s_1,\dots,s_n}\brc{z}+ \sum_{k_1=1}^{s_1} \frac{\brc{-1}^{k_1-1}}{\brc{k_1-1}!\brc{s_1-k_1}!} \, \text{cLi}_{k_1,s_2,\dots,s_n} \log\brc{-z}^{s_1-k_1} +\\
&\sum_{i=2}^{n-1}\sum_{s_i\geq k_1\geq \dots\geq k_i\geq 1} \frac{\brc{-1}^{k_1}}{\brc{k_i-1}!\brc{s_i-k_1}!} \prod_{j=1}^{i-1} \tbinom{s_j-1+k_j-k_{j+1}}{s_j-1}\, \text{cLi}_{k_i,s_{i+1},\dots,s_n} \times\\
&\log(\brc{-1}^i z)^{s_i-k_1} \text{ELi}_{s_1+k_1-k_2,\dots,s_{i-1}+k_{i-1}-k_i}\brc{z}+\\
&\sum_{s_n\geq k_1\geq\dots\geq k_n\geq 0} \frac{\brc{-1}^{k_1+1}}{k_n!\brc{s_n-k_1}!}\prod_{j=1}^{n-1} \tbinom{s_j-1+k_j-k_{j+1}}{s_j-1} \, \Gamma^{\brc{k_n}}\brc{1}  \times\\
&\log\brc{\brc{-1}^n z}^{s_n-k_1} \text{ELi}_{s_1+k_1-k_2,\dots,s_{n-1}+k_{n-1}-k_n}\brc{z},
\end{aligned}
\end{equation}
and the constants $\text{cLi}_{k_1,s_2,\dots,s_n}$ are defined as
\eqlb{eq:cLi_gen_def}{k_1\geq1:\quad \text{cLi}_{k_1,s_2,\dots,s_n}= \int_{0}^{\infty} \rme^{-t} \log\brc{t}^{k_1-1}\, \text{EL}_{s_2,\dots,s_n}\brc{t} \frac{\rmd t}{t}.}
We delegate the proof of (\ref{eq:ELi_sub_gen}) to the Appendix \ref{app:reln_proof} and focus here on computing the constants $\text{cLi}_{s_1,\dots,s_n}$.

At level $2$, we have:
\eq{\text{cLi}_{m,n}= \int_{0}^{\infty} \rme^{-t} \log\brc{t}^{m-1}\, \text{EL}_{n}\brc{t} \frac{\rmd t}{t}.}
Substituting the Taylor expansion of $\text{EL}_{n}\brc{t}$,
\eq{\text{EL}_{n}\brc{t}=\sum_{k=1}^{\infty} \frac{t^k}{k^n k!},}
and using the integral form of the gamma function
\eq{\Gamma\brc{k}= \int_{0}^{\infty} \rme^{-t} \, t^{k-1} \rmd t,\quad
\Gamma^{(m-1)}\brc{k}= \int_{0}^{\infty} \rme^{-t} \log\brc{t}^{m-1} t^{k-1} \rmd t,}
we get
\eqlb{eq:cLi_l2_Gamma}{\text{cLi}_{m,n}= \sum_{k=1}^{\infty} \frac{1}{k^{n+1}} \,\frac{\Gamma^{(m-1)}\brc{k}}{\Gamma\brc{k}}.}
The ratio $\Gamma^{(m-1)}\brc{k}/ \Gamma\brc{k}$ can be expressed in terms of polygamma functions. For some first values of $m\geq 1$ one has:
\eq{\frac{\Gamma^{\brc{1}}\brc{k}}{\Gamma\brc{k}} = \psi\brc{k},}
\eq{\frac{\Gamma^{\brc{2}}\brc{k}}{\Gamma\brc{k}} = \psi^{\brc{1}}\brc{k}+\psi\brc{k}^2,}
\eq{\frac{\Gamma^{\brc{3}}\brc{k}}{\Gamma\brc{k}} = \psi^{\brc{2}}\brc{k}+3\,\psi\brc{k}\psi^{\brc{1}}\brc{k}+\psi\brc{k}^3.}
The following recurrence relation can be used to compute the latter ratios for general values of $m$:
\eqlb{eq:Gamma_rec}{\frac{\Gamma^{\brc{m}}\brc{x}}{\Gamma\brc{x}} = \partial_x \brc{\frac{\Gamma^{\brc{m-1}}\brc{x}}{\Gamma\brc{x}}}
+ \psi\brc{x} \, \frac{\Gamma^{\brc{m-1}}\brc{x}}{\Gamma\brc{x}}.}
Since $k$ is a positive integer, the polygamma functions can be written as the following sums:
\eq{\psi\brc{k}= - \gamma + \sum_{k_1=1}^{k-1}\frac{1}{k_1},}
\eq{l\geq 1:\quad \psi^{\brc{l}}\brc{k}=\brc{-1}^{l+1} l! \sum_{k_1\geq k} \frac{1}{k_1^{l+1}}.}
Thus, the constants $\text{cLi}_{m,n}$ can be described using multiple zeta values or MZVs (the same is valid for level $n\geq2$ constants).
For $m=1,2,3$, we have
\eq{\text{cLi}_{1,n}= \zeta\brc{n+1},\quad \text{cLi}_{2,n}= \zeta\brc{n+1,1} -\gamma \, \zeta\brc{n+1},}
\eq{\text{cLi}_{3,n}= \zeta\brc{n+3}+\zeta\brc{2,n+1}+\zeta\brc{n+1,2} +2\, \zeta\brc{n+1,1,1} -2 \,\gamma\, \zeta\brc{n+1,1} +
\gamma^2 \, \zeta\brc{n+1}.}
For higher values of $m$, we use (\ref{eq:Gamma_rec}) to express the ratios of gamma functions in terms of the polygamma functions. Since each ratio is a polynomial in $\psi\brc{k}$ and $\psi^{\brc{l}}\brc{k}$, the computation of (\ref{eq:cLi_l2_Gamma}) boils down to calculating the following sums:
\eqlb{eq:cpsi_gen}{\sum_{k=1}^{\infty} \frac{\psi\brc{k}^r}{k^{n+1}} \psi^{\brc{l_1}}\brc{k}\dots \psi^{\brc{l_i}}\brc{k}}
with $r\geq 0$ and $l_j\geq 1$, $j=1,\dots,i$ (some or all values $l_j$ may coincide).
First, let's consider the case with just $\psi\brc{k}^r$:
\eq{\sum_{k=1}^{\infty} \frac{\psi\brc{k}^r}{k^{n+1}} = \sum_{k=1}^{\infty}\sum_{j=0}^{r} \frac{1}{k^{n+1}}
\binom{r}{j} \brc{-\gamma}^{r-j}  \brc{\sum_{k_1=1}^{k-1}\frac{1}{k_1}}^j.}
The latter can be simplified using the sum over all ordered partitions of $j\geq 1$ and the multinomial coefficients:
\eq{\brc{\sum_{k_1=1}^{k-1}\frac{1}{k_1}}^j=\sum_{k_1=1}^{k-1}\dots \sum_{k_j=1}^{k-1} \frac{1}{k_1}\dots \frac{1}{k_j} =
\sum_{\brc{v_1,\dots,v_i}\in \ops\brc{j}} \sum_{k>k_1>\dots>k_i\geq 1 } \binom{j}{v_1\dots v_i}
\frac{1}{k_1^{v_1} \dots k_i^{v_i}}.}
This gives the following answer in terms of MZVs:
\eqlb{eq:cpsi_dig_res}{\sum_{k=1}^{\infty} \frac{\psi\brc{k}^r}{k^{n+1}} = \sum_{j=0}^{r} \sum_{\ops\brc{j}} \brc{-\gamma}^{r-j}\binom{r}{j}
\binom{j}{\ops\brc{j}} \zeta\brc{n+1, \ops\brc{j}},}
where for $j=0$ we have
\eq{\binom{0}{\ops\brc{0}} \zeta\brc{n+1, \ops\brc{0}}= \zeta\brc{n+1}.}
Next, we consider the case with a single polygamma function $\psi^{\brc{l}}\brc{k}$, $l\geq 1$:
\eq{\sum_{k=1}^{\infty} \frac{\psi^{\brc{l}}\brc{k}}{k^{n+1}} =\brc{-1}^{l+1} l! \sum_{k=1}^{\infty}
\sum_{k_1= k}^{\infty} \frac{1}{k^{n+1} k_1^{l+1}}=\brc{-1}^{l+1} l! \bsq{\zeta\brc{n+l+2}+\zeta\brc{l+1,n+1}}.}
In the case where there is an array of polygamma functions $\psi^{\brc{l_j}}\brc{k}$
\eqlb{eq:cpsi_polyg}{\sum_{k=1}^{\infty} \frac{1}{k^{n+1}} \psi^{\brc{l_1}}\brc{k}\dots \psi^{\brc{l_i}}\brc{k},}
we will use the stuffle product $\star$ for indices (see \cite{Wald}). If we have two vectors $\mathbf{s_1}$ and $\mathbf{s_2}$, then the stuffle product is a multiset of vectors obtained by merging $\mathbf{s_1}$ and $\mathbf{s_2}$ in a specific way: each element of $\mathbf{s_1}$ can be put before or after any element in $\mathbf{s_2}$ and can also be added to any element in $\mathbf{s_2}$ with a condition that the initial ordering of elements in $\mathbf{s_1}$ is preserved in the final string. For example, consider $\mathbf{s_1}=\brc{1}$ and $\mathbf{s_2}=\brc{2,3}$:
\eq{\brc{1}\star \brc{2,3} = \bfi{\brc{1,2,3},\brc{2,1,3},\brc{2,3,1},\brc{2,4},\brc{3,3}}.}
The stuffle product is defined to deal with the product of multiple polylogarythm functions and in our case it can be used as follows:
\eq{\sum_{k_1=k}^{\infty}\dots \sum_{k_i=k}^{\infty} \frac{1}{k_1^{l_1+1}}\dots \frac{1}{k_i^{l_i+1}} =
\sum_{\brc{v_1,\dots,v_j}\in \brc{l_1+1}\star\dots\star\brc{l_i+1}} \sum_{k_1>\dots>k_j\geq k} \frac{1}{k_1^{v_1}}\dots \frac{1}{k_j^{v_j}}.}
Substituting this back into (\ref{eq:cpsi_polyg}) gives us an expression in terms of MZVs:
\begin{equation}\label{eq:cpsi_polyg_res}
\begin{aligned}
\sum_{k=1}^{\infty} \frac{1}{k^{n+1}} \psi^{\brc{l_1}}\brc{k}\dots \psi^{\brc{l_i}}\brc{k} = \brc{-1}^{\sum_{j=1}^{i}\brc{l_j+1}}
l_1!\dots l_i! \sum_{\textbf{v}\in \brc{l_1+1}\star\dots\star\brc{l_i+1}} \bsq{\zeta\brc{\textbf{v}\oplus n+1} +\zeta\brc{\textbf{v}, n+1}}.
\end{aligned}
\end{equation}
Now, the two results (\ref{eq:cpsi_dig_res}) and (\ref{eq:cpsi_polyg_res}) can be combined to get the answer for the general expression (\ref{eq:cpsi_gen}) for $i\geq 1$:
\begin{equation}\label{eq:cpsi_n_res}
\begin{aligned}
&\sum_{k=1}^{\infty} \frac{\psi\brc{k}^r}{k^{n+1}} \psi^{\brc{l_1}}\brc{k}\dots \psi^{\brc{l_i}}\brc{k}=\\
&\brc{-1}^{\sum_{j=1}^{i}\brc{l_j+1}} l_1!\dots l_i! \sum_{j=0}^{r} \sum_{\ops\brc{j}}
\sum_{\textbf{v}\in \brc{l_1+1}\star\dots\star\brc{l_i+1}} \brc{-\gamma}^{r-j}\binom{r}{j} \binom{j}{\ops\brc{j}} \bsq{\zeta\brc{\textbf{v}\oplus n+1, \ops\brc{j}}+\zeta\brc{\textbf{v}, n+1, \ops\brc{j}}}.
\end{aligned}
\end{equation}
The latter computation can be further optimized by collecting all coincident $l_j$'s and using the following:
\eq{\sum_{k=1}^{\infty} \frac{\psi^{\brc{l}}\brc{k}^r}{k^{n+1}}  = \brc{-1}^{r\brc{l+1}} \brc{l!}^r
\sum_{\textbf{v}\in\ops\brc{r}} \binom{r}{\textbf{v}}\bsq{\zeta\brc{\brc{l+1}\cdot\textbf{v}\oplus n+1}+\zeta\brc{\brc{l+1}\cdot\textbf{v},n+1}},}
where $\brc{l+1}\cdot\textbf{v}$ is a product of a scalar and a vector.

In the most general case (\ref{eq:cLi_gen_def}) with $n\geq3$, the constants are defined in terms of the dressed functions $\text{EL}_{s_2,\dots,s_n}\brc{z}$. Since we want the final result to be represented in terms of MZVs, we will define another set of constants in terms of the undressed functions $el_{s_2,\dots,s_m}\brc{z}$:
\eqlb{eq:cli_def}{ cli_{s_1,\dots,s_m}= \int_{0}^{\infty} \rme^{-t} \log\brc{t}^{s_1-1}\, el_{s_2,\dots,s_m}\brc{t} \frac{\rmd t}{t}.}
Due to relations (\ref{eq:EL_rel_even}) and (\ref{eq:EL_rel_odd}) between the dressed and undressed functions, the constants $\text{cLi}_{s_1,\dots,s_n}$ can be written as the following linear combinations of the constants $cli_{s_1,\dots,s_m}$:
\eq{\text{cLi}_{s_0,r_1,s_1,\dots,r_k,s_k} = \sum_{\ops\brc{s_1+1}}\dots\sum_{\ops\brc{s_k+1}} cli_{s_0,\brc{r_1-1}\oplus\ops\brc{s_1+1}\oplus\dots\oplus\brc{r_k-1}\oplus \ops\brc{s_k+1}},}
\eq{\text{cLi}_{s_0,r_1,s_1,\dots,r_k,s_k,r_{k+1}} = \sum_{\ops\brc{s_1+1}}\dots\sum_{\ops\brc{s_k+1}} cli_{s_0,\brc{r_1-1}\oplus\ops\brc{s_1+1}\oplus\dots\oplus\brc{r_k-1}\oplus \ops\brc{s_k+1}\oplus r_{k+1}}.}
Substituting the Taylor expansion (\ref{eq:polyExp_def}) into (\ref{eq:cli_def}), we get
\eq{cli_{s_0,s_1,\dots,s_n}=\sum_{k_1>k_2>\dots>k_n \geq 1} \frac{1}{k_1^{s_1+1} k_2^{s_2} \dots k_n^{s_n}} \,\frac{\Gamma^{(s_0-1)}\brc{k_1}}{\Gamma\brc{k_1}}.}
Thus, we will again express the ratio $\Gamma^{(s_0-1)}\brc{k_1}/ \Gamma\brc{k_1}$ in terms of polygamma functions and compute the following sums:
\eq{ \sum_{k_1>k_2>\dots>k_n \geq 1} \frac{\psi\brc{k_1}^r}{k_1^{s_1+1} k_2^{s_2} \dots k_n^{s_n}}\,
\psi^{\brc{l_1}}\brc{k_1}\dots \psi^{\brc{l_i}}\brc{k_1},}
where $l_j\geq 1$ and some or all values $l_j$ may coincide. This computation is similar to the one for level $2$ constants. The main difference is when considering powers of the digamma function $\psi\brc{k_1}^r$, the stuffle product with indices $\brc{s_2,\dots,s_n}$ needs to be introduced:
\eqlb{eq:cpsi_sn_dig_res}{\sum_{k_1>k_2>\dots>k_n \geq 1} \frac{\psi\brc{k_1}^r}{k_1^{s_1+1} k_2^{s_2} \dots k_n^{s_n}} =
\sum_{j=0}^{r} \sum_{\ops\brc{j}}\sum_{\textbf{u}\in \ops\brc{j}\star\brc{s_2,\dots,s_n}} \brc{-\gamma}^{r-j}\binom{r}{j}
\binom{j}{\ops\brc{j}} \zeta\brc{s_1+1, \textbf{u}},}
where, for $j=0$, we have $\ops\brc{0}\star\brc{s_2,\dots,s_n}=\brc{s_2,\dots,s_n}$. If only polygamma functions $\psi^{\brc{l_j}}\brc{k_1}$ with $l_j\geq 1$ are present, the result (\ref{eq:cpsi_polyg_res}) stays virtually the same:
\eqlb{eq:cpsi_sn_polyg_res}{\sum_{k_1>k_2>\dots>k_n \geq 1} \frac{\psi^{\brc{l_1}}\brc{k_1}\dots \psi^{\brc{l_i}}\brc{k_1}}
{k_1^{s_1+1} k_2^{s_2} \dots k_n^{s_n}} =  \brc{-1}^{\sum_{j=1}^{i}\brc{l_j+1}} l_1!\dots l_i! }
\eqn{\times \sum_{\textbf{v}\in \brc{l_1+1}\star\dots\star\brc{l_i+1}} \bsq{\zeta\brc{\textbf{v}\oplus s_1+1,s_2,\dots,s_n} +\zeta\brc{\textbf{v}, s_1+1,s_2,\dots,s_n}}.}
Combining (\ref{eq:cpsi_sn_dig_res}) and (\ref{eq:cpsi_sn_polyg_res}), we arrive at
\begin{equation}
\begin{aligned}
&\sum_{k_1>k_2>\dots>k_n \geq 1} \frac{\psi\brc{k_1}^r}{k_1^{s_1+1} k_2^{s_2} \dots k_n^{s_n}}\,
\psi^{\brc{l_1}}\brc{k_1}\dots \psi^{\brc{l_i}}\brc{k_1}= \brc{-1}^{\sum_{j=1}^{i}\brc{l_j+1}} l_1!\dots l_i! \sum_{j=0}^{r} \sum_{\ops\brc{j}} \sum_{\textbf{u}\in \ops\brc{j}\star\brc{s_2,\dots,s_n}}\\
&\times \sum_{\textbf{v}\in \brc{l_1+1}\star\dots\star\brc{l_i+1}}  \brc{-\gamma}^{r-j}\binom{r}{j} \binom{j}{\ops\brc{j}} \bsq{\zeta\brc{\textbf{v}\oplus s_1+1,\textbf{u}} +\zeta\brc{\textbf{v}, s_1+1,\textbf{u}}}.
\end{aligned}
\end{equation}

\section{Asymptotic expansions for \text{ELi} functions}\label{sec:asymptotics}

In this section, we are going to prove that all multiple polyexponential integrals $\text{ELi}_{s_1,\dots,s_n}\brc{z}$ behave at most like $1/z$ when $z\rightarrow -\infty$.
The starting point is the asymptotic series for the exponential integral:
\begin{equation}\label{asymptoticELi1}
\mathrm{ELi}_1(z)=\frac{\mathrm{e}^z}{z}\left(\sum_{k=0}^{n-1}\frac{k!}{z^k}+R_n(z)\right),
\end{equation}
where $R_n(z)$ is a remainder, which is explicitly given by
\begin{equation}
R_n(z)=n!\,z\,\mathrm{e}^{-z}\int_{-\infty}^z\frac{\mathrm{e}^t}{t^{n+1}}\,\mathrm{d}t.
\end{equation}
Then, we use the following Ansatz for the generic asymptotics of the $\mathrm{ELi}$ functions
\begin{equation}\label{asymptoticansatz}
\begin{aligned}
&\mathrm{ELi}_{s_1,\dots,s_{2n+1}}(z)=\mathrm{e}^z\sum_{j=1}^N\frac{c^{(s_1,\dots,s_{2n+1})}_j}{z^j}+\mathcal{O}\left(\frac{\mathrm{e}^z}{z^{N+1}}\right),\\
&\mathrm{ELi}_{s_1,\dots,s_{2n}}(z)=\sum_{j=1}^N\frac{c^{(s_1,\dots,s_{2n})}_j}{z^j}+\mathcal{O}\left(\frac{1}{z^{N+1}}\right).
\end{aligned}
\end{equation}
From \eqref{asymptoticELi1}, we have
\begin{equation}
c^{(1)}_j=
\Gamma(j).
\end{equation}
The coefficients of the $\mathrm{ELi}$ functions of higher weight satisfy the following recursive relations\footnote{As a remark on the notation, when writing $s_2,\dots,s_{2n}$ all the intermediate weights $s_k$, $2< k< 2n$ are present in the middle.}:
\begin{equation}\label{coeffs>1}
s>1\colon\quad c^{(s)}_j=\begin{cases}
0\quad\quad &\text{for}\ j< s,\\
1\quad\quad &\text{for}\ j= s,\\
c^{(s-1)}_{j-1}+(j-1)\,c^{(s)}_{j-1}\quad &\text{for}\ j>s,
\end{cases}
\end{equation}
\begin{equation}\label{coeffeven1}
n\ge 1\colon\quad c^{(1,s_2,\dots,s_{2n})}_j=
(-1)^j\,\frac{c_j^{(s_2,\dots,s_{2n})}}{j},
\end{equation}
\begin{equation}\label{coeffodd1}
n\ge 1\colon\quad c^{(1,s_2,\dots,s_{2n+1})}_j=
(-1)^j\,c_{j-1}^{(s_2,\dots,s_{2n+1})}+(j-1)\,c_{j-1}^{(1,s_2,\dots,s_{2n+1})},
\end{equation}
\begin{equation}\label{coeffs>1even}
n\ge 1,\ s_1>1\colon\quad c^{(s_1,s_2,\dots,s_{2n})}_j=
-\frac{c_j^{(s_1-1,s_2,\dots,s_{2n})}}{j},
\end{equation}
\begin{equation}\label{coeffs>1odd}
n\ge 1,\ s_1>1\colon\quad c^{(s_1,s_2,\dots,s_{2n+1})}_j=
c_{j-1}^{(s_1-1,s_2,\dots,s_{2n+1})}+(j-1)\,c_{j-1}^{(s_1,s_2,\dots,s_{2n+1})}.
\end{equation}
The above recursion relations can be obtained from the following relations between the $\mathrm{ELi}$ functions:
\begin{equation}
\begin{aligned}
&z\,\frac{\mathrm{d}}{\mathrm{d}z}\,\mathrm{ELi}_{1,s_2,\dots,s_n}(z)=-\mathrm{e}^z\,\mathrm{ELi}_{s_2,\dots,s_n}(-z),\\
&z\,\frac{\mathrm{d}}{\mathrm{d}z}\,\mathrm{ELi}_{s_1,s_2,\dots,s_n}(z)=\mathrm{ELi}_{s_1-1,s_2,\dots,s_n}(z),\quad\text{for}\ s_1\in\mathbb{Z}_{>1},
\end{aligned}
\end{equation}
where all the $\mathrm{ELi}$ functions vanish as $z\to-\infty$, therefore no integration constants appear in the asymptotic expansions.

Let us start by proving \eqref{coeffs>1} by induction on $s>1$. For $s=2$, by using the Ansatz \eqref{asymptoticansatz}, we have
\begin{equation}
z\,\frac{\mathrm{d}}{\mathrm{d}z}\,\mathrm{ELi}_2(z)=\mathrm{e}^z\left(\sum_{j\ge 1}\frac{c^{(2)}_{j}}{z^{j-1}}-\sum_{j\ge 1}j\,\frac{c^{(2)}_j}{z^j}\right).
\end{equation}
On the other hand, the recursive formula for $\mathrm{ELi}_2(z)$ gives
\begin{equation}
z\,\frac{\mathrm{d}}{\mathrm{d}z}\,\mathrm{ELi}_2(z)=\mathrm{ELi}_1(z)=\mathrm{e}^z\sum_{j\ge 1}\frac{(j-1)!}{z^j}.
\end{equation}
Comparing the coefficients in the powers of $1/z$, we find $c_1^{(2)}=0$, $c_2^{(2)}=1$, and
\begin{equation}
c_j^{(2)}=(j-1)\,c_j^{(2)}+(j-2)!\quad\text{for}\quad j\ge 2.
\end{equation}

Suppose the claim holds for every $n\le s-1$, for some $s>2$. By using the Ansatz \eqref{asymptoticansatz},
\begin{equation}
z\,\frac{\mathrm{d}}{\mathrm{d}z}\,\mathrm{ELi}_s(z)=z\,\frac{\mathrm{d}}{\mathrm{d}z}\left(\sum_{j\ge 1}\mathrm{e}^z\frac{c^{(s)}_j}{z^j}\right)=\mathrm{e}^z\left(\sum_{j\ge 1}\frac{c^{(s)}_{j}}{z^{j-1}}-\sum_{j\ge 1}j\,\frac{c^{(s)}_j}{z^j}\right)
\end{equation}
and the relation
\begin{equation}
z\,\frac{\mathrm{d}}{\mathrm{d}z}\,\mathrm{ELi}_s(z)=\mathrm{ELi}_{s-1}(z),
\end{equation}
along with the inductive hypothesis for $\mathrm{ELi}_{s-1}(z)$, we find, comparing the powers of $1/z$,
\begin{equation}
c_j^{(s)}=\begin{cases}
0\quad &\text{for}\quad j<s,\\
c_{j-1}^{(s-1)}=1\quad &\text{for}\quad j=s,\\
c_{j-1}^{(s-1)}+(j-1)c_{j-1}^{(s)}\quad &\text{for}\quad j>s.
\end{cases}
\end{equation}

For the proof of \eqref{coeffeven1}, we start with
\begin{equation}
z\,\frac{\mathrm{d}}{\mathrm{d}z}\,\mathrm{ELi}_{1,s_2,\dots,s_{2n}}(z)=-\mathrm{e}^z\,\mathrm{ELi}_{s_2,\dots,s_{2n}}(-z)
\end{equation}
and use the Ansatz \eqref{asymptoticansatz} to get
\begin{equation}
\begin{aligned}
&z\,\frac{\mathrm{d}}{\mathrm{d}z}\,\mathrm{ELi}_{1,s_2,\dots,s_{2n}}(z)=z\,\frac{\mathrm{d}}{\mathrm{d}z}\left(\sum_{j\ge 1}\frac{c_j^{(1,s_2,\dots,s_{2n})}}{z^j}\right)=\sum_{j\ge 1}-j\,\frac{c_j^{(1,s_2,\dots,s_{2n})}}{z^j},\\
&-\mathrm{e}^z\,\mathrm{ELi}_{s_2,\dots,s_{2n}}(-z)=-\mathrm{e}^z\left(\mathrm{e}^{-z}\sum_{j\ge 1}\frac{c_j^{(s_2,\dots,s_{2n})}}{(-z)^j}\right)=\sum_{j\ge 1}(-1)^{j+1}\frac{c_j^{(s_2,\dots,s_{2n})}}{z^j}.
\end{aligned}
\end{equation}
Comparing the powers of $1/z$, we conclude
\begin{equation}
c_j^{(1,s_2,\dots,s_{2n})}=\frac{(-1)^j}{j}\,c_j^{(s_2,\dots,s_{2n})}.
\end{equation}

Analogously, for the proof of \eqref{coeffodd1},  \eqref{coeffs>1even}, and \eqref{coeffs>1odd} we use \eqref{asymptoticansatz} together with the following relations:
\begin{equation}
\begin{aligned}
z\,\frac{\mathrm{d}}{\mathrm{d}z}\,\mathrm{ELi}_{1,s_2,\dots,s_{2n+1}}(z) &=-\mathrm{e}^z\,\mathrm{ELi}_{s_2,\dots,s_{2n+1}}(-z),\\
z\,\frac{\mathrm{d}}{\mathrm{d}z}\,\mathrm{ELi}_{s_1,s_2,\dots,s_{2n}}(z) &=\mathrm{ELi}_{s_1-1,s_2,\dots,s_{2n}}(z),\\
z\,\frac{\mathrm{d}}{\mathrm{d}z}\,\mathrm{ELi}_{s_1,s_2,\dots,s_{2n+1}}(z) &=\mathrm{ELi}_{s_1-1,s_2,\dots,s_{2n+1}}(z).
\end{aligned}
\end{equation}

The solutions to the recursion relations can be written in terms of the following \emph{multiple harmonic numbers}:
\begin{equation}
H_m^{(s_1,\dots,s_n)}\equiv\sum_{j_1=1}^{m}\frac{1}{j_1^{s_1}}\sum_{j_2=1}^{j_1-1}\frac{1}{j_2^{s_2}}\dots\sum_{j_n=1}^{j_{n-1}-1}\frac{1}{j_n^{s_n}},\quad m,s_1,\dots,s_n\in\mathbb{Z}_{\ge 1}.
\end{equation}
We will also need the following two properties of the multiple harmonic numbers that follow from the definition. First, 
\begin{equation}\label{property1}
H_m^{(s_1,\dots,s_n)}=0\quad\text{if}\quad n>m,
\end{equation}
since at least one of the sums is empty.
The second property is
\begin{equation}\label{property2}
H_m^{(s_1,s_2,\dots,s_n)}=\frac{1}{m^{s_1}}\,H_{m-1}^{(s_2,\dots,s_n)}+H_{m-1}^{(s_1,s_2,\dots,s_n)},
\end{equation}
which follows by splitting the most external sum.

The recursion relation \eqref{coeffs>1} is solved by
\begin{equation}
c_j^{(s)}=\Gamma(j)\,H_{j-1}^{(\overbrace{1,\dots,1}^{s-1})}.
\end{equation}
The fact that $\Gamma(j)\,H_{j-1}^{(\overbrace{1,\dots,1}^{s-1})}=0$ for $j<s$ follows from \eqref{property1}. Let us prove by induction on $s> 1$ that for $j=s$ we have 
\begin{equation}
\Gamma(s)\,H_{s-1}^{(\overbrace{1,\dots,1}^{s-1})}=1.
\end{equation}
For $s=2$, this is trivial.
Suppose by inductive hypothesis $\Gamma(s-1)\,H_{s-2}^{(\overbrace{1,\dots,1}^{s-2})}=1$ for some $s>2$. Then, using \eqref{property1} and \eqref{property2},
\begin{equation}
\Gamma(s)\,H_{s-1}^{(\overbrace{1,\dots,1}^{s-1})}=\Gamma(s)\left(\frac{1}{s-1}\,H_{s-2}^{(\overbrace{1,\dots,1}^{s-2})}+H_{s-2}^{(\overbrace{1,\dots,1}^{s-1})}\right)=\Gamma(s)\left(\frac{1}{s-1}\,\frac{1}{\Gamma(s-1)}+0\right)=1.
\end{equation}
Finally, let us prove the claim for $j>s$. For $s=2$, we want to see that 
\begin{equation}
\Gamma(j)\,H_{j-1}^{(1)}=\Gamma(j-1)+(j-1)\,\Gamma(j-1)\,H_{j-2}^{(1)}.
\end{equation}
Dividing for $\Gamma(j)$, the above identity reads
\begin{equation}
\sum_{k=1}^{j-1}\frac{1}{k}=\frac{1}{j-1}+\sum_{k=1}^{j-2}\frac{1}{k},
\end{equation}
which holds for the property \eqref{property2}.

For $s>2$, let us see that
\begin{equation}
\Gamma(j)\,H_{j-1}^{(\overbrace{1,\dots,1}^{s-1})}=\Gamma(j-1)\,H_{j-2}^{(\overbrace{1,\dots,1}^{s-2})}+(j-1)\,\Gamma(j-1)\,H_{j-2}^{(\overbrace{1,\dots,1}^{s-1})}.
\end{equation}
Dividing for $\Gamma(j)$, the above identity reads
\begin{equation}
H_{j-1}^{(\overbrace{1,\dots,1}^{s-1})}=\frac{1}{j-1}\,H_{j-2}^{(\overbrace{1,\dots,1}^{s-2})}+H_{j-2}^{(\overbrace{1,\dots,1}^{s-1})},
\end{equation}
which holds due to \eqref{property2}.

For the cases with higher weights, the solutions to the recurrence relations are given by
\begin{equation}\label{solutionrecurrence12}
\begin{aligned}
c_j^{(s_1,\dots,s_{2n})}&=
\begin{cases}
0\quad &j\le\sum_{m=1}^{n}s_{2m}-1,\\
\frac{(-1)^{j-1+\sum_{m=1}^{n}s_{2m-1}}}{j^{s_1}}\,\Gamma(j)\,H_{j-1}^{(\overbrace{1,\dots,1}^{s_2-1},s_3+1,\dots,\overbrace{1,\dots,1}^{s_{2n}-1})}\quad &j>\sum_{m=1}^{n}s_{2m}-1,
\end{cases}
\end{aligned}
\end{equation}
\begin{equation}\label{solutionrecurrence34}
\begin{aligned}
c_j^{(s_1,\dots,s_{2n+1})}&=
\begin{cases}
0\quad &j\le\sum_{m=0}^{n}s_{2m+1}-1,\\
(-1)^{\sum_{m=1}^{n}s_{2m}}\,\Gamma(j)\,H_{j-1}^{(\overbrace{1,\dots,1}^{s_1-1},s_2+1,\dots,\overbrace{1,\dots,1}^{s_{2n+1}-1})}\quad &j>\sum_{m=0}^{n}s_{2m+1}-1,
\end{cases}
\end{aligned}
\end{equation}
with the convention that
\begin{equation}
\overbrace{1,\dots,1}^{0}=\emptyset.
\end{equation}
The cases in which the coefficients are zero follow from the property \eqref{property1}. The recurrence relations \eqref{coeffeven1} and \eqref{coeffs>1even} are satisfied by these solutions by straightforward substitution. For the recurrence relations \eqref{coeffodd1} and \eqref{coeffs>1odd}, instead, we use the property \eqref{property2}.



The solutions \eqref{solutionrecurrence12} and \eqref{solutionrecurrence34} satisfy the recurrence relation \eqref{coeffodd1} if and only if
\begin{equation}
\begin{aligned}
&(-1)^{\sum_{m=1}^{n}s_{2m}}\,\Gamma(j)\,H_{j-1}^{(s_2+1,\overbrace{1,\dots,1}^{s_{3}-1},\dots,\overbrace{1,\dots,1}^{s_{2n+1}-1})}=\\
&(-1)^j\,\frac{(-1)^{j-2+\sum_{m=1}^{n}s_{2m}}}{(j-1)^{s_2}}\,\Gamma(j-1)\,H_{j-2}^{(\overbrace{1,\dots,1}^{s_3-1},s_4+1,\dots,\overbrace{1,\dots,1}^{s_{2n+1}-1})}+(j-1)\,(-1)^{\sum_{m=1}^{n}s_{2m}}\,\Gamma(j-1)\,H_{j-2}^{(s_2+1,\overbrace{1,\dots,1}^{s_{3}-1},\dots,\overbrace{1,\dots,1}^{s_{2n+1}-1})}.
\end{aligned}
\end{equation}
Dividing both sides by $(-1)^{\sum_{m=1}^{n}s_{2m}}\,\Gamma(j)$, the equality becomes
\begin{equation}
\begin{aligned}
&H_{j-1}^{(s_2+1,\overbrace{1,\dots,1}^{s_{3}-1},\dots,\overbrace{1,\dots,1}^{s_{2n+1}-1})}=\frac{1}{(j-1)^{s_2+1}}\,H_{j-2}^{(\overbrace{1,\dots,1}^{s_3-1},s_4+1,\dots,\overbrace{1,\dots,1}^{s_{2n+1}-1})}+H_{j-2}^{(s_2+1,\overbrace{1,\dots,1}^{s_{3}-1},\dots,\overbrace{1,\dots,1}^{s_{2n+1}-1})},
\end{aligned}
\end{equation}
which is true because of \eqref{property2}.


The solutions \eqref{solutionrecurrence34} 
satisfy the recurrence relation \eqref{coeffs>1odd} if and only if
\begin{equation}
\begin{aligned}
&(-1)^{\sum_{m=1}^{n}s_{2m}}\,\Gamma(j)\,H_{j-1}^{(\overbrace{1,\dots,1}^{s_1-1},s_2+1,\dots,\overbrace{1,\dots,1}^{s_{2n+1}-1})}=\\
&(-1)^{\sum_{m=1}^{n}s_{2m}}\,\Gamma(j-1)\,H_{j-2}^{(\overbrace{1,\dots,1}^{s_1-2},s_2+1,\dots,\overbrace{1,\dots,1}^{s_{2n+1}-1})}+(j-1)\,(-1)^{\sum_{m=1}^{n}s_{2m}}\,\Gamma(j-1)\,H_{j-2}^{(\overbrace{1,\dots,1}^{s_1-1},s_2+1,\dots,\overbrace{1,\dots,1}^{s_{2n+1}-1})}.
\end{aligned}
\end{equation}
By dividing both sides by $(-1)^{\sum_{m=1}^{n}s_{2m}}\,\Gamma(j)$, the equality becomes
\begin{equation}
\begin{aligned}
&H_{j-1}^{(\overbrace{1,\dots,1}^{s_1-1},s_2+1,\dots,\overbrace{1,\dots,1}^{s_{2n+1}-1})}=&\frac{1}{j-1}\,H_{j-2}^{(\overbrace{1,\dots,1}^{s_1-2},s_2+1,\dots,\overbrace{1,\dots,1}^{s_{2n+1}-1})}+H_{j-2}^{(\overbrace{1,\dots,1}^{s_1-1},s_2+1,\dots,\overbrace{1,\dots,1}^{s_{2n+1}-1})},
\end{aligned}
\end{equation}
which is true by \eqref{property2}.

We conclude by inserting the above solutions of the recurrence relations in the Ansatz \eqref{asymptoticansatz}:
\begin{equation}
\boxed{
\begin{aligned}
&\mathrm{ELi}_{s_1,\dots,s_{2n+1}}(z)=\mathrm{e}^z\sum_{j=\sum_{m=0}^{n}s_{2m+1}}^N\frac{(-1)^{\sum_{m=1}^{n}s_{2m}}}{z^j}\,\Gamma(j)\,H_{j-1}^{(\overbrace{1,\dots,1}^{s_1-1},s_2+1,\dots,\overbrace{1,\dots,1}^{s_{2n+1}-1})}+\mathcal{O}\left(\frac{\mathrm{e}^z}{z^{N+1}}\right),\\
&\mathrm{ELi}_{s_1,\dots,s_{2n}}(z)=\sum_{j=\sum_{m=1}^{n}s_{2m}}^N\frac{(-1)^{j-1+\sum_{m=1}^{n}s_{2m-1}}}{j^{s_1}\,z^j}\,\Gamma(j)\,H_{j-1}^{(\overbrace{1,\dots,1}^{s_2-1},s_3+1,\dots,\overbrace{1,\dots,1}^{s_{2n}-1})}+\mathcal{O}\left(\frac{1}{z^{N+1}}\right).  
\end{aligned}}
\end{equation}

\section{Quadratic identities for \text{EL} functions}\label{sec:quadid}

In this section, we study quadratic relations satisfied by the dressed multiple polyexponential functions, which are inspired by the relations satisfied by the multiple polylogarithms.
The first relevant identity reads
\begin{equation}\label{EL1EL1}
\mathrm{EL}_1(z)\mathrm{EL}_1(-z)+\mathrm{EL}_{1,1}(z)+\mathrm{EL}_{1,1}(-z)+\mathrm{EL}_2(z)+\mathrm{EL}_2(-z)=0.
\end{equation}
This can be proven by considering the derivative and noticing that the left-hand side vanishes for $z=0$. We have
\begin{equation}
\begin{aligned}
&z\frac{\mathrm{d}}{\mathrm{d}z}\left[\mathrm{EL}_1(z)\mathrm{EL}_1(-z)+\mathrm{EL}_{1,1}(z)+\mathrm{EL}_{1,1}(-z)+\mathrm{EL}_2(z)+\mathrm{EL}_2(-z)\right]=\\
&=\left(\mathrm{e}^z-1\right)\,\mathrm{EL}_1(-z)+\left(\mathrm{e}^{-z}-1\right)\,\mathrm{EL}_1(z)-\mathrm{e}^z\,\mathrm{EL}_1(-z)-\mathrm{e}^{-z}\,\mathrm{EL}_1(z)+\mathrm{EL}_1(z)+\mathrm{EL}_1(-z)=0.
\end{aligned}
\end{equation}

Then, for every $n\ge 2$, the following identities hold: 
\begin{equation}\label{EL1ELn}
\begin{aligned}
\mathrm{EL}_1(z)\mathrm{EL}_n(-z)+(n+1)\mathrm{EL}_{n+1}(z)+\mathrm{EL}_{n,1}(z)+\sum_{j=1}^{n}\mathrm{EL}_{n+1-j,j}(z)+\left(z\mapsto -z\right)=0.
\end{aligned}
\end{equation}
We remark that this identity has a similar counterpart on the multiple polylogarithm side:
\begin{equation}
-\mathrm{Li}_1(z)\mathrm{Li}_n(z)+\mathrm{Li}_{n,1}(z)+\sum_{j=1}^{n}\mathrm{Li}_{n+1-j,j}(z)=0.
\end{equation}
We prove the identity \eqref{EL1ELn} by induction on $n\ge 2$. For the base step $n=2$, we have
\begin{equation}
\begin{aligned}
&z\frac{\mathrm{d}}{\mathrm{d}z}\left[\mathrm{EL}_1(z)\mathrm{EL}_2(-z)+\mathrm{EL}_2(z)\mathrm{EL}_1(-z)\right]=(\mathrm{e}^z-1)\mathrm{EL}_2(-z)+(\mathrm{e}^{-z}-1)\mathrm{EL}_2(z)+2\mathrm{EL}_1(z)\mathrm{EL}_1(-z)=\\
&=\mathrm{e}^z\mathrm{EL}_2(-z)+\mathrm{e}^{-z}\mathrm{EL}_2(z)-2\mathrm{EL}_{1,1}(z)-2\mathrm{EL}_{1,1}(-z)-3\mathrm{EL}_2(z)-3\mathrm{EL}_2(-z),
\end{aligned}
\end{equation}
where passing to the last line we used the identity \eqref{EL1EL1}. Since
\begin{equation}
z\frac{\mathrm{d}}{\mathrm{d}z}\mathrm{EL}_{1,2}(z)=-\mathrm{e}^z\mathrm{EL}_2(-z),
\end{equation}
we conclude that
\begin{equation}
\begin{aligned}
\mathrm{EL}_1(z)\mathrm{EL}_2(-z)+3\,\mathrm{EL}_{3}(z)+2\,\mathrm{EL}_{2,1}(z)+\mathrm{EL}_{1,2}(z)+\left(z\mapsto -z\right)=0,
\end{aligned}
\end{equation}
and the base step is proved.

For the inductive step, the inductive hypothesis reads
\begin{equation}
\begin{aligned}
\mathrm{EL}_1(z)\mathrm{EL}_{n-1}(-z)+n\,\mathrm{EL}_{n}(z)+\mathrm{EL}_{n-1,1}(z)+\sum_{j=1}^{n-1}\mathrm{EL}_{n-j,j}(z)+\left(z\mapsto -z\right)=0.
\end{aligned}
\end{equation}
Let us consider
\begin{equation}
\begin{aligned}
&z\frac{\mathrm{d}}{\mathrm{d}z}\left\{\mathrm{EL}_1(z)\mathrm{EL}_n(-z)+(n+1)\mathrm{EL}_{n+1}(z)+\mathrm{EL}_{n,1}(z)+\sum_{j=1}^{n}\mathrm{EL}_{n+1-j,j}(z)+\left(z\mapsto -z\right)\right\}=\\
&=(\mathrm{e}^z-1)\mathrm{EL}_n(-z)+\mathrm{EL}_1(z)\mathrm{EL}_{n-1}(-z)+\\
&\quad+(n+1)\mathrm{EL}_{n}(z)+\mathrm{EL}_{n-1,1}(z)+\sum_{j=1}^{n-1}\mathrm{EL}_{n-j,j}(z)-\mathrm{e}^z\,\mathrm{EL}_{n}(-z)+\left(z\mapsto -z\right)=\\
&=\mathrm{EL}_1(z)\mathrm{EL}_{n-1}(-z)+n\,\mathrm{EL}_{n}(z)+\mathrm{EL}_{n-1,1}(z)+\sum_{j=1}^{n-1}\mathrm{EL}_{n-j,j}(z)+\left(z\mapsto -z\right),
\end{aligned}
\end{equation}
which is equal to zero by the inductive hypothesis. Therefore, since the derivative of the non-linear relation vanishes and is equal to 0 at $z=0$, the induction proof is complete.

For $1\le m\le n$, we can write the generic identities using the following Ansatz\footnote{We remark that if $m=n$, the first part of the identity is given by $2\mathrm{EL}_m(z)\mathrm{EL}_m(-z)$.}:
\begin{equation}\label{genericnonlinearEL}
\begin{aligned}
\mathrm{EL}_m(z)\mathrm{EL}_n(-z)+\alpha_0^{(m,n)}\mathrm{EL}_{n+m}(z)+\sum_{j=1}^{n}\alpha_j^{(m,n)}\mathrm{EL}_{n+m-j,j}(z)+\left(z\mapsto -z\right)=0,
\end{aligned}
\end{equation}
for some integer coefficients $\alpha_j^{(m,n)}$, with $0\le j\le n$. Taking the derivative, we can find recursion relations satisfied by these coefficients.

Applying the operator $z\frac{\mathrm{d}}{\mathrm{d}z}$ to \eqref{genericnonlinearEL}, we get
\begin{equation}
\begin{aligned}
&\mathrm{EL}_{m}(z)\mathrm{EL}_{n-1}(-z)+\mathrm{EL}_{m-1}(z)\mathrm{EL}_n(-z)+\\
&+\alpha_0^{(m,n)}\mathrm{EL}_{n+m-1}(z)+\sum_{j=1}^{n}\alpha_j^{(m,n)}\mathrm{EL}_{n+m-j-1,j}(z)+\left(z\mapsto -z\right)=0.
\end{aligned}
\end{equation}
Using the same Ansatz \eqref{genericnonlinearEL} for the first terms
\begin{equation}
\begin{aligned}
&\mathrm{EL}_{m}(z)\mathrm{EL}_{n-1}(-z)+\alpha_0^{(m,n-1)}\mathrm{EL}_{n+m-1}(z)+\sum_{j=1}^{n-1}\alpha_j^{(m,n-1)}\mathrm{EL}_{n+m-j-1,j}(z)+\left(z\mapsto -z\right)=0,\\
&\mathrm{EL}_{m-1}(z)\mathrm{EL}_n(-z)+\alpha_0^{(m-1,n)}\mathrm{EL}_{n+m-1}(z)+\sum_{j=1}^{n}\alpha_j^{(m-1,n)}\mathrm{EL}_{n+m-j-1,j}(z)+\left(z\mapsto -z\right)=0,
\end{aligned}
\end{equation}
we conclude
\begin{equation}\label{recursionnonlinearEL}
\begin{aligned}
0\le j\le n-1\colon\quad &\alpha_j^{(m,n)}=\alpha_j^{(m,n-1)}+\alpha_j^{(m-1,n)},\\
j=n\colon\quad &\alpha_n^{(m,n)}=\alpha_n^{(m-1,n)}.
\end{aligned}
\end{equation}
From the second relation, since $\alpha_n^{(1,n)}=1$ for all $n\ge 2$, we conclude that 
\begin{equation}\label{consequencesecondrecursionrel}
\begin{aligned}
\alpha_n^{(m,n)}&=1\ \text{for\ all\ } m< n,\\
\alpha_m^{(m,m)}&=2\alpha_m^{(m-1,m)}=2.
\end{aligned}
\end{equation} 
The other recursions with $0\le j\le n-1$ must be solved by double induction on $(m,n)$ for all $n\ge m$.

The starting point of the induction is $m=1$, and, using the identity \eqref{EL1ELn}, we know that for all $n\ge 1$
\begin{equation}\label{alpha1n}
\begin{aligned}
\alpha_0^{(1,n)}=n+1,\quad \alpha_1^{(1,n)}=2,\quad \alpha_j^{(1,n)}=1\ \text{for}\ 2\le j\le n.
\end{aligned}
\end{equation}
Then, suppose that for some $m>1$ we know the values of $\alpha_j^{(m-1,n)}$ for all $n\ge m-1$. The solutions for $\alpha_j^{(m,n)}$ can be found by solving the recurrence relation
\begin{equation}
A[n]=A[n-1]+\alpha_j^{(m-1,n)}
\end{equation}
for $A[n]$ with initial condition
\begin{equation}
A[m]=2\,\alpha_j^{(m-1,m)},\quad\text{for\ all\ }j\le m,
\end{equation}
which follows from 
\begin{equation}\label{alphamm}
\alpha_j^{(m,m)}=2\,\alpha_j^{(m-1,m)},
\end{equation}
and with initial condition
\begin{equation}
A[j]=1,
\end{equation}
for $m<j\le n$, which follows from \eqref{consequencesecondrecursionrel}.

Let us show explicitly the solution for the case $m=2$. The starting point is
\begin{equation}\label{initialconditionsm=2}
\alpha_j^{(2,2)}=2\,\alpha_j^{(1,2)}=\begin{cases}
6\quad j=0,\\
4\quad j=1,\\
2\quad j=2.
\end{cases}
\end{equation}
Solving the recurrence relation \eqref{recursionnonlinearEL} with $m=2$, using the solutions \eqref{alpha1n} and with the initial conditions given in \eqref{initialconditionsm=2}, we find
\begin{equation}
\alpha_0^{(2,n)}=\binom{n+2}{2},\quad \alpha_1^{(2,n)}=2n,\quad \alpha_2^{(2,n)}=n,\quad 3\le j\le n\colon\ \alpha_j^{(2,n)}=n-j+1.
\end{equation}
For example, we have for $m=2$ and $n=2,3,4$:
\begin{equation}
\begin{aligned}
&\mathrm{EL}_2(z)\mathrm{EL}_2(-z)+6\,\mathrm{EL}_{4}(z)+4\,\mathrm{EL}_{3,1}(z)+2\,\mathrm{EL}_{2,2}(z)+\left(z\mapsto -z\right)=0,\\
&\mathrm{EL}_2(z)\mathrm{EL}_3(-z)+10\,\mathrm{EL}_{5}(z)+6\,\mathrm{EL}_{4,1}(z)+3\,\mathrm{EL}_{3,2}(z)+\mathrm{EL}_{2,3}(z)+\left(z\mapsto -z\right)=0,\\
&\mathrm{EL}_2(z)\mathrm{EL}_4(-z)+15\,\mathrm{EL}_{6}(z)+8\,\mathrm{EL}_{5,1}(z)+4\,\mathrm{EL}_{4,2}(z)+2\,\mathrm{EL}_{3,3}(z)+\mathrm{EL}_{2,4}(z)+\left(z\mapsto -z\right)=0.
\end{aligned}
\end{equation}


As happens for the case $j=n$ (see \eqref{consequencesecondrecursionrel}), we can also write the generic solution for $j=0,1,2$ and for all $m<j\le n$. These are given in terms of a single binomial coefficient:
\begin{align}
\alpha_0^{(m,n)}&=\binom{m+n}{m},\label{alphaj=0}\\
\alpha_1^{(m,n)}&=2\binom{m+n-2}{m-1},\label{alphaj=1}\\
\alpha_2^{(m,n)}&=\binom{m+n-2}{m-1},\label{alphaj=2}\\
\alpha_j^{(m,n)}&=\binom{n-j+m-1}{m-1},\quad\text{for\ } m<j\le n\label{alphaj>m},
\end{align}
which can be proved by double induction on $m\le n$. 
For the remaining cases $3\le j\le m\le n$, the $\alpha_j^{(m,n)}$ cannot be written in terms of a single binomial coefficient, as can be seen from the following case $j=m$:
\begin{equation}
\alpha_m^{(m,n)}=1+\binom{n-1}{m-1},\quad m\ge 3.\label{alphaj=m}
\end{equation}

Let us present the identities for $m=3,4$.
For $m=3$, we start with $n=3$:
\begin{equation}
\mathrm{EL}_3(z)\mathrm{EL}_3(-z)+20\mathrm{EL}_6(z)+12\mathrm{EL}_{5,1}(z)+6\mathrm{EL}_{4,2}(z)+2\mathrm{EL}_{3,3}(z)+(z\mapsto -z)=0,
\end{equation}
and the general identity for $n\geq3$ is given by
\begin{equation}
\begin{aligned}
&\alpha_1^{(3,n)}=n(n+1),\\
&\alpha_2^{(3,n)}=\binom{n+1}{2},\\
&\alpha_3^{(3,n)}=\frac{n^2-3 n+4}{2}=1+\binom{n-1}{2},\\
&\alpha_j^{(3,n)}=\binom{n-j+2}{2},\quad 4\le j\le n.
\end{aligned}
\end{equation}
For $m=4$, we start with $n=4$:
\begin{equation}
\mathrm{EL}_4(z)\mathrm{EL}_4(-z)+70\mathrm{EL}_8(z)+40\mathrm{EL}_{7,1}(z)+20\mathrm{EL}_{6,2}(z)+8\mathrm{EL}_{5,3}(z)+2\mathrm{EL}_{4,4}(z)+(z\mapsto -z)=0,
\end{equation}
and the general identity for $n\geq4$ is given by \eqref{alphaj=0}--\eqref{alphaj>m} and
\begin{equation}
\begin{aligned}
&\alpha_3^{(4,n)}=\frac{n \left(n^2-3 n+8\right)}{6}=2+\binom{n-2}{1}+\binom{n}{3},\\
&\alpha_4^{(4,n)}=\frac{n \left(n^2-6 n+11\right)}{6}=1+\binom{n-1}{3}.
\end{aligned}
\end{equation}

\section{Conclusion}
In this paper, we described a new set of special functions called multiple polyexponential integrals. These functions provide a complete perturbative description for the solutions of certain second-order ODEs around an irregular singularity of rank one. We also discuss two other sets of functions characterized by their Taylor series expansions around the origin.

The first application of multiple polyexponential integrals is to the linear perturbations of asymptotically flat black holes. In particular, they allow us to obtain analytic small-frequency expansions of the elements of the scattering matrix for the Schwarzschild black hole \cite{Aminov:2024mul}. The difference between the resulting small-frequency expansion and the other known methods (the MST and the instanton methods \cite{Mano:1996mf,Mano:1996vt,Mano:1996gn,Aminov:2020yma,Bautista:2023sdf}) is that there are no additional parameters. The MST method uses an additional non-integer parameter $\nu$ that regularizes the expansion. The instanton approach uses the cutoff $\Lambda\propto \omega$ as the perturbative parameter, while the frequency $\omega$ is considered to be finite. We already explored the connection between our small-frequency expansion and the instanton approach in \cite{Aminov:2024mul}. It would also be interesting to understand the connection with the MST method and the corresponding EFT developed recently in \cite{Ivanov:2022qqt,Saketh:2023bul,Ivanov:2024sds}.

From a mathematical perspective, our analysis raises some important related questions. In Sec.~\ref{sec:quadid} we studied a class of quadratic identities satisfied by the dressed multiple polyexponential functions. It would be interesting to generalize these identities to products of multiple polyexponential functions with level greater than one, which might lead to a more involved class of identities. In particular, some algebraic structures may underlie these relations as it happens for the harmonic algebra in the case of shuffle and stuffle products for multiple polylogarithms and multiple zeta values \cite{Wald,HOFFMAN1997477}.

Another relevant direction of study is the second-order ODEs with higher-rank irregular singularities. One possible example is the biconfluent Heun equation, which has an irregular singularity of rank two at infinity. The solution of this equation involves $\rme^{z^2}$ in a leading order when certain parameters are considered small. This poses the following question: what is the complete set of special functions describing the corresponding perturbative expansion? It is possible that multiple polyexponential integrals or their generalizations appear in this context too.

\appendix

\section{Derivative identity}\label{sec:App_diff}
In this appendix, we prove the identity
\begin{equation}\label{eq:derivative_id_el}
z\frac{\mathrm{d}}{\mathrm{d}z}el_{1,s_2,\dots,s_n}(z)=-el_{s_2,\dots,s_n}(z)-(-1)^n\mathrm{e}^z\sum_{\mathrm{op}(s_2)}\dots\sum_{\mathrm{op}(s_n)}el_{\mathrm{op}(s_2),\dots,\mathrm{op}(s_n)}(-z).
\end{equation}
First, using the definition (\ref{eq:polyExp_def}), we have
\begin{equation}\label{eq:lhs_derivative_id}
\begin{aligned}
&z\frac{\mathrm{d}}{\mathrm{d}z}el_{1,s_2,\dots,s_n}(z)+el_{s_2,\dots,s_n}(z)=\sum_{k_1>k_2>\dots>k_n\ge 1}\frac{z^{k_1}}{k_1! k_2^{s_2}\dots k_n^{s_n}}+\sum_{k_2>\dots>k_{n}\ge 1}\frac{z^{k_2}}{k_2! k_2^{s_2} k_3^{s_3}\dots k_n^{s_n}}=\\
&=\sum_{k_1\ge k_2>k_3>\dots>k_n\ge 1}\frac{z^{k_1}}{k_1! k_2^{s_2} k_3^{s_3}\dots k_n^{s_n}} = \sum_{k_1=1}^{\infty} H_{k_1}^{(s_2,\dots,s_n)} \, \frac{z^{k_1}}{k_1!},
\end{aligned}
\end{equation}
where the \emph{multiple harmonic numbers} were introduced:
\begin{equation}\label{multipleharmonicnumber}
H_m^{(s_2,\dots,s_n)}\equiv\sum_{m\ge k_2> k_3>\dots> k_n\ge 1}\frac{1}{k_2^{s_2}\dots k_n^{s_n}}.
\end{equation}
These numbers satisfy the following identity
\begin{equation}\label{fromidentity4.70}
\begin{aligned}
H_m^{(s_2,\dots,s_n)}=\sum_{m\ge k_1\ge k_2> k_3>\dots > k_n\ge 1}
\frac{(-1)^{k_1+k_2}}{k_1k_2^{s_2-1}k_3^{s_3}\dots k_n^{s_n}}\binom{m}{k_1}\binom{k_1}{k_2},
\end{aligned}
\end{equation}
which is instrumental to the proof of the derivative relation \eqref{eq:derivative_id_el}.
To prove \eqref{fromidentity4.70}, we define an additional set of \emph{multiple harmonic star numbers}:
\begin{equation}
{}^*H_m^{(s_2,\dots,s_n)}\equiv\sum_{m\ge k_2\ge k_3\ge\dots\ge k_n\ge 1}\frac{1}{k_2^{s_2}\dots k_n^{s_n}},
\end{equation}
where we sum over all $k_j$ with $j=2,\dots,n$.
These satisfy a very similar identity:
\begin{equation}\label{identity4.70}
{}^*H_m^{(s_2,\dots,s_n)}=\sum_{m\ge k_1\ge k_2\ge\dots\ge k_n\ge 1}\frac{(-1)^{k_1+k_2}}{k_1k_2^{s_2-1}k_3^{s_3}\dots k_n^{s_n}}\binom{m}{k_1}\binom{k_1}{k_2},
\end{equation}
which can be proven using the integral representation. In the simplest case when $n=2$, we have
\eq{{}^*H_m^{(s_2)}\equiv \sum_{k_2=1}^{m} \frac{1}{k_2^{s_2}}=H_m^{(s_2)},}
where $H_m^{(s_2)}$ is the $m$-th generalized harmonic number of order $s_2$. Taking $s_2=1$, we get $m$-th harmonic number with the following integral representation attributed to Euler:
\eq{H_m^{(1)}=\int_0^1\frac{1-x^m}{1-x}\,\rmd x.}
Substituting a new variable $x=1-u$, we get the following well-known identity:
\eq{H_m^{(1)}= \int_0^1\frac{1-\brc{1-u}^m}{u}\,\rmd u= \sum_{k_1=1}^{m} \frac{\brc{-1}^{k_1-1}}{k_1} \binom{m}{k_1}.}
Adding one more summation of the form
\eq{\sum_{k_2=1}^{k_1} \brc{-1}^{k_2} \binom{k_1}{k_2} =-1,}
gives us a particular case of (\ref{identity4.70}):
\eq{H_m^{(1)}= \sum_{k_1=1}^{m} \sum_{k_2=1}^{k_1} \frac{\brc{-1}^{k_1+k_2}}{k_1} \binom{m}{k_1} \binom{k_1}{k_2} .}
The same can be done in the slightly more general case with $s_2>1$. The integral representation of the corresponding generalized harmonic number is 
\eq{H_m^{(s_2)}= \int_{0}^{1}\frac{\rmd \,x_{s_2}}{x_{s_2}}\int_{0}^{x_{s_2}}\frac{\rmd \,x_{s_2-1}}{x_{s_2-1}} \dots \int_{0}^{x_{2}}\frac{1-x_{1}^m}{1-x_{1}} \rmd x_1.}
Defining $u=1-x_1$, we get for the innermost integral:
\eq{\int_{1-x_2}^{1}\frac{1-\brc{1-u}^m}{u} \rmd u = \sum_{k_1=1}^{m} \frac{\brc{-1}^{k_1-1}}{k_1} \binom{m}{k_1} \brc{1-\brc{1-x_2}^{k_1}}.}
The next integral is then
\eq{\int_{0}^{x_3}\brc{1-\brc{1-x_2}^{k_1}}\frac{\rmd x_2}{x_2} =
\int_{0}^{x_3}\sum_{k_2=1}^{k_1} \binom{k_1}{k_2} \brc{-x_2}^{k_2-1} \rmd x_2 =\sum_{k_2=1}^{k_1} \frac{\brc{-1}^{k_2-1}}{k_2} \binom{k_1}{k_2} x_3^{k_2}.}
Further integrals will increase the power of $k_2$ in the denominator of the latter expression until it becomes $s_2-1$. Putting everything together, we get the general identity for the case $n=2$:
\eqlb{eq:mH_id_level1}{H_m^{(s_2)}\equiv {}^*H_m^{(s_2)}= \sum_{m\ge k_1\ge k_2\ge 1} \frac{\brc{-1}^{k_1+k_2}}{k_1\,k_2^{s_2-1}} \binom{m}{k_1} \binom{k_1}{k_2} .}
Next, we consider the $n=3$ case, where the following definition is useful:
\begin{equation}
\begin{aligned}
\label{eq:mH_s2_def}
\mathcal{H}_m^{(s_2)} \brc{y}\equiv\,& \int_{y}^{1}\frac{\rmd \,x_{s_2}}{x_{s_2}}\int_{0}^{x_{s_2}}\frac{\rmd \,x_{s_2-1}}{x_{s_2-1}} \dots \int_{0}^{x_{2}}\frac{1-x_{1}^m}{1-x_{1}} \rmd x_1 \\
=\,&\sum_{m\ge k_1\ge k_2\ge 1} \frac{\brc{-1}^{k_1+k_2}}{k_1\,k_2^{s_2-1}} \binom{m}{k_1} \binom{k_1}{k_2} \brc{1-y^{k_2}},
\end{aligned}
\end{equation}
such that
\eq{H_m^{(s_2)}=\mathcal{H}_m^{(s_2)} \brc{0}.}
Now, the integral representation of ${}^*H_m^{(s_2,s_3)}$ can be written as
\eq{{}^*H_m^{(s_2,s_3)}=\int_{0}^{1}\frac{\rmd \,y_{s_3}}{y_{s_3}}\int_{0}^{y_{s_3}}\frac{\rmd \,y_{s_3-1}}{y_{s_3-1}} \dots \int_{0}^{y_{2}} \mathcal{H}_m^{(s_2)} \brc{y_1} \frac{\rmd y_1}{1-y_{1}} .}
The innermost integral can be computed  using the second line of (\ref{eq:mH_s2_def}): 
\eq{\int_{0}^{y_{2}} \frac{1- y_1^{k_2}}{1-y_1} \rmd y_1=\sum_{k_3=1}^{k_2} \frac{y_2^{k_3}}{k_3} .}
Taking further integrals increases the power of $k_3$ in the denominator, which gives us the next identity:
\eqlb{eq:mH_id_level2}{{}^*H_m^{(s_2,s_3)}= \sum_{m\ge k_1\ge k_2\ge k_3\ge 1} \frac{\brc{-1}^{k_1+k_2}}{k_1\,k_2^{s_2-1} k_3^{s_3}} \binom{m}{k_1} \binom{k_1}{k_2} .}
The corresponding level $2$ function $\mathcal{H}_m^{(s_2,s_3)} \brc{x}$ is defined as
\begin{equation}
\begin{aligned}\mathcal{H}_m^{(s_2,s_3)} \brc{x}\equiv & \int_{x}^{1}\frac{\rmd \,y_{s_3}}{y_{s_3}}\int_{0}^{y_{s_3}}\frac{\rmd \,y_{s_3-1}}{y_{s_3-1}} \dots \int_{0}^{y_{2}} \mathcal{H}_m^{(s_2)} \brc{y_1} \frac{\rmd y_1}{1-y_{1}}\\
=& \sum_{m\ge k_1\ge k_2\ge k_3\ge 1} \frac{\brc{-1}^{k_1+k_2}}{k_1\,k_2^{s_2-1} k_3^{s_3}} \binom{m}{k_1} \binom{k_1}{k_2} \brc{1-x^{k_3}}.
\end{aligned}
\end{equation}

The proof of the generic identity \eqref{identity4.70} proceeds inductively in the same way. The integral representation of the level $n-1$ generalized harmonic number ${}^*H_m^{(s_2,s_3,\dots,s_n)}$ is
\begin{equation}\label{recursiveidentityharmonic}
{}^*H_m^{(s_2,s_3,\dots,s_n)}=\int_0^1\frac{\rmd x_{s_n}}{x_{s_n}}\int_0^{x_{s_n}}\frac{\rmd x_{s_n-1}}{x_{s_n-1}}\dots\int_0^{x_{2}}\mathcal{H}_{m}^{(s_2,s_3,\dots,s_{n-1})}(x_1)\frac{\rmd x_{1}}{1-x_{1}},
\end{equation}
where the level $n-2$ function $\mathcal{H}_{m}^{(s_2,s_3,\dots,s_{n-1})}(x)$ is defined as
\begin{equation}
\begin{aligned}
\mathcal{H}_{m}^{(s_2,s_3,\dots,s_{n-1})}(x)&=\int_{x}^1\frac{\rmd y_{s_{n-1}}}{y_{s_{n-1}}}\int_0^{y_{s_{n-1}}} \frac{\rmd y_{s_{n-1}-1}}{y_{s_{n-1}-1}}\dots\int_0^{y_2} \mathcal{H}_{m}^{(s_2,s_3,\dots,s_{n-2})}\brc{y_1} \frac{\rmd y_1}{1-y_{1}}.
\end{aligned}
\end{equation}
By inductive hypothesis, $\mathcal{H}_{m}^{(s_2,s_3,\dots,s_{n-1})}(x)$ satisfies the following identity
\begin{equation}\label{curlyHgeneric}
\begin{aligned}
\mathcal{H}_{m}^{(s_2,s_3,\dots,s_{n-1})}(x)&=\sum_{m\ge k_1\ge k_2\ge\dots\ge k_{n-1}\ge 1}\frac{(-1)^{k_2+k_3}}{k_1\,k_2^{s_2-1}k_3^{s_3}\dots k_{n-1}^{s_{n-1}}}\binom{m}{k_1}\binom{k_1}{k_2}\brc{1-x^{k_{n-1}}}.
\end{aligned}
\end{equation}
Inserting \eqref{curlyHgeneric} in \eqref{recursiveidentityharmonic}, we can compute the innermost integral as
\begin{equation}
\int_0^{x_{2}}\frac{1-x_1^{k_{n-1}}}{1-x_{1}}\rmd x_{1}=\sum_{k_n=1}^{k_{n-1}}\frac{x_{2}^{k_n}}{k_n}.
\end{equation}
Hence, \eqref{recursiveidentityharmonic} becomes
\begin{equation}
\begin{aligned}
&{}^*H_m^{(s_2,s_3,\dots,s_n)}=\\
&\sum_{m\ge k_1\ge\dots\ge k_{n-1}\ge k_n\ge 1}\frac{(-1)^{k_2+k_3}}{k_1\,k_2^{s_2-1}k_3^{s_3}\dots k_{n-1}^{s_{n-1}}\,k_n}\binom{m}{k_1}\binom{k_1}{k_2}\int_0^1\frac{\rmd x_{s_n}}{x_{s_n}}\int_0^{x_{s_n}}\frac{\rmd x_{s_n-1}}{x_{s_n-1}}\dots\int_0^{x_{3}}x_{2}^{k_n-1}\rmd x_{2}.
\end{aligned}
\end{equation}
The last $s_n-1$ integrations increase the power of $k_n$, until it becomes $s_n$, concluding the proof.

The next step is to rewrite \eqref{identity4.70} in terms of the multiple harmonic numbers \eqref{multipleharmonicnumber}.
For level one multiple harmonic numbers $H_m^{(s_2)}$, the identity is the same (\ref{eq:mH_id_level1}). For level $2$, we have:
\eq{{}^*H_m^{(s_2,s_3)}=H_m^{(s_2,s_3)}+H_m^{(s_2+s_3)} .}
The same decomposition rule applies to the right-hand side of (\ref{eq:mH_id_level2}):
\begin{equation}
\begin{aligned}
{}^*H_m^{(s_2,s_3)}=\sum_{m\ge k_1\ge k_2\ge k_3\ge 1} \frac{\brc{-1}^{k_1+k_2}}{k_1\,k_2^{s_2-1} k_3^{s_3}} \binom{m}{k_1} \binom{k_1}{k_2} = & \sum_{m\ge k_1\ge k_2> k_3\ge 1} \frac{\brc{-1}^{k_1+k_2}}{k_1\,k_2^{s_2-1} k_3^{s_3}} \binom{m}{k_1} \binom{k_1}{k_2} \\
&+\sum_{m\ge k_1\ge k_2\ge 1} \frac{\brc{-1}^{k_1+k_2}}{k_1\,k_2^{s_2+s_3-1}}\binom{m}{k_1} \binom{k_1}{k_2}.
\end{aligned}
\end{equation}
Using (\ref{eq:mH_id_level1}), we can cancel the level one contribution and get the following result:
\eq{H_m^{(s_2,s_3)}=\sum_{m\ge k_1\ge k_2> k_3\ge 1} \frac{\brc{-1}^{k_1+k_2}}{k_1\,k_2^{s_2-1} k_3^{s_3}} \binom{m}{k_1} \binom{k_1}{k_2}. }
Continuing the same argument for multiple harmonic star numbers of level $n>2$, we cancel all multiple harmonic numbers with levels lower than $n$ to get
the general identity \eqref{fromidentity4.70}.
Applying \eqref{fromidentity4.70} to \eqref{eq:lhs_derivative_id}, we get
\begin{equation}\label{eq:identity_strict_harmonic}
\begin{aligned}
\sum_{k_1\ge k_2>k_3>\dots>k_n\ge 1}\frac{z^{k_1}}{k_1! k_2^{s_2} k_3^{s_3}\dots k_n^{s_n}}=&\sum_{m\ge 1}\frac{z^m}{m!}H_{m}^{(s_2,\dots,s_n)}\\
=&\sum_{m=1}^{\infty}\frac{z^{m}}{m!}\sum_{m\ge k_1\ge k_2> k_3>\dots > k_n\ge 1}
\frac{(-1)^{k_1+k_2}}{k_1k_2^{s_2-1}k_3^{s_3}\dots k_n^{s_n}}\binom{m}{k_1}\binom{k_1}{k_2}.
\end{aligned}
\end{equation}
By changing the order of summation between indices $m$ and $k_1$, we can factor out the exponential:
\begin{equation}\label{derivativestart}
\begin{aligned}
\sum_{k_1=1}^{\infty}\frac{\brc{-z}^{k_1}}{k_1! k_1}\sum_{m= k_1}^{\infty} \frac{z^{m-k_1}}{\brc{m-k_1}!} \sum_{k_2=1}^{k_1}\dots\sum_{k_n=1}^{k_{n-1}-1}
\frac{(-1)^{k_2}}{k_2^{s_2-1}k_3^{s_3}\dots k_n^{s_n}}\binom{k_1}{k_2}&\\
=\mathrm{e}^z\sum_{k_1=1}^{\infty}\frac{(-z)^{k_1}}{k_1! k_1}\sum_{k_2=1}^{k_1}\frac{(-1)^{k_2}}{k_2^{s_2-1}}\binom{k_1}{k_2}\sum_{k_2>k_3>\dots>k_n\ge 1}\frac{1}{k_3^{s_3}\dots k_n^{s_n}}.&
\end{aligned}
\end{equation}
Now, we are very close to proving the derivative identity. To see this, consider again the case with $n=2$:
\begin{equation}
z\frac{\mathrm{d}}{\mathrm{d}z}el_{1,s_2}(z)+el_{s_2}(z)=\rme^z\,\sum_{k_1=1}^{\infty}\frac{(-z)^{k_1}}{k_1! k_1}\sum_{k_2=1}^{k_1}\frac{(-1)^{k_2}}{k_2^{s_2-1}}\binom{k_1}{k_2}.
\end{equation}
The following identity will bring us to the desired result (\ref{eq:derivative_id_el}):
\eqlb{eq:sumId2_n2}{s_2>1:\quad \sum_{k_2=1}^{k_1}\frac{(-1)^{k_2}}{k_2^{s_2-1}}\binom{k_1}{k_2} =-
\sum_{k_2=1}^{k_1}\sum_{k_3=1}^{k_2}\dots \sum_{k_{s_2}=1}^{k_{s_2-1}}\frac{1}{k_2\,k_3\dots k_{s_2}}.}
For $s_2=1$, we have simply
\eqlb{eq:binomialcoeff_id}{\sum_{k_2=1}^{k_1}(-1)^{k_2}\binom{k_1}{k_2} = -1,}
which completes the proof of the first derivative identity:
\eq{z\frac{\mathrm{d}}{\mathrm{d}z}el_{1,1}(z)= - el_{1}(z) - \rme^z el_{1}\brc{-z}.}
For $s_2>1$, let's use the following integral representation:
\eq{\sum_{k_2=1}^{k_1}\frac{(-1)^{k_2}}{k_2^{s_2-1}}\binom{k_1}{k_2} = -\int_{0}^{1} \frac{\rmd y_{s_2-1}}{y_{s_2-1}}\int_{0}^{y_{s_2-1}} \frac{\rmd y_{s_2-2}}{y_{s_2-2}}\dots\int_{0}^{y_2} \frac{1-\brc{1-y_1}^{k_1}}{y_1} \,\rmd y_1.}
Changing every variable of integration  $y_i$ to $y_i=1-x_i$, we get
\eq{\sum_{k_2=1}^{k_1}\frac{(-1)^{k_2}}{k_2^{s_2-1}}\binom{k_1}{k_2} = - \int_{0}^{1} \frac{\rmd x_{s_2-1}}{1-x_{s_2-1}} \int_{x_{s_2-1}}^{1} \frac{\rmd x_{s_2-2}}{1-x_{s_2-2}}\dots \int_{x_2}^{1} \frac{1-x_1^{k_1}}{1-x_1} \,\rmd x_1.}
The innermost integral is
\eq{\int_{x_2}^{1} \frac{1-x_1^{k_1}}{1-x_1} \,\rmd x_1 = \sum_{k_2=1}^{k_1} \frac{1-x_2^{k_2}}{k_2},}
resulting in the next integral being of the same form with $k_1$ replaced by $k_2$:
\eq{ \sum_{k_2=1}^{k_1} \frac{1}{k_2} \int_{x_3}^{1} \frac{1-x_2^{k_2}}{1-x_2}\,\rmd x_2 =\sum_{k_2=1}^{k_1} \sum_{k_3=1}^{k_2} \frac{1-x_3^{k_3}}{k_2\, k_3}.}
Repeating this until we reach the last integral, gives us (\ref{eq:sumId2_n2}).
Going back to the derivative identity with $n=2$, we have
\begin{equation}
z\frac{\mathrm{d}}{\mathrm{d}z}el_{1,s_2}(z)=-el_{s_2}(z) -\rme^z\,\sum_{k_1\ge k_2\ge \dots\ge k_{s_2}\ge 1} \frac{1}{k_1\,k_2\dots k_{s_2}}\frac{(-z)^{k_1}}{k_1!},
\end{equation}
which can be rewritten using the sum over all ordered partitions of $s_2$:
\begin{equation}
z\frac{\mathrm{d}}{\mathrm{d}z}el_{1,s_2}(z)=-el_{s_2}(z) -\rme^z\,\sum_{\ops\brc{s_2}} el_{\ops\brc{s_2}}\brc{-z}.
\end{equation}
This completes the proof in the case of $n=2$. In the case $n=3$, the following identity is needed:
\eqlb{eq:sumId2_n3}{\sum_{k_2=1}^{k_1}\sum_{k_3=1}^{k_2-1}\frac{(-1)^{k_2}}{k_2^{s_2-1} k_3^{s_3}}\binom{k_1}{k_2} =\sum_{k_1\ge k_2\ge \dots\ge k_{s_2}>k_{s_2+1}\ge k_{s_2+2}\ge \dots\ge k_{s_2+s_3}\ge 1} \frac{1}{k_1\,k_2\dots k_{s_2}}\,\frac{1}{k_{s_2+1}\,k_{s_2+2}\dots k_{s_2+s_3}},}
where on the right-hand side, we do not sum over the first index $k_1$, and we have a strict inequality between $k_{s_2}$ and $k_{s_2+1}$. First, we introduce a level one function similar to $\mathcal{H}^{s_2}_m \brc{x}$:
\begin{equation}
\label{eq:def_G_l1}
\begin{aligned}
\mathcal{G}^{s_2}_{k_1} \brc{x}\equiv &  - \int_{x}^{1} \frac{\rmd y_{s_2-1}}{y_{s_2-1}}\int_{0}^{y_{s_2-1}} \frac{\rmd y_{s_2-2}}{y_{s_2-2}}\dots\int_{0}^{y_2} \frac{1-\brc{1-y_1}^{k_1}}{y_1} \,\rmd y_1\\
=& -
\sum_{k_2=1}^{k_1}\sum_{k_3=1}^{k_2}\dots \sum_{k_{s_2}=1}^{k_{s_2-1}}\frac{\brc{1-x}^{k_{s_2}}}{k_2\,k_3\dots k_{s_2}},
\end{aligned}
\end{equation}
such that
\eq{\sum_{k_2=1}^{k_1}\frac{(-1)^{k_2}}{k_2^{s_2-1}}\binom{k_1}{k_2} =\mathcal{G}^{s_2}_{k_1} \brc{0}.}
Then, the following integral representation can be written:
\eq{\sum_{k_2=1}^{k_1}\sum_{k_3=1}^{k_2}\frac{(-1)^{k_2}}{k_2^{s_2-1} k_3^{s_3}}\binom{k_1}{k_2} = \int_0^1\frac{\rmd y_{s_3}}{y_{s_3}}\int_{0}^{y_{s_3}}\frac{\rmd y_{s_3-1}}{y_{s_3-1}}\dots\int_{0}^{y_2}\mathcal{G}^{s_2}_{k_1}\brc{y_1} \frac{\rmd y_1}{1-y_1}.}
Again, we change all the variables of integration as $y_i=1-x_i$, which gives
\eq{\sum_{k_2=1}^{k_1}\sum_{k_3=1}^{k_2}\frac{(-1)^{k_2}}{k_2^{s_2-1} k_3^{s_3}}\binom{k_1}{k_2} = \int_0^1\frac{\rmd x_{s_3}}{1-x_{s_3}}\int_{x_{s_3}}^{1}\frac{\rmd x_{s_3-1}}{1-x_{s_3-1}}\dots\int_{x_2}^{1}\mathcal{G}^{s_2}_{k_1}\brc{1-x_1} \frac{\rmd x_1}{x_1}.}
Using the second line of (\ref{eq:def_G_l1}), we take the innermost integral:
\eq{\int_{x_2}^{1}\mathcal{G}^{s_2}_{k_1}\brc{1-x_1} \frac{\rmd x_1}{x_1} = -\sum_{k_2=1}^{k_1}\sum_{k_3=1}^{k_2}\dots \sum_{k_{s_2}=1}^{k_{s_2-1}}\frac{1-x_2^{k_{s_2}}}{k_2\,k_3\dots k_{s_2}^2}.}
Further integrals will add summations over new indices $k_{s_2+1}$, $\dots$, $k_{s_2+s_3-1}$, resulting in the following identity:
\eq{\sum_{k_2=1}^{k_1}\sum_{k_3=1}^{k_2}\frac{(-1)^{k_2}}{k_2^{s_2-1} k_3^{s_3}}\binom{k_1}{k_2} =- \sum_{k_1\ge k_2\ge \dots\ge k_{s_2+s_3-1}\ge 1} \frac{1}{k_1\,k_2\dots k_{s_2}^2}\,\frac{1}{k_{s_2+1}\,k_{s_2+2}\dots k_{s_2+s_3-1}}.}
Subtracting this from (\ref{eq:sumId2_n3}) gives us an already proven identity (\ref{eq:sumId2_n2}) with $s_2$ replaced by $s_2+s_3$, which, in turn, proves (\ref{eq:sumId2_n3}). To complete the proof of the derivative identity with $n=3$, we substitute (\ref{eq:sumId2_n3}) into (\ref{derivativestart}) and get via (\ref{eq:identity_strict_harmonic})
\begin{equation}
\begin{aligned}
z\frac{\mathrm{d}}{\mathrm{d}z} & el_{1,s_2,s_3}(z)+el_{s_2,s_3}(z)= \sum_{k_1\ge k_2>k_3>\dots>k_n\ge 1}\frac{z^{k_1}}{k_1! k_2^{s_2} k_3^{s_3}\dots k_n^{s_n}}\\
=&\,\rme^{z}\,\sum_{k_1\ge k_2\ge \dots\ge k_{s_2}>k_{s_2+1}\ge k_{s_2+2}\ge \dots\ge k_{s_2+s_3}\ge 1} \frac{1}{k_1\,k_2\dots k_{s_2}} \,\frac{1}{k_{s_2+1}\,k_{s_2+2}\dots k_{s_2+s_3}} \,\frac{(-z)^{k_1}}{k_1!},
\end{aligned}
\end{equation}
where we do sum over the first index $k_1$. The latter can be conveniently rewritten using ordered partitions $\ops\brc{s_2}$ and $\ops\brc{s_3}$:
\eq{z\frac{\mathrm{d}}{\mathrm{d}z} el_{1,s_2,s_3}(z)= - el_{s_2,s_3}(z) +\rme^z \sum_{\ops\brc{s_2}}\sum_{\ops\brc{s_3}}el_{\ops\brc{s_2},\ops\brc{s_3}}(-z).}

For general values of $n$, we want to prove that, for every $k_1\ge 1$, 
\begin{equation}\label{derivativedifficultpart}
\begin{aligned}
&\sum_{k_2=1}^{k_1}\frac{(-1)^{k_2}}{k_2^{s_2-1}}\binom{k_1}{k_2}\sum_{k_3=1}^{k_2-1}\frac{1}{k_3^{s_3}}\sum_{k_4=1}^{k_3-1}\frac{1}{k_4^{s_4}}\dots\sum_{k_n=1}^{k_{n-1}-1}\frac{1}{k_n^{s_n}}=\\
&-\sum_{l_{2,2}=1}^{k_1}\sum_{l_{2,3}=1}^{l_{2,2}}\dots\sum_{l_{2,s_2}=1}^{l_{2,s_2-1}}\frac{1}{l_{2,2}\,l_{2,3}\dots l_{2,s_2}}\sum_{k_3=1}^{l_{2,s_2}-1}\frac{(-1)^{k_3}}{k_3^{s_3}}\binom{l_{2,s_2}-1}{k_3}\sum_{k_4=1}^{k_3-1}\frac{1}{k_4^{s_4}}\dots\sum_{k_n=1}^{k_{n-1}-1}\frac{1}{k_n^{s_n}}.
\end{aligned}
\end{equation}
This can be done inductively using the integral representation as before. Or, equivalently, we can  use the \emph{Hockey-Stick identity}
\begin{equation}
\sum_{l=j}^n\binom{l-1}{j-1}=\binom{n}{j}, 
\end{equation}
which for any sequence $(a_j)_{j\in\mathbb{N}}$ and for any $k\ge 1$ gives
\begin{equation}
\begin{aligned}
\sum_{j=1}^n(-1)^{j-1}\binom{n}{j}\frac{a_j}{j^k}=\sum_{j=1}^n\sum_{l=j}^n(-1)^{j-1}\binom{l-1}{j-1}\frac{a_j}{j^k}=\sum_{j=1}^n\sum_{l=j}^n(-1)^{j-1}\binom{l}{j}\frac{1}{l}\frac{a_j}{j^{k-1}}.
\end{aligned}
\end{equation}
Changing the order of summation results in
\eqlb{stepreductionsums}{\sum_{j=1}^n(-1)^{j-1}\binom{n}{j}\frac{a_j}{j^k}=\sum_{l=1}^n\frac{1}{l}\sum_{m=1}^l(-1)^{m-1}\binom{l}{m}\frac{a_m}{m^{k-1}}.}
Iterating the argument, we have
\begin{equation}\label{iterativesequence}
\sum_{j=1}^n(-1)^{j-1}\binom{n}{j}\frac{a_j}{j^k}=\sum_{l_1=1}^n\frac{1}{l_1}\sum_{l_2=1}^{l_1}\frac{1}{l_2}\dots\sum_{l_k=1}^{l_{k-1}}\frac{1}{l_k}\sum_{m=1}^{l_k}(-1)^{m-1}\binom{l_k}{m}a_m.
\end{equation}
Let us apply this formula with $j=k_2, 
n=k_{1}, k=s_2-1$, and $(a_j)_{j\in\mathbb{N}}=\brc{H_{j-1}^{(s_3,\dots,s_n)}}_{j\in\mathbb{N}}$ 
to the left-hand side of \eqref{derivativedifficultpart}. We have
\begin{equation}
\begin{aligned}
&\sum_{k_2=1}^{k_{1}}\frac{(-1)^{k_2}}{k_2^{s_2-1}}\binom{k_{1}}{k_2}H_{k_2-1}^{(s_3,\dots,s_n)}=-\sum_{l_{2,2}=1}^{k_{1}}\frac{1}{l_{2,2}}\dots\sum_{l_{2,s_2}=1}^{l_{2,s_2-1}}\frac{1}{l_{2,s_2}}\sum_{m=1}^{l_{2,s_2}}(-1)^{m-1}\binom{l_{2,s_2}}{m}H_{m-1}^{(s_3,\dots,s_n)}.
\end{aligned}
\end{equation}
The identity \eqref{derivativedifficultpart} is proved if, for every $l,s_3,\dots,s_n\ge 1$, 
\begin{equation}\label{eq:mididentity}
\sum_{m=1}^{l}(-1)^{m-1}\binom{l}{m}H_{m-1}^{(s_3,\dots,s_n)}=\sum_{m=1}^{l}(-1)^{m-1}\binom{l}{m}\sum_{j=1}^{m-1}\frac{1}{j^{s_3}}H_{j-1}^{(s_4,\dots,s_n)}=\sum_{m=1}^{l-1}\frac{(-1)^m}{m^{s_3}}\binom{l-1}{m}H_{m-1}^{(s_4,\dots,s_n)}.
\end{equation}
We can prove that for every sequence $(b_j)_{j\in\mathbb{N}}$
\begin{equation}\label{usefulforbinomialtransform}
\sum_{m=1}^l(-1)^{m-1}\binom{l}{m}\sum_{j=1}^{m-1}b_j=\sum_{m=1}^{l-1}(-1)^m\binom{l-1}{m}b_m.
\end{equation}
By swapping the sums on the left-hand side, we have indeed
\begin{equation}
\sum_{m=1}^l(-1)^{m-1}\binom{l}{m}\sum_{j=1}^{m-1}b_j=\sum_{j=1}^{l-1}b_j\sum_{m=j+1}^{l}(-1)^{m-1}\binom{l}{m}=\sum_{j=1}^{l-1}(-1)^j\binom{l-1}{j}b_j,
\end{equation}
where in the last line we used the partial sum identity
\eq{\sum_{m=0}^j\brc{-1}^m\binom{l}{m}=\brc{-1}^j\binom{l-1}{j}.}
The proof of \eqref{eq:mididentity} then follows from \eqref{usefulforbinomialtransform} by considering the sequence
\begin{equation}
b_j=\frac{1}{j^{s_3}}H_{j-1}^{(s_4,\dots,s_n)}.
\end{equation}

Coming back to \eqref{derivativestart}, using \eqref{derivativedifficultpart}, we have proved that
\begin{equation}
\begin{aligned}
&\mathrm{e}^z\sum_{k_1=1}^{\infty}\frac{(-z)^{k_1}}{k_1!\,k_1}\sum_{k_2=1}^{k_1}\frac{(-1)^{k_2}}{k_2^{s_2-1}}\binom{k_1}{k_2}\sum_{k_3=1}^{k_2-1}\dots\sum_{k_n=1}^{k_{n-1}-1}\frac{1}{k_3^{s_3}\dots k_n^{s_n}}=\\
-&\mathrm{e}^z\sum_{k_1=1}^{\infty}\frac{(-z)^{k_1}}{k_1!}\sum_{l_{2,2}=1}^{k_1}\dots\sum_{l_{2,s_2}=1}^{l_{2,s_2-1}}\frac{1}{k_1\, l_{2,2}\dots l_{2,s_2}}\sum_{k_3=1}^{l_{2,s_2}-1}\frac{(-1)^{k_3}}{k_3^{s_3}}\binom{l_{2,s_2}-1}{k_3}\sum_{k_4=1}^{k_3-1}\dots\sum_{k_n=1}^{k_{n-1}-1}\frac{1}{k_4^{s_4}\dots k_n^{s_n}}.
\end{aligned}
\end{equation}
Applying recursively \eqref{derivativedifficultpart}, we turn every summation over $k_j$ with $j=3,\dots,n$ into the summation over the set of indices $l_{j,1}$, $\dots$, $l_{j,s_j}$. The result can be written as
\begin{equation}
\begin{aligned}
&\mathrm{e}^z\sum_{k_1=1}^{\infty}\sum_{k_2=1}^{k_1}\frac{(-z)^{k_1}}{k_1!\, k_1}\frac{(-1)^{k_2}}{k_2^{s_2-1}}\binom{k_1}{k_2}\sum_{k_2>k_3>\dots>k_n\ge 1}\frac{1}{k_3^{s_3}\dots k_n^{s_n}}=\\
&=(-1)^{n-1}\mathrm{e}^z\sum_{k_1=1}^{\infty}\frac{(-z)^{k_1}}{k_1!}\sum_{k_1\ge l_{2,2}\ge\dots\ge l_{2,s_2}\ge 1}\frac{1}{k_1\, l_{2,2}\dots l_{2,s_2}}\times\\
&\times\sum_{l_{2,s_2}>l_{3,1}\ge l_{3,2}\ge\dots\ge l_{3,s_3}\ge 1}\frac{1}{l_{3,1} l_{3,2}\dots l_{3,s_3}}\dots\sum_{k_n=1}^{l_{n-1,s_{n-1}}-1}\frac{(-1)^{k_n}}{k_n^{s_n}}\binom{l_{n-1,s_{n-1}}-1}{k_n}=\\
&=(-1)^{n-1}\mathrm{e}^z\sum_{k_1=1}^{\infty}\frac{(-z)^{k_1}}{k_1!}\sum_{k_1\ge l_{2,2}\ge\dots\ge l_{2,s_2}\ge 1}\frac{1}{k_1\, l_{2,2}\dots l_{2,s_2}}\times\\
&\times\sum_{l_{2,s_2}>l_{3,1}\ge l_{3,2}\ge\dots\ge l_{3,s_3}\ge 1}\frac{1}{l_{3,1} l_{3,2}\dots l_{3,s_3}}\dots\sum_{l_{n-1,s_{n-1}}>l_{n,1}\ge l_{n,2}\ge\dots\ge l_{n,s_n}\ge 1}\frac{1}{l_{n,1}\dots l_{n,s_n}},
\end{aligned}
\end{equation}
where passing to the last line we used the previous argument \eqref{iterativesequence} with $(a_j)_{j\in\mathbb{N}}=(1)_{j\in\mathbb{N}}$, and the identity \eqref{eq:binomialcoeff_id}.
Therefore, the proof is finished noting that, for every $j=2,\dots,n$, the sums
\begin{equation}
\sum_{l_{j,1}\ge l_{j,2}\ge\dots\ge l_{j,s_j}\ge 1}\frac{1}{l_{j,1}\,l_{j,2}\dots l_{j,s_j}}
\end{equation}
with $l_{2,1}\equiv k_1$
produce the sums over the ordered partitions of $s_j$:
\begin{equation}
\begin{aligned}
&(-1)^{n-1}\mathrm{e}^z\sum_{k_1=1}^{\infty}\frac{(-z)^{k_1}}{k_1!}\sum_{k_1\ge l_{2,2}\ge\dots\ge l_{2,s_2}\ge 1}\frac{1}{k_1\, l_{2,2}\dots l_{2,s_2}}\times\\
&\sum_{l_{2,s_2}>l_{3,1}\ge l_{3,2}\ge\dots\ge l_{3,s_3}\ge 1}\frac{1}{l_{3,1} l_{3,2}\dots l_{3,s_3}}\dots\sum_{l_{n-1,s_{n-1}}>l_{n,1}\ge l_{n,2}\ge\dots\ge l_{n,s_n}\ge 1}\frac{1}{l_{n,1}\dots l_{n,s_n}}=\\
&(-1)^{n-1}\mathrm{e}^z\sum_{\mathrm{op}(s_2)}\sum_{\mathrm{op}(s_3)}\dots\sum_{\mathrm{op}(s_n)}el_{\mathrm{op}(s_2),\mathrm{op}(s_3),\dots,\mathrm{op}(s_n)}(-z).
\end{aligned}
\end{equation}

\section{Relation between \texorpdfstring{$\mathrm{EL}$}{} and \texorpdfstring{$el$}{} functions by induction}
\label{app:rel_Taylor}

In this appendix, we want to prove that, for every $k\ge 1$,
\begin{equation}\label{ELeven}
\begin{aligned}
\mathrm{EL}_{r_1,s_1,\dots,r_k,s_k}(z)=\sum_{\mathrm{op}(s_1+1)}\dots\sum_{\mathrm{op}(s_k+1)}el_{(r_1-1)\oplus\mathrm{op}(s_1+1)\oplus\dots\oplus(r_k-1)\oplus\mathrm{op}(s_k+1)}(z),
\end{aligned}
\end{equation}
and
\begin{equation}\label{ELodd}
\begin{aligned}
\mathrm{EL}_{r_1,s_1,\dots,r_{k},s_{k},r_{k+1}}(z)=\sum_{\mathrm{op}(s_1+1)}\dots\sum_{\mathrm{op}(s_{k}+1)}el_{(r_1-1)\oplus\mathrm{op}(s_1+1)\oplus\dots\oplus(r_{k}-1)\oplus\mathrm{op}(s_{k}+1)\oplus r_{k+1}}(z).
\end{aligned}
\end{equation}
The proof is done by induction on $k$. The base step, $k=1$, has already been proved in Sec.~\ref{sec:mpf_def} and Sec.~\ref{sec:taylor}.

For the inductive step, we begin from the case $r_1=1$. Once the proof is complete for $r_1=1$, the proof for higher $r_1$ can be obtained by integrating and increasing the first weight on both sides of the equality.

Supposing the identities \eqref{ELeven} and \eqref{ELodd} to hold for every $1\le k\le n-1$, we want to prove
\begin{equation}\label{firstidentitygenericEL}
\begin{aligned}
\mathrm{EL}_{1,s_1,\dots,r_{n},s_{n}}(z)=\sum_{\mathrm{op}(s_1+1)}\dots\sum_{\mathrm{op}(s_{n}+1)}el_{\mathrm{op}(s_1+1)\oplus\dots\oplus(r_{n}-1)\oplus\mathrm{op}(s_{n}+1)}(z),
\end{aligned}
\end{equation}
and
\begin{equation}\label{secondidentitygenericEL}
\begin{aligned}
\mathrm{EL}_{1,s_1,\dots,r_{n},s_{n},r_{n+1}}(z)=\sum_{\mathrm{op}(s_1+1)}\dots\sum_{\mathrm{op}(s_n+1)}el_{\mathrm{op}(s_1+1)\oplus\dots\oplus(r_n-1)\oplus\mathrm{op}(s_n+1)\oplus r_{n+1}}(z).
\end{aligned}
\end{equation}
We will prove \eqref{firstidentitygenericEL} by induction on $s_1\ge 1$. The same steps we take to prove \eqref{firstidentitygenericEL} can be used to get \eqref{secondidentitygenericEL}. In particular, the key identity to obtain \eqref{firstidentitygenericEL} will be \eqref{identity_harmonic_appendix}, while the analogous for \eqref{secondidentitygenericEL} would be the identity \eqref{keyidentityharmonic2}. We will see at the end of the appendix that these identities can be proved simultaneously by induction.

For the base step on $s_1$, let us prove
\begin{equation}\label{id:EL_11n}
\mathrm{EL}_{1,1,r_2,s_2\dots,r_{n},s_{n}}(z)=\sum_{\mathrm{op}(2)}\sum_{\mathrm{op}(s_{2}+1)}\dots\sum_{\mathrm{op}(s_{n}+1)}el_{\mathrm{op}(2)\oplus(r_2-1)\oplus\dots\oplus(r_{n}-1)\oplus\mathrm{op}(s_{n}+1)}(z).
\end{equation}
The derivative of the right-hand side is given by the following identity:
\begin{equation}\label{basestepELidentity}
\begin{aligned}
&z\frac{\mathrm{d}}{\mathrm{d}z}\sum_{\mathrm{op}(2)}\sum_{\mathrm{op}(s_{2}+1)}\dots\sum_{\mathrm{op}(s_{n}+1)}el_{\mathrm{op}(2)\oplus(r_2-1)\oplus \mathrm{op}(s_{2}+1)\oplus\dots\oplus(r_{n}-1)\oplus\mathrm{op}(s_{n}+1)}(z)=\\
&-\mathrm{e}^z\sum_{\text{op}(r_2+1)}\dots\sum_{\text{op}(r_n+1)}el_{\text{op}(r_2+1)\oplus\dots\oplus\text{op}(r_n+1)\oplus s_n}(-z).
\end{aligned}
\end{equation}
To prove it, we use the Taylor series of the undressed multiple polyexponential functions on the left-hand side of the equality to get
\begin{equation}
\begin{aligned}
&z\frac{\mathrm{d}}{\mathrm{d}z}\sum_{\mathrm{op}(2)}\sum_{\mathrm{op}(s_{2}+1)}\dots\sum_{\mathrm{op}(s_{n}+1)}el_{\mathrm{op}(2)\oplus(r_2-1)\oplus \mathrm{op}(s_{2}+1)\oplus\dots\oplus(r_{n}-1)\oplus\mathrm{op}(s_{n}+1)}(z)=\\
&\sum_{k=1}^{\infty}\frac{z^k}{k!}\,{}^*H_k^{(r_2+1,\overbrace{1,\dots,1}^{s_2-1},r_3+1,\dots,\overbrace{1,\dots,1}^{s_n})},
\end{aligned}
\end{equation}
For the right-hand side of (\ref{basestepELidentity}), we have
\begin{equation}
\begin{aligned}
&-\mathrm{e}^z\sum_{\text{op}(r_2+1)}\dots\sum_{\text{op}(r_n+1)}el_{\text{op}(r_2+1)\oplus\dots\oplus\text{op}(r_n+1)\oplus s_n}(-z)=\\
&\sum_{m\ge 0}\frac{z^m}{m!}\sum_{l=1}^{\infty}\frac{(-z)^{l-1}}{l!\,l}\,{}^*H_l^{(\overbrace{1,\dots,1}^{r_2-1},s_2+1,\dots,\overbrace{1,\dots,1}^{r_n-1},s_n+1)}=\\
&\sum_{k\ge 1}\frac{z^k}{k!}\sum_{l= 1}^k\frac{(-1)^{l-1}}{l}\binom{k}{l}\,{}^*H_l^{(\overbrace{1,\dots,1}^{r_2-1},s_2+1,\dots,\overbrace{1,\dots,1}^{r_n-1},s_n+1)}.
\end{aligned}
\end{equation}
Then, \eqref{basestepELidentity} is a consequence of the general identity
\begin{equation}\label{basestepharmonicidentity}
k\ge 1:\quad {}^*H_k^{(r_2+1,\overbrace{1,\dots,1}^{s_2-1},\dots,\overbrace{1,\dots,1}^{s_n})}=\sum_{l= 1}^k\frac{(-1)^{l-1}}{l}\binom{k}{l}\,{}^*H_l^{(\overbrace{1,\dots,1}^{r_2-1},s_2+1,\dots,\overbrace{1,\dots,1}^{r_n-1},s_n+1)},
\end{equation}
which we prove at the end of the Appendix.

Assuming that \eqref{basestepELidentity} is proven, we rewrite the right-hand side using the inductive hypothesis:
\begin{equation}
\begin{aligned}
z\frac{\mathrm{d}}{\mathrm{d}z}\sum_{\mathrm{op}(2)}\sum_{\mathrm{op}(s_{2}+1)}\dots\sum_{\mathrm{op}(s_{n}+1)} & el_{\mathrm{op}(2)\oplus(r_2-1)\oplus \mathrm{op}(s_{2}+1)\oplus\dots\oplus(r_{n}-1)\oplus\mathrm{op}(s_{n}+1)}(z)\\
&=-\mathrm{e}^z\,\text{EL}_{1,r_2,s_2\dots,r_{n},s_{n}}(-z)=z\frac{\mathrm{d}}{\mathrm{d}z}\text{EL}_{1,1,r_2,s_2\dots,r_{n},s_{n}}(z).
\end{aligned}
\end{equation}
This proves that \eqref{id:EL_11n} is true up to an integration constant. By definition of the dressed multiple polyexponential functions, this constant is equal to zero, which concludes the proof of the base step on $s_1$.

Suppose now that \eqref{firstidentitygenericEL} holds for every $1\le s_1\le N-1$, and let's prove the thesis for $s_1=N$. Consider the derivative of the right-hand side of \eqref{firstidentitygenericEL}:
\begin{equation}\label{inductivestepidentitygenericEL}
\begin{aligned}
&z\frac{\mathrm{d}}{\mathrm{d}z}\sum_{\mathrm{op}(N+1)}\dots\sum_{\mathrm{op}(s_{n}+1)}el_{\mathrm{op}(N+1)\oplus(r_2-1)\oplus\dots\oplus(r_{n}-1)\oplus\mathrm{op}(s_{n}+1)}(z)=\\
&z\frac{\mathrm{d}}{\mathrm{d}z}\sum_{\mathrm{op}(N)}\dots\sum_{\mathrm{op}(s_{n}+1)}\left[el_{\mathrm{op}(N), r_2\oplus\dots\oplus(r_{n}-1)\oplus\mathrm{op}(s_{n}+1)}(z)+el_{\mathrm{op}(N)\oplus r_2\oplus\dots\oplus(r_{n}-1)\oplus\mathrm{op}(s_{n}+1)}(z)\right],
\end{aligned}
\end{equation}
where the sum over $\ops\brc{N+1}$ was split into two sums over $\ops\brc{N}$ using \eqref{eq:op_dec2}. By inductive hypothesis on $s_1$, the derivative of the second term reads
\begin{equation}
\begin{aligned}
&z\frac{\mathrm{d}}{\mathrm{d}z}\sum_{\mathrm{op}(N)}\dots\sum_{\mathrm{op}(s_{n}+1)}el_{\mathrm{op}(N)\oplus r_2\oplus\dots\oplus(r_{n}-1)\oplus\mathrm{op}(s_{n}+1)}(z)=z\frac{\mathrm{d}}{\mathrm{d}z}\,\text{EL}_{1,N-1,r_2+1,s_2,\dots,r_n,s_n}(z)=\\
&-\mathrm{e}^z\,\text{EL}_{N-1,r_2+1,s_2,\dots,r_n,s_n}(-z).
\end{aligned}
\end{equation}
By inductive hypothesis on the level of the dressed multiple polyexponential function, we then have
\begin{equation}
\begin{aligned}
-\mathrm{e}^z\,\text{EL}_{N-1,r_2+1,s_2,\dots,r_n,s_n}(-z)=-\mathrm{e}^z\sum_{\text{op}(r_2+2)}\dots\sum_{\text{op}(r_n+1)}el_{(N-2)\oplus\text{op}(r_2+2)\oplus\dots\oplus\text{op}(r_n+1)\oplus s_n}(-z).
\end{aligned}
\end{equation}
For the derivative of the first sum on the right-hand side of \eqref{inductivestepidentitygenericEL}, we use 
the following identity analogous to \eqref{eq:obs_el1n} and \eqref{eq:threeblockderivativeel}:
\begin{equation}\label{eq:nblockderivativeel}
\begin{aligned}
&z\frac{\mathrm{d}}{\mathrm{d}z}\sum_{\mathrm{op}(N)}\sum_{\text{op}(s_2+1)}\dots\sum_{\text{op}(s_n+1)}el_{\mathrm{op}(N),r_2\oplus\text{op}(s_2+1)\oplus\dots\oplus(r_n-1)\oplus\text{op}(s_n+1)}(z)=\\
&\mathrm{e}^z\sum_{\text{op}(r_2+1)}\dots\sum_{\text{op}(r_n+1)}el_{N-1,\text{op}(r_2+1)\oplus\dots\oplus\text{op}(r_n+1)\oplus s_n}(-z).
\end{aligned}
\end{equation}
For the proof of \eqref{eq:nblockderivativeel}, let us look at the Taylor series expansion of the left-hand side
\begin{equation}
\begin{aligned}
&z\frac{\mathrm{d}}{\mathrm{d}z}\sum_{\mathrm{op}(N)}\sum_{\text{op}(s_2+1)}\dots\sum_{\text{op}(s_n+1)}el_{\mathrm{op}(N),r_2\oplus\text{op}(s_2+1)\oplus\dots\oplus(r_n-1)\oplus\text{op}(s_n+1)}(z)=\\
&\sum_{k_1\ge k_2\ge\dots\ge k_N\ge 1}\frac{z^{k_1}}{k_1!\,k_2\dots k_N}\,{}^*H_{k_N-1}^{(r_2+1,\overbrace{1,\dots,1}^{s_2-1},\dots,r_n+1,\overbrace{1,\dots,1}^{s_n})},
\end{aligned}
\end{equation}
and of the right-hand side
\begin{equation}
\begin{aligned}
&\mathrm{e}^z\sum_{\text{op}(r_2+1)}\dots\sum_{\text{op}(r_n+1)}el_{N-1,\text{op}(r_2+1)\oplus\dots\oplus\text{op}(r_n+1)\oplus s_n}(-z)=\\
&\sum_{m\ge 0}\frac{z^m}{m!}\sum_{l\ge 1}\frac{(-z)^l}{l!\,l^{N-1}}\,{}^*H_{l-1}^{(\overbrace{1,\dots,1}^{r_2},s_2+1,\dots,\overbrace{1,\dots,1}^{r_n-1},s_n+1)}=\\
&\sum_{k_1=1}^{\infty}\frac{z^{k_1}}{k_1!}\sum_{l=1}^{k_1}\frac{(-1)^{l}}{l^{N-1}}\binom{k_1}{l}\,{}^*H_{l-1}^{(\overbrace{1,\dots,1}^{r_2},s_2+1,\dots,\overbrace{1,\dots,1}^{r_n-1},s_n+1)}=\\
&\sum_{k_1=1}^{\infty}\frac{z^{k_1}}{k_1!}\sum_{k_1\ge k_2\ge\dots\ge k_N}\frac{1}{k_2\dots k_N}\sum_{j=1}^{k_N}(-1)^j\binom{k_N}{j}\,{}^*H_{j-1}^{(\overbrace{1,\dots,1}^{r_2},s_2+1,\dots,\overbrace{1,\dots,1}^{r_n-1},s_n+1)},
\end{aligned}
\end{equation}
where to pass to the last line we used the identity \eqref{iterativesequence}.
The identity is then a consequence of the general identity
\begin{equation}\label{identity_harmonic_appendix}
k\ge 2:\quad {}^*H_{k-1}^{(r_2+1,\overbrace{1,\dots,1}^{s_2-1},\dots,\overbrace{1,\dots,1}^{s_n})}=\sum_{j=1}^k(-1)^{j}\binom{k}{j}\,{}^*H_{j-1}^{(\overbrace{1,\dots,1}^{r_2},s_2+1,\dots,\overbrace{1,\dots,1}^{r_n-1},s_n+1)},
\end{equation}
the proof of which we postpone to the end of the Appendix.

Putting all the results together, we get from \eqref{inductivestepidentitygenericEL}
\begin{equation}
\begin{aligned}
&z\frac{\mathrm{d}}{\mathrm{d}z}\sum_{\mathrm{op}(N+1)}\dots\sum_{\mathrm{op}(s_{n}+1)}el_{\mathrm{op}(N+1)\oplus (r_2-1)\oplus\dots\oplus(r_{n}-1)\oplus\mathrm{op}(s_{n}+1)}(z)=\\
&-\mathrm{e}^z\sum_{\text{op}(r_2+2)}\dots\sum_{\text{op}(r_n+1)}el_{(N-2)\oplus\text{op}(r_2+2)\oplus\dots\oplus\text{op}(r_n+1)\oplus s_n}(-z)\\
&+\mathrm{e}^z\sum_{\mathrm{op}(r_2+1)}\dots\sum_{\mathrm{op}(r_n+1)}el_{N-1,\mathrm{op}(r_2+1)\oplus(s_2-1)\oplus\dots\oplus\mathrm{op}(r_n+1)\oplus s_n}(-z)=\\
&-\mathrm{e}^z\sum_{\text{op}(r_2+1)}\dots\sum_{\text{op}(r_n+1)}el_{(N-1)\oplus\text{op}(r_2+1)\oplus\dots\oplus\text{op}(r_n+1)\oplus s_n}(-z)=-\mathrm{e}^z\,\text{EL}_{N,r_2,\dots,r_n,s_n}(-z).
\end{aligned}
\end{equation}
Since
\begin{equation}
z\frac{\mathrm{d}}{\mathrm{d}z}\text{EL}_{1,N,r_2,\dots,r_n,s_n}(z)=-\mathrm{e}^z\,\text{EL}_{N,r_2,\dots,r_n,s_n}(-z),
\end{equation}
\eqref{firstidentitygenericEL} is proven.

We conclude the appendix by proving the identities \eqref{basestepharmonicidentity} and \eqref{identity_harmonic_appendix} in the reverse order. 
By taking the binomial transform, the identity \eqref{identity_harmonic_appendix} is equivalent to the identities of the form
\begin{equation}\label{keyidentityharmonic1}
\begin{aligned}
{}^*H_{N-1}^{(\overbrace{1,\dots,1}^{r_2-1},s_2+1,\dots,\overbrace{1,\dots,1}^{r_n-1},s_n+1)}=\sum_{j=1}^N(-1)^j\binom{N}{j}\,{}^*H_{j-1}^{(r_2,\overbrace{1,\dots,1}^{s_2-1},r_3+1,\dots,r_n+1,\overbrace{1,\dots,1}^{s_n})}
\end{aligned}
\end{equation}
or
\begin{equation}\label{keyidentityharmonic2}
\begin{aligned}
{}^*H_{N-1}^{(\overbrace{1,\dots,1}^{r_2-1},s_2+1,\dots,\overbrace{1,\dots,1}^{r_n-1},s_n+1,\overbrace{1,\dots,1}^{r_{n+1}})}=\sum_{j=1}^N(-1)^j\binom{N}{j}\,{}^*H_{j-1}^{(r_2,\overbrace{1,\dots,1}^{s_2-1},r_3+1,\dots,r_n+1,\overbrace{1,\dots,1}^{s_n-1},r_{n+1}+1)}.
\end{aligned}
\end{equation}
The two identities \eqref{keyidentityharmonic1} and \eqref{keyidentityharmonic2} are analogous and can be proved by induction on the number of (arrays of) indices of the form $\overbrace{(1,\dots,1)}^s$ or $r>1$. 

The base step is 
\begin{equation}\label{basestepharmonicbinomial}
{}^*H_{N-1}^{\overbrace{(1,\dots,1)}^{n}}=\sum_{j=1}^{N}(-1)^j\binom{N}{j}\,{}^*H_{j-1}^{(n)}.
\end{equation}
We can manipulate the right-hand side as follows:
\begin{equation}
\begin{aligned}
&\sum_{j=1}^{N}(-1)^j\binom{N}{j}\,{}^*H_{j-1}^{(n)}=\sum_{j=1}^{N}\sum_{l=1}^{j-1}(-1)^j\binom{N}{j}\frac{1}{l^{n}}=\\
&=\sum_{l=1}^{N-1}\frac{1}{l^{n}}\sum_{j=l+1}^{N}(-1)^j\binom{N}{j}=\sum_{l=1}^{N-1}(-1)^{l-1}\binom{N-1}{l}\frac{1}{l^{n}}=\\
&=\sum_{m_1}^{N-1}\frac{1}{m_1}\dots\sum_{m_{n}=1}^{m_{n-1}}\frac{1}{m_{n}}={}^*H_{N-1}^{\overbrace{(1,\dots,1)}^{n}},
\end{aligned}
\end{equation}
where passing to the last line we used the identity \eqref{iterativesequence}. 

Thanks to the binomial transform, the identity \eqref{basestepharmonicbinomial} is equivalent to
\begin{equation}\label{harmonicidentitylevel1}
{}^*H_{N-1}^{(n)}=\sum_{j=1}^{N}(-1)^j\binom{N}{j}\,{}^*H_{j-1}^{\overbrace{(1,\dots,1)}^{n}}.
\end{equation}
The next step, with two arrays of indices, would be to prove the identity
\begin{equation}
{}^*H_{N-1}^{(s_2,\overbrace{1,\dots,1}^{n})}=\sum_{j=1}^{N}(-1)^j\binom{N}{j}\,{}^*H_{j-1}^{(\overbrace{1,\dots,1}^{s_2-1},n+1)},
\end{equation}
which, by taking the binomial transform, is equivalent to
\begin{equation}\label{harmonicidentitylevel2}
{}^*H_{N-1}^{(\overbrace{1,\dots,1}^{s_2-1},n+1)}=\sum_{j=1}^{N}(-1)^j\binom{N}{j}\,{}^*H_{j-1}^{(s_2,\overbrace{1,\dots,1}^{n})}.
\end{equation}
To prove \eqref{harmonicidentitylevel2}, we use the identity \eqref{harmonicidentitylevel1}. The left-hand side of \eqref{harmonicidentitylevel2} is then equal to
\begin{equation}
\begin{aligned}
&\sum_{k_1=1}^{N-1}\dots\sum_{k_{s_2-1}=1}^{k_{s_2-2}}\frac{1}{k_1\dots k_{s_2-1}}\,{}^*H_{k_{s_2-1}}^{(n+1)}=\sum_{k_1=1}^{N-1}\dots\sum_{k_{s_2-1}=1}^{k_{s_2-2}}\frac{1}{k_1\dots k_{s_2-1}}\,\sum_{j=1}^{k_{s_2-1}+1}(-1)^j\binom{k_{s_2-1}+1}{j}\,{}^*H_{j-1}^{(\overbrace{1,\dots,1}^{n+1})}=\\
&\sum_{k_1=1}^{N-1}\dots\sum_{k_{s_2-1}=1}^{k_{s_2-2}}\frac{1}{k_1\dots k_{s_2-1}}\,\sum_{m=0}^{k_{s_2-1}}(-1)^{m-1}\binom{k_{s_2-1}+1}{m+1}\,{}^*H_{m}^{(\overbrace{1,\dots,1}^{n+1})},
\end{aligned}
\end{equation}
where passing to the last line we defined $m=j-1$.
We now "bring down" the first index $1$ of the multiple harmonic star number and then change the order of the last two sums:
\begin{equation}
\begin{aligned}
&\sum_{k_1=1}^{N-1}\dots\sum_{k_{s_2-1}=1}^{k_{s_2-2}}\frac{1}{k_1\dots k_{s_2-1}}\,\sum_{m=0}^{k_{s_2-1}}(-1)^{m-1}\binom{k_{s_2-1}+1}{m+1}\sum_{j=1}^{m}\frac{1}{j}\,{}^*H_{j}^{(\overbrace{1,\dots,1}^{n})}=\\
&\sum_{k_1=1}^{N-1}\dots\sum_{k_{s_2-1}=1}^{k_{s_2-2}}\frac{1}{k_1\dots k_{s_2-1}}\,\sum_{j=1}^{k_{s_2-1}}\frac{(-1)^{j-1}}{j}\binom{k_{s_2-1}}{j}\,{}^*H_{j}^{(\overbrace{1,\dots,1}^{n})},
\end{aligned}
\end{equation}
where passing to the last line we also used the partial sum identity.
Renaming $j$ as $k_{s_2}$ and  using \eqref{iterativesequence}, we have
\begin{equation}\label{lhs_level2}
\begin{aligned}
&\sum_{k_1=1}^{N-1}\dots\sum_{k_{s_2-1}=1}^{k_{s_2-2}}\frac{1}{k_1\dots k_{s_2-1}}\,\sum_{k_{s_2}=1}^{k_{s_2-1}}\frac{(-1)^{k_{s_2}-1}}{k_{s_2}}\binom{k_{s_2-1}}{k_{s_2}}\,{}^*H_{k_{s_2}}^{(\overbrace{1,\dots,1}^{n})}=\\
&\sum_{k_1=1}^{N-1}\dots\sum_{k_{s_2}=1}^{k_{s_2-1}}\frac{1}{k_1\dots k_{s_2}}\,\sum_{l=1}^{k_{s_2}}(-1)^{l-1}\binom{k_{s_2}}{l}\,{}^*H_{l}^{(\overbrace{1,\dots,1}^{n})}.
\end{aligned}
\end{equation}
The right-hand side of \eqref{harmonicidentitylevel2} is, swapping the summations and using the partial sum identity,
\begin{equation}\label{rhs_level2}
\begin{aligned}
&\sum_{j=1}^{N}(-1)^j\binom{N}{j}\,{}^*H_{j-1}^{(s_2,\overbrace{1,\dots,1}^{n})}=\sum_{j=1}^{N}(-1)^j\binom{N}{j}\sum_{l=1}^{j-1}\frac{1}{l^{s_2}}\,{}^*H_l^{(\overbrace{1,\dots,1}^{n})}=\\
&\sum_{l=1}^{N-1}\frac{(-1)^{l-1}}{l^{s_2}}\binom{N-1}{l}\,{}^*H_l^{(\overbrace{1,\dots,1}^{n})}.
\end{aligned}
\end{equation}
The equality between \eqref{lhs_level2} and \eqref{rhs_level2} then follows by applying \eqref{iterativesequence}, which completes the proof of \eqref{harmonicidentitylevel2}.

The generic inductive step can be proved by following the same reasoning. 
The key observation is that the operator
\begin{equation}
\sum_{j=1}^N(-1)^j\binom{N}{j}
\end{equation}
applied to a multiple harmonic star number of the form ${}^*H_{N-1}^{(\alpha_1,\beta_1,\dots,\alpha_n,\beta_n)}$, where $\alpha_i=(\overbrace{1,\dots,1}^{s_i})$ and $\beta_i=r_i\in\mathbb{Z}_{>1}$, transforms it into ${}^*H_{N-1}^{(\beta'_1,\alpha'_1,\dots,\beta'_n,\alpha'_n)}$ with  $\alpha'_i=(\overbrace{1,\dots,1}^{s'_i})$ and $\beta'_i=r'_i\in\mathbb{Z}_{>1}$. The same is true for ${}^*H_{N-1}^{(\alpha_1,\beta_1,\dots,\alpha_n,\beta_n,\alpha_{n+1})}$, which transforms into ${}^*H_{N-1}^{(\beta'_1,\alpha'_1,\dots,\beta'_n,\alpha'_n,\beta'_{n+1})}$.
For the inductive step with a generic number $n>2$ of (arrays of) indices, it is sufficient to prove first
\begin{equation}\label{inductiveidentityharmonic}
{}^*H_{N-1}^{(\overbrace{1,\dots,1}^{r_2-1},s_2+1,\dots,s_n+1,\overbrace{1,\dots,1}^{r_{n+1}})}=\sum_{j=1}^{N}(-1)^j\binom{N}{j}\,{}^*H_{j-1}^{(r_2,\overbrace{1,\dots,1}^{s_2-1},\dots,\overbrace{1,\dots,1}^{s_n-1},r_{n+1}+1)}
\end{equation}
and then
\begin{equation}\label{inductiveidentityharmonicn+1}
{}^*H_{N-1}^{(\overbrace{1,\dots,1}^{r_2-1},s_2+1,\dots,\overbrace{1,\dots,1}^{r_{n+1}-1},s_{n+1}+1)}=\sum_{j=1}^{N}(-1)^j\binom{N}{j}\,{}^*H_{j-1}^{(r_2,\overbrace{1,\dots,1}^{s_2-1},\dots,r_{n+1}+1,\overbrace{1,\dots,1}^{s_{n+1}})}
\end{equation}
from \eqref{inductiveidentityharmonic}. As will be clear from the proof of \eqref{inductiveidentityharmonic}, the result \eqref{inductiveidentityharmonicn+1} can be obtained straightforwardly.

To prove \eqref{inductiveidentityharmonic}, we start by rewriting the left-hand side using the inductive hypothesis:
\begin{equation}
\begin{aligned}
&{}^*H_{N-1}^{(\overbrace{1,\dots,1}^{r_2-1},s_2+1,\dots,s_n+1,\overbrace{1,\dots,1}^{r_{n+1}})}=\sum_{k_1=1}^{N-1}\dots\sum_{k_{r_2-1}=1}^{k_{r_2-2}}\frac{1}{k_1\dots k_{r_2-1}}\,{}^*H_{k_{r_2-1}}^{(s_2+1,\dots,s_n+1,\overbrace{1,\dots,1}^{r_{n+1}})}=\\
&\sum_{k_1=1}^{N-1}\dots\sum_{k_{r_2-1}=1}^{k_{r_2-2}}\frac{1}{k_1\dots k_{r_2-1}}\sum_{j=1}^{k_{r_2-1}+1}(-1)^j\binom{k_{r_2-1}+1}{j}\,{}^*H_{j-1}^{(\overbrace{1,\dots,1}^{s_2},r_3+1,\dots,r_{n+1}+1)}=\\
&\sum_{k_1=1}^{N-1}\dots\sum_{k_{r_2-1}=1}^{k_{r_2-2}}\frac{1}{k_1\dots k_{r_2-1}}\sum_{m=0}^{k_{r_2-1}}(-1)^{m-1}\binom{k_{r_2-1}+1}{m+1}\,{}^*H_{m}^{(\overbrace{1,\dots,1}^{s_2},r_3+1,\dots,r_{n+1}+1)},
\end{aligned}
\end{equation}
where passing to the last line we defined $m=j-1$.
We now "bring down" the first index $1$ of the multiple harmonic star number:
\begin{equation}\label{step1n}
\begin{aligned}
&\sum_{k_1=1}^{N-1}\dots\sum_{k_{r_2-1}=1}^{k_{r_2-2}}\frac{1}{k_1\dots k_{r_2-1}}\sum_{m=0}^{k_{r_2-1}}(-1)^{m-1}\binom{k_{r_2-1}+1}{m+1}\sum_{j=1}^m\frac{1}{j}\,{}^*H_{j}^{(\overbrace{1,\dots,1}^{s_2-1},r_3+1,\dots,r_{n+1}+1)}=\\
&\sum_{k_1=1}^{N-1}\dots\sum_{k_{r_2-1}=1}^{k_{r_2-2}}\frac{1}{k_1\dots k_{r_2-1}}\sum_{j=1}^{k_{r_2-1}}\frac{(-1)^{j-1}}{j}\binom{k_{r_2-1}}{j}\,{}^*H_{j}^{(\overbrace{1,\dots,1}^{s_2-1},r_3+1,\dots,r_{n+1}+1)},
\end{aligned}
\end{equation}
where passing to the second line we swapped the summations and used the partial sum identity. Using identity \eqref{iterativesequence}, we have
\begin{equation}
\begin{aligned}
&\sum_{k_1=1}^{N-1}\dots\sum_{k_{r_2-1}=1}^{k_{r_2-2}}\frac{1}{k_1\dots k_{r_2-1}}\sum_{j=1}^{k_{r_2-1}}\frac{(-1)^{j-1}}{j}\binom{k_{r_2-1}}{j}\,{}^*H_{j}^{(\overbrace{1,\dots,1}^{s_2-1},r_3+1,\dots,r_{n+1}+1)}=\\
&\sum_{k_1=1}^{N-1}\dots\sum_{k_{r_2}=1}^{k_{r_2-1}}\frac{1}{k_1\dots k_{r_2}}\sum_{l=1}^{k_{r_2}}(-1)^{l-1}\binom{k_{r_2}}{l}\,{}^*H_{l}^{(\overbrace{1,\dots,1}^{s_2-1},r_3+1,\dots,r_{n+1}+1)}.
\end{aligned}
\end{equation}
Finally, we can rewrite the right-hand side of \eqref{inductiveidentityharmonic} as follows
\begin{equation}\label{step2n}
\begin{aligned}
&\sum_{j=1}^{N}(-1)^j\binom{N}{j}\,{}^*H_{j-1}^{(r_2,\overbrace{1,\dots,1}^{s_2-1},\dots,\overbrace{1,\dots,1}^{s_n-1},r_{n+1}+1)}=\sum_{j=1}^{N}(-1)^j\binom{N}{j}\sum_{l=1}^{j-1}\frac{1}{l^{r_2}}\,{}^*H_{l}^{(\overbrace{1,\dots,1}^{s_2-1},r_3+1,\dots,r_{n+1}+1)}=\\
&=\sum_{l=1}^{N-1}\frac{1}{l^{r_2}}\,{}^*H_{l}^{(\overbrace{1,\dots,1}^{s_2-1},r_3+1,\dots,r_{n+1}+1)}\sum_{j=l+1}^{N}(-1)^j\binom{N}{j}=\\
&=\sum_{l=1}^{N-1}\frac{(-1)^{l-1}}{l^{r_2}}\binom{N-1}{l}\,{}^*H_{l}^{(\overbrace{1,\dots,1}^{s_2-1},r_3+1,\dots,r_{n+1}+1)}.
\end{aligned}
\end{equation}
The equality between the right-hand side and left-hand side of \eqref{inductiveidentityharmonic}  then follows from the identity \eqref{iterativesequence}.

To prove \eqref{inductiveidentityharmonicn+1}, we rewrite the left-hand side using \eqref{inductiveidentityharmonic}:
\begin{equation}
\begin{aligned}
&{}^*H_{N-1}^{(\overbrace{1,\dots,1}^{r_2-1},s_2+1,\dots,s_n+1,\overbrace{1,\dots,1}^{r_{n+1}-1},s_{n+1}+1)}=\sum_{k_1=1}^{N-1}\dots\sum_{k_{r_2-1}=1}^{k_{r_2-2}}\frac{1}{k_1\dots k_{r_2-1}}\,{}^*H_{k_{r_2-1}}^{(s_2+1,\dots,s_n+1,\overbrace{1,\dots,1}^{r_{n+1}-1},s_{n+1}+1)}=\\
&\sum_{k_1=1}^{N-1}\dots\sum_{k_{r_2-1}=1}^{k_{r_2-2}}\frac{1}{k_1\dots k_{r_2-1}}\sum_{j=1}^{k_{r_2-1}+1}(-1)^j\binom{k_{r_2-1}+1}{j}\,{}^*H_{j-1}^{(\overbrace{1,\dots,1}^{s_2},r_3+1,\dots,r_{n+1}+1,\overbrace{1,\dots,1}^{s_{n+1}})}.
\end{aligned}
\end{equation}
The rest of the proof follows the same steps as in (\ref{step1n}--\ref{step2n}), where the set of indices $\overbrace{1,\dots,1}^{s_{n+1}}$ is added to each multiple star harmonic number. 

We finish the Appendix by proving identity \eqref{basestepharmonicidentity}. We can rewrite the left-hand side by using the binomial transform of \eqref{keyidentityharmonic1}:
\begin{equation}
\begin{aligned}
{}^*H_k^{(r_2+1,\overbrace{1,\dots,1}^{s_2-1},\dots,\overbrace{1,\dots,1}^{s_n})}=\,&\sum_{j=1}^{k+1}(-1)^j\binom{k+1}{j}\,{}^*H_{j-1}^{(\overbrace{1,\dots,1}^{r_2},s_2+1,\dots,s_n+1)}=\\
&\sum_{m=0}^{k}(-1)^{m-1}\binom{k+1}{m+1}\,{}^*H_{m}^{(\overbrace{1,\dots,1}^{r_2},s_2+1,\dots,s_n+1)},
\end{aligned}
\end{equation}
where we defined $m=j-1$. We bring down the first index $1$ of the multiple harmonic star number, swap the summations, and use the partial sum identity:
\begin{equation}
\begin{aligned}
&\sum_{m=0}^{k}(-1)^{m-1}\binom{k+1}{m+1}\,{}^*H_{m}^{(\overbrace{1,\dots,1}^{r_2},s_2+1,\dots,s_n+1)}=\\
&\sum_{m=0}^{k}(-1)^{m-1}\binom{k+1}{m+1}\sum_{l=1}^m\frac{1}{l}\,{}^*H_{l}^{(\overbrace{1,\dots,1}^{r_2-1},s_2+1,\dots,s_n+1)}=\\
&\sum_{l=1}^{k}\frac{(-1)^{l-1}}{l}\binom{k}{l}\,{}^*H_{l}^{(\overbrace{1,\dots,1}^{r_2-1},s_2+1,\dots,s_n+1)},
\end{aligned}
\end{equation}
which is equal to the right-hand side of \eqref{basestepharmonicidentity}.

\section{Proof of the relations for integrals of levels 2 and 3}
\label{app:rel23_proof}
To prove (\ref{eq:ELimn_sub}), we apply an induction with a base case (\ref{eq:ELi1n_sub}). Assuming that (\ref{eq:ELimn_sub}) is correct for a fixed value of $m$, we derive the relation for $\text{ELi}_{m+1,n}\brc{z}$:
\eq{\text{ELi}_{m+1,n}\brc{z}=\int_{-\infty}^{z}\text{ELi}_{m,n}\brc{t} \frac{\rmd t}{t}.}
Using (\ref{eq:ELimn_sub}), we split the integral into two parts:
\begin{itemize}
\item[1.] $\displaystyle \int_{-\infty}^{z}\brc{ \text{EL}_{m,n}\brc{t} + \sum_{k=1}^{m} \frac{\brc{-1}^{k-1}}{\brc{k-1}!\brc{m-k}!} \,
\text{cLi}_{k,n} \log\brc{-t}^{m-k}} \frac{\rmd t}{t}$,
\item[2.] $\displaystyle \int_{-\infty}^{z} \,\sum_{k=0}^{n}\sum_{l=k}^{n}
\frac{\brc{-1}^{l+1}}{k!\brc{n-l}!} \, \binom{l+m-k-1}{m-1} \Gamma^{\brc{k}}\brc{1} \log\brc{t}^{n-l} \text{ELi}_{l+m-k}\brc{t}
\frac{\rmd t}{t}$.
\end{itemize}
In the first case, we split the integration domain and integrate by parts:
\eq{m\geq 2:\quad \int_{-\infty}^{-1}\brc{ \text{EL}_{m,n}\brc{t} + \sum_{k=1}^{m} \frac{\brc{-1}^{k-1}}{\brc{k-1}!\brc{m-k}!} \,
\text{cLi}_{k,n} \log\brc{-t}^{m-k}} \rmd \log\brc{-t} =}
\eq{=- \int_{-\infty}^{-1}\log\brc{-t} \brc{ \text{EL}_{m-1,n}\brc{t} + \sum_{k=1}^{m-1} \frac{\brc{-1}^{k-1}}{\brc{k-1}!\brc{m-k-1}!} \,
\text{cLi}_{k,n} \log\brc{-t}^{m-k-1}} \frac{\rmd t}{t},}
where we took into account that the first expression in brackets tends to zero as $t\rightarrow-\infty$ (and the same is true for the second expression). Continuing integrating by parts, we arrive at
\eqlb{eq:ELimn_sub_proof1}{m\geq 1:\quad\frac{\brc{-1}^m}{m!} \int_{-\infty}^{-1}\log\brc{-t}^m \brc{-\frac{\rme^t}{t}} \text{EL}_{n}\brc{-t} \rmd t =
\frac{\brc{-1}^m}{m!} \int_{1}^{\infty} \rme^{-t} \log\brc{t}^{m}\, \text{EL}_{n}\brc{t} \frac{\rmd t}{t}.}
Next, we integrate from $-1$ to $z$:
\eqlb{eq:ELimn_sub_proof2}{\int_{-1}^{z} \text{EL}_{m,n}\brc{t} \frac{\rmd t}{t}= \text{EL}_{m+1,n}\brc{z} + \int_{-1}^{0} \text{EL}_{m,n}\brc{t} \rmd \log\brc{-t},}
\eq{\int_{-1}^{z}\, \sum_{k=1}^{m} \frac{\brc{-1}^{k-1} \text{cLi}_{k,n}}{\brc{k-1}!\brc{m-k}!} \, \log\brc{-t}^{m-k} \frac{\rmd t}{t} = \sum_{k=1}^{m} \frac{\brc{-1}^{k-1}\text{cLi}_{k,n}}{\brc{k-1}!\brc{m+1-k}!} \, \log\brc{-z}^{m+1-k}.}
The integral on the right-hand side of (\ref{eq:ELimn_sub_proof2}) can be rewritten to match the integral in (\ref{eq:ELimn_sub_proof1}):
\eqlb{eq:ELimn_sub_proof3}{\int_{-1}^{0} \text{EL}_{m,n}\brc{t} \rmd \log\brc{-t} =
-\frac{1}{2} \int_{-1}^{0} \text{EL}_{m-1,n}\brc{t} \rmd \log\brc{-t}^2= \dots=
\frac{\brc{-1}^{m-1}}{m!} \int_{-1}^{0} \text{EL}_{1,n}\brc{t} \rmd \log\brc{-t}^m=}
\eqn{=\frac{\brc{-1}^{m}}{m!} \int_{-1}^{0} \log\brc{-t}^m \brc{-\frac{\rme^t}{t}} \text{EL}_{n}\brc{-t} \rmd t=
\frac{\brc{-1}^m}{m!} \int_{0}^{1} \rme^{-t} \log\brc{t}^{m}\, \text{EL}_{n}\brc{t} \frac{\rmd t}{t}.}
Thus, the overall result for the first part is
\eq{\text{EL}_{m+1,n}\brc{z} + \sum_{k=1}^{m+1} \frac{\brc{-1}^{k-1}}{\brc{k-1}!\brc{m+1-k}!} \,
\text{cLi}_{k,n} \log\brc{-z}^{m+1-k},}
where we introduced a new constant $\text{cLi}_{m+1,n}$ that matches the general definition (\ref{eq:cLimn_def}).

In the second case, we integrate each term of the sum by parts until the power of the logarithm function is reduced to zero:
\eqn{\int_{-\infty}^{z} \log\brc{t}^{n-l} \text{ELi}_{l+m-k}\brc{t} \frac{\rmd t}{t} = \int_{-\infty}^{z} \log\brc{t}^{n-l} \rmd \, \text{ELi}_{l+m-k+1}\brc{t}= \log\brc{z}^{n-l} \text{ELi}_{l+m-k+1}\brc{z} }
\eqlb{eq:ELimn_sub_proof4}{- \brc{n-l} \int_{-\infty}^{z} \log\brc{t}^{n-l-1} \rmd \, \text{ELi}_{l+m-k+2}\brc{t} =\dots =}
\eqn{=\sum_{j=l}^{n} \brc{-1}^{j-l} \frac{\brc{n-l}!}{\brc{n-j}!} \log\brc{z}^{n-j} \text{ELi}_{j+m-k+1}\brc{z}.}
Substituting this back into the sum, we get
\eq{\sum_{k=0}^{n}\sum_{l=k}^{n}\sum_{j=l}^{n}
\frac{\brc{-1}^{j+1}}{k!\brc{n-j}!} \, \binom{l+m-k-1}{m-1} \Gamma^{\brc{k}}\brc{1} \log\brc{z}^{n-j} \text{ELi}_{j+m-k+1}\brc{z}.}
Now, we can exchange the summation order considering $k\leq l \leq j \leq n$:
\eq{\sum_{k=0}^{n}\sum_{l=k}^{n}\sum_{j=l}^{n}\dots \longrightarrow \sum_{k=0}^{n}\sum_{j=k}^{n}\sum_{l=k}^{j}\dots \,.}
Then, the sum over index $l$ can be isolated and computed using the Hockey-Stick identity:
\eqlb{eq:hockey_id}{\sum_{l=k}^{j} \binom{l+m-k-1}{m-1}=\binom{j+m-k}{m}.}
Finally, combining the results for both parts, we arrive at (\ref{eq:ELimn_sub}) with an index $m$ shifted by $1$:
\eq{\text{ELi}_{m+1,n}\brc{z}= \text{EL}_{m+1,n}\brc{z} + \sum_{k=1}^{m+1} \frac{\brc{-1}^{k-1}}{\brc{k-1}!\brc{m+1-k}!} \,
\text{cLi}_{k,n} \log\brc{-z}^{m+1-k}}
\eq{+\sum_{k=0}^{n}\sum_{j=k}^{n} \frac{\brc{-1}^{j+1}}{k!\brc{n-j}!} \, \binom{j+m-k}{m} \Gamma^{\brc{k}}\brc{1} \log\brc{z}^{n-j} \text{ELi}_{j+m+1-k}\brc{z},}
which ends the proof by induction.

The proof of (\ref{eq:ELi_lmn_sub}) is essentially the same. We assume that (\ref{eq:ELi_lmn_sub}) is correct for a fixed value of $l$ and derive the relation for $l+1$ using the definition:
\eqlb{eq:ELi_lmn_sub_proof0}{\text{ELi}_{l+1,m,n}\brc{z}=\int_{-\infty}^{z}\text{ELi}_{l,m,n}\brc{t} \frac{\rmd t}{t}.}
The three relevant integrals in this case are
\begin{itemize}
\item[1.] $\displaystyle \int_{-\infty}^{z}\brc{\text{EL}_{l,m,n}\brc{t}+ \sum_{k=1}^{l} \frac{\brc{-1}^{k-1}}{\brc{k-1}!\brc{l-k}!} \,
\text{cLi}_{k,m,n} \log\brc{-t}^{l-k} } \frac{\rmd t}{t}$,
\item[2.] $\displaystyle \int_{-\infty}^{z} \log\brc{t}^{m-j} \text{ELi}_{l+j-k}\brc{t} \frac{\rmd t}{t}$,
\item[3.] $\displaystyle \int_{-\infty}^{z} \log\brc{-t}^{n-i} \text{ELi}_{l+i-j,m+j-k}\brc{t} \frac{\rmd t}{t}$.
\end{itemize}
Taking the first integral from $-\infty$ to $-1$, in complete analogy with (\ref{eq:ELimn_sub_proof1}), we get
\eqlb{eq:ELi_lmn_sub_proof1}{\int_{-\infty}^{-1}\brc{\text{EL}_{l,m,n}\brc{t}+ \sum_{k=1}^{l} \frac{\brc{-1}^{k-1}}{\brc{k-1}!\brc{l-k}!} \,
\text{cLi}_{k,m,n} \log\brc{-t}^{l-k} } \frac{\rmd t}{t} = \frac{\brc{-1}^l}{l!} \int_{1}^{\infty} \rme^{-t} \log\brc{t}^{l}\, \text{EL}_{m,n}\brc{t} \frac{\rmd t}{t}.}
The remaining part of the first integral is
\eqlb{eq:ELi_lmn_sub_proof2}{\int_{-1}^{z}\brc{\text{EL}_{l,m,n}\brc{t}+ \sum_{k=1}^{l} \frac{\brc{-1}^{k-1}}{\brc{k-1}!\brc{l-k}!} \,
\text{cLi}_{k,m,n} \log\brc{-t}^{l-k} } \frac{\rmd t}{t} =\text{EL}_{l+1,m,n}\brc{z}+\int_{-1}^{0} \text{EL}_{l,m,n}\brc{t}\frac{\rmd t}{t} }
\eqn{+ \sum_{k=1}^{l} \frac{\brc{-1}^{k-1}}{\brc{k-1}!\brc{l+1-k}!} \, \text{cLi}_{k,m,n} \log\brc{-z}^{l+1-k}.}
Following (\ref{eq:ELimn_sub_proof3}), we rewrite the integral on the right-hand side of (\ref{eq:ELi_lmn_sub_proof2}) to match the result in (\ref{eq:ELi_lmn_sub_proof1}):
\eq{\int_{-1}^{0} \text{EL}_{l,m,n}\brc{t}\frac{\rmd t}{t}= \frac{\brc{-1}^l}{l!} \int_{0}^{1} \rme^{-t} \log\brc{t}^{l}\, \text{EL}_{m,n}\brc{t} \frac{\rmd t}{t}.}
Since, by definition,
\eq{\text{cLi}_{l+1,m,n}= \int_{0}^{\infty} \rme^{-t} \log\brc{t}^{l}\, \text{EL}_{m,n}\brc{t} \frac{\rmd t}{t},}
the overall result for the first integral is
\eq{\text{EL}_{l+1,m,n}\brc{z}+ \sum_{k=1}^{l+1} \frac{\brc{-1}^{k-1}}{\brc{k-1}!\brc{l+1-k}!} \, \text{cLi}_{k,m,n} \log\brc{-z}^{l+1-k}.}

The second and third integrals follow directly from (\ref{eq:ELimn_sub_proof4}):
\eq{\int_{-\infty}^{z} \log\brc{t}^{m-j} \text{ELi}_{l+j-k}\brc{t} \frac{\rmd t}{t}=\sum_{i=j}^{m} \brc{-1}^{i-j} \frac{\brc{m-j}!}{\brc{m-i}!} \log\brc{z}^{m-i} \text{ELi}_{l+i-k+1}\brc{z},}
\eqlb{eq:ELi_4_second}{\int_{-\infty}^{z} \log\brc{-t}^{n-i} \text{ELi}_{l+i-j,m+j-k}\brc{t} \frac{\rmd t}{t} = \sum_{i_1=i}^{n} \brc{-1}^{i_1-i} \frac{\brc{n-i}!}{\brc{n-i_1}!} \log\brc{-z}^{n-i_1} \text{ELi}_{l+i_1-j+1,m+j-k}\brc{z}.}
Putting these back into (\ref{eq:ELi_lmn_sub_proof0}) gives
\eqlb{eq:ELi_lmn_sub_proof3}{\sum_{k=1}^{m} \sum_{j=k}^{m}\sum_{i=j}^{m} \frac{\brc{-1}^{i}}{\brc{k-1}!\brc{m-i}!} \, \tbinom{l+j-k-1}{l-1} \text{cLi}_{k,n} \log\brc{z}^{m-i} \text{ELi}_{l+i-k+1}\brc{z}}
and
\eqlb{eq:ELi_lmn_sub_proof4}{\sum_{k=0}^{n}\sum_{j=k}^{n}\sum_{i=j}^{n}\sum_{i_1=i}^{n}
\frac{\brc{-1}^{i_1+1}}{k!\brc{n-i_1}!} \, \tbinom{l+i-j-1}{l-1} \tbinom{m+j-k-1}{m-1} \Gamma^{\brc{k}}\brc{1} \log\brc{-z}^{n-i_1} \text{ELi}_{l+i_1-j+1,m+j-k}\brc{z}.}
Using the Hockey-Stick identity (\ref{eq:hockey_id}), we can sum over index $j$ in (\ref{eq:ELi_lmn_sub_proof3}) and over the index $i$ in (\ref{eq:ELi_lmn_sub_proof4}), which results in
\eqlb{eq:ELi_HS1}{\sum_{k=1}^{m}\sum_{i=k}^{m} \frac{\brc{-1}^{i}}{\brc{k-1}!\brc{m-i}!} \, \tbinom{l+i-k}{l} \text{cLi}_{k,n} \log\brc{z}^{m-i} \text{ELi}_{l+i-k+1}\brc{z}}
and
\eqlb{eq:ELi_HS2}{\sum_{k=0}^{n}\sum_{j=k}^{n}\sum_{i_1=j}^{n}
\frac{\brc{-1}^{i_1+1}}{k!\brc{n-i_1}!} \, \tbinom{l+i_1-j}{l} \tbinom{m+j-k-1}{m-1} \Gamma^{\brc{k}}\brc{1} \log\brc{-z}^{n-i_1} \text{ELi}_{l+i_1-j+1,m+j-k}\brc{z}.}
Collecting all three parts together and renaming some of the summation indices, we arrive at
\eqn{\text{ELi}_{l+1,m,n}\brc{z}=\text{EL}_{l+1,m,n}\brc{z}+ \sum_{k=1}^{l+1} \frac{\brc{-1}^{k-1}}{\brc{k-1}!\brc{l+1-k}!} \, \text{cLi}_{k,m,n} \log\brc{-z}^{l+1-k}}
\eq{+\sum_{k=1}^{m}\sum_{j=k}^{m} \frac{\brc{-1}^{j}}{\brc{k-1}!\brc{m-j}!} \, \tbinom{l+j-k}{l} \text{cLi}_{k,n} \log\brc{z}^{m-j} \text{ELi}_{l+j-k+1}\brc{z}}
\eqn{+\sum_{k=0}^{n}\sum_{j=k}^{n}\sum_{i=j}^{n}
\frac{\brc{-1}^{i+1}}{k!\brc{n-i}!} \, \tbinom{l+i-j}{l} \tbinom{m+j-k-1}{m-1} \Gamma^{\brc{k}}\brc{1} \log\brc{-z}^{n-i} \text{ELi}_{l+i-j+1,m+j-k}\brc{z},}
which finalizes the proof of (\ref{eq:ELi_lmn_sub}).

\section{Proof of the general relation for integrals of level \texorpdfstring{$n\geq 2$}{}}
\label{app:reln_proof}

Let us prove the general relation \eqref{eq:ELi_sub_gen} for a fixed $n>2$ by induction. The first cases $n=2$ and $n=3$ were proven in the previous Appendix. We start from the base case $s_1=1$:
\begin{equation}
\text{ELi}_{1,s_2,\dots,s_n}\brc{z}=-\int_{-\infty}^z\frac{\mathrm{e}^t}{t}\text{ELi}_{s_2,\dots,s_n}\brc{-t}\mathrm{d}t.
\end{equation}
By induction on the level $n$,
\begin{equation}\label{eq:ELi_inductive_s1=1}
\begin{aligned}
\text{ELi}_{s_2,\dots,s_n}\brc{z} =\,& \text{EL}_{s_2,\dots,s_n}\brc{z}+ \sum_{k_1=1}^{s_2} \frac{\brc{-1}^{k_1-1}}{\brc{k_1-1}!\brc{s_2-k_1}!} \, \text{cLi}_{k_1,s_3,\dots,s_n} \log\brc{-z}^{s_2-k_1}+\\
&\sum_{i=3}^{n-1}\sum_{s_i\geq k_1\geq \dots\geq k_{i-1}\geq 1} \frac{\brc{-1}^{k_1}}{\brc{k_{i-1}-1}!\brc{s_i-k_1}!} \prod_{j=1}^{i-2} \tbinom{s_j-1+k_j-k_{j+1}}{s_j-1}\, \text{cLi}_{k_{i-1},s_{i+1},\dots,s_n}  \times\\
& \log(\brc{-1}^{i-1} z)^{s_i-k_1}\text{ELi}_{s_2+k_1-k_2,\dots,s_{i-1}+k_{i-2}-k_{i-1}}\brc{z}+\\
&\sum_{s_n\geq k_1\geq\dots\geq k_{n-1}\geq 0} \frac{\brc{-1}^{k_1+1}}{k_{n-1}!\brc{s_n-k_1}!}\prod_{j=1}^{n-2} \tbinom{s_j-1+k_j-k_{j+1}}{s_j-1} \, \Gamma^{\brc{k_{n-1}}}\brc{1} \times\\
& \log\brc{\brc{-1}^{n-1} z}^{s_n-k_1}\text{ELi}_{s_2+k_1-k_2,\dots,s_{n-1}+k_{n-2}-k_{n-1}}\brc{z}.
\end{aligned}
\end{equation}
We split the integral into three kinds of integrals:
\begin{itemize}
\item[1.] $\displaystyle -\int_{-\infty}^z\frac{\mathrm{e}^t}{t}\text{EL}_{s_2,\dots,s_n}\brc{-t}\mathrm{d}t$,
\item[2.] $\displaystyle -\int_{-\infty}^z\frac{\mathrm{e}^t}{t}\sum_{k_1=1}^{s_2}\frac{(-1)^{k_1-1}}{(k_1-1)!\,(s_2-k_1)!}\text{cLi}_{k_1,s_3,\dots,s_n}\log(t)^{s_2-k_1}\mathrm{d}t$,
\item[3.] $\displaystyle -\int_{-\infty}^z\frac{\mathrm{e}^t}{t}\log\left((-1)^{i}t\right)^{s_i-k_1}\text{ELi}_{s_2+k_1-k_2,\dots,s_{i-1}+k_{i-2}-k_{i-1}}\brc{-t}\mathrm{d}t$,\quad $i\in\{3,\dots,n\}$.
\end{itemize}
The first integral is
\begin{equation}
-\int_{-\infty}^z\frac{\mathrm{e}^t}{t}\text{EL}_{s_2,\dots,s_n}\brc{-t}\mathrm{d}t=\text{EL}_{1,s_2,\dots,s_n}\brc{z}+\text{cLi}_{1,s_2,\dots,s_n},
\end{equation}
where $\text{cLi}_{1,s_2,\dots,s_n}$ matches the definition in \eqref{eq:cLi_gen_def}.
For the second integral, we integrate by parts until the power of the logarithm vanishes, leading to
\begin{equation}
\sum_{k_1=1}^{s_2}\sum_{k_2=k_1}^{s_2}\frac{(-1)^{k_2}}{(k_1-1)!\,(s_2-k_2)!}\text{cLi}_{k_1,s_3,\dots,s_n}\log(z)^{s_2-k_2}\text{ELi}_{k_2-k_1+1}\brc{z}.
\end{equation}
The third integral can be computed in complete analogy with (\ref{eq:ELimn_sub_proof4}):
\begin{equation}
\begin{aligned}
&\int_{-\infty}^z\log\left((-1)^{i}t\right)^{s_i-k_1}\mathrm{d}\,\text{ELi}_{1,s_2+k_1-k_2,\dots,s_{i-1}+k_{i-2}-k_{i-1}}\brc{t}=\\
&\sum_{l=k_1}^{s_i}(-1)^{l-k_1}\frac{(s_i-k_1)!}{(s_i-l)!}\log\left((-1)^{i}z\right)^{s_i-l}\text{ELi}_{l-k_1+1,s_2+k_1-k_2,\dots,s_{i-1}+k_{i-2}-k_{i-1}}\brc{z}.
\end{aligned}
\end{equation}
Reintroducing the summations and using the Hockey-Stick identity \eqref{eq:hockey_id}, we have for the integral of the second through fifth line of \eqref{eq:ELi_inductive_s1=1}:
\begin{equation}
\begin{aligned}
&\sum_{i=3}^{n-1}\sum_{s_i\geq k_1\geq \dots\geq k_{i-1}\geq 1} \frac{\brc{-1}^{k_1}}{\brc{k_{i-1}-1}!\brc{s_i-k_1}!} \prod_{j=1}^{i-2} \tbinom{s_j-1+k_j-k_{j+1}}{s_j-1}\, \text{cLi}_{k_{i-1},s_{i+1},\dots,s_n}\times\\
&\sum_{l=k_1}^{s_i}(-1)^{l-k_1}\frac{(s_i-k_1)!}{(s_i-l)!}\log\left((-1)^{i}z\right)^{s_i-l}\text{ELi}_{l-k_1+1,s_2+k_1-k_2,\dots,s_{i-1}+k_{i-2}-k_{i-1}}\brc{z}+\\
&\sum_{s_n\geq k_1\geq\dots\geq k_{n-1}\geq 0} \frac{\brc{-1}^{k_1+1}}{k_{n-1}!\brc{s_n-k_1}!}\prod_{j=1}^{n-2} \tbinom{s_j-1+k_j-k_{j+1}}{s_j-1} \, \Gamma^{\brc{k_{n-1}}}\brc{1}\times\\
&\sum_{l=k_1}^{s_n}(-1)^{l-k_1}\frac{(s_i-k_1)!}{(s_i-l)!}\log\left((-1)^{n}z\right)^{s_n-l}\text{ELi}_{l-k_1+1,s_2+k_1-k_2,\dots,s_{n-1}+k_{n-2}-k_{n-1}}\brc{z}=\\
&\sum_{i=3}^{n-1}\sum_{s_i\geq k_1\geq \dots\geq k_i\geq 1} \frac{\brc{-1}^{k_1}}{\brc{k_i-1}!\brc{s_i-k_1}!} \prod_{j=1}^{i-1} \tbinom{s_j-1+k_j-k_{j+1}}{s_j-1}\, \text{cLi}_{k_i,s_{i+1},\dots,s_n}\times\\
&\log\left((-1)^{i}z\right)^{s_i-k_1}\text{ELi}_{k_1-k_2+1,s_2+k_2-k_3,\dots,s_{i-1}+k_{i-1}-k_i}\brc{z}+\\
&\sum_{s_n\geq k_1\geq\dots\geq k_{n}\geq 0} \frac{\brc{-1}^{k_1+1}}{k_{n}!\brc{s_n-k_1}!}\prod_{j=1}^{n-1} \tbinom{s_j-1+k_j-k_{j+1}}{s_j-1} \, \Gamma^{\brc{k_{n}}}\brc{1}\times\\
&\log\left((-1)^{n}z\right)^{s_n-k_1}\text{ELi}_{k_1-k_2+1,s_2+k_2-k_3,\dots,s_{n-1}+k_{n-1}-k_{n}}\brc{z},
\end{aligned}
\end{equation}
where we renamed the indices as
\begin{equation}
(l,k_1,\dots,k_{n-1})\mapsto (k_1,\dots,k_{n}).
\end{equation}
Putting together the results from the three integrals, swapping the indices $k_1\leftrightarrow k_2$ for the result of the second integral, we conclude the proof of the base step on $s_1$ at a given level $n$:
\begin{equation}
\begin{aligned}
\text{ELi}_{1,s_2,\dots,s_n}\brc{z}=&\,\text{EL}_{1,s_2,\dots,s_n}\brc{z}+\text{cLi}_{1,s_2,\dots,s_n}+\\
&\sum_{i=2}^{n-1}\sum_{s_i\geq k_1\geq \dots\geq k_i\geq 1} \frac{\brc{-1}^{k_1}}{\brc{k_i-1}!\brc{s_i-k_1}!} \prod_{j=1}^{i-1} \tbinom{s_j-1+k_j-k_{j+1}}{s_j-1}\, \text{cLi}_{k_i,s_{i+1},\dots,s_n}\times\\
&\log\left((-1)^{i}z\right)^{s_i-k_1}\text{ELi}_{k_1-k_2+1,s_2+k_2-k_3,\dots,s_{i-1}+k_{i-1}-k_i}\brc{z}+\\
&\sum_{s_n\geq k_1\geq\dots\geq k_{n}\geq 0} \frac{\brc{-1}^{k_1+1}}{k_{n}!\brc{s_n-k_1}!}\prod_{j=1}^{n-1} \tbinom{s_j-1+k_j-k_{j+1}}{s_j-1} \, \Gamma^{\brc{k_{n}}}\brc{1}\times\\
&\log\left((-1)^{n}z\right)^{s_n-k_1}\text{ELi}_{k_1-k_2+1,s_2+k_2-k_3,\dots,s_{n-1}+k_{n-1}-k_{n}}\brc{z}.
\end{aligned}
\end{equation}

Let us pass to the case $s_1>1$. Suppose by induction that the thesis holds for a fixed $s_1$ and let us prove it for $s_1+1$. By definition
\begin{equation}
\text{ELi}_{s_1+1,s_2,\dots,s_n}\brc{z}=\int_{-\infty}^z\text{ELi}_{s_1,s_2,\dots,s_n}\brc{t}\frac{\mathrm{d}t}{t}.
\end{equation}
Using \eqref{eq:ELi_sub_gen}, we split the integral into two kinds of integrals:
\begin{itemize}
\item[1.] $\displaystyle \int_{-\infty}^z \biggl(\text{EL}_{s_1,\dots,s_n}\brc{t}+\sum_{k_1=1}^{s_1}\frac{(-1)^{k_1-1}}{(k_1-1)!\,(s_1-k_1)!}\text{cLi}_{k_1,s_{2},\dots,s_n}\,\log(-z)^{s_1-k_1}\biggr)\frac{\mathrm{d}t}{t}$,
\item[2.] $\displaystyle \int_{-\infty}^z\log\left((-1)^{i}t\right)^{s_i-k_1}\text{ELi}_{s_1+k_1-k_2,s_2+k_2-k_3,\dots,s_{i-1}+k_{i-1}-k_{i}}\brc{t}\frac{\mathrm{d}t}{t}$,\quad $i\in\{2,\dots,n\}$.
\end{itemize}
We further split the first integral into two parts:
\begin{equation}
\begin{aligned}
&\int_{-\infty}^{-1} \biggl(\text{EL}_{s_1,\dots,s_n}\brc{t}+\sum_{k_1=1}^{s_1}\frac{(-1)^{k_1-1}}{(k_1-1)!\,(s_1-k_1)!}\text{cLi}_{k_1,s_{2},\dots,s_n}\,\log(-z)^{s_1-k_1}\biggr)\frac{\mathrm{d}t}{t}=\\
&\frac{(-1)^{s_1}}{s_1!}\int_1^{\infty}\mathrm{e}^{-t}\log(t)^{s_1}\text{EL}_{s_2,\dots,s_n}(t)\frac{\mathrm{d}t}{t},
\end{aligned}
\end{equation}
and
\begin{equation}
\begin{aligned}
&\int_{-1}^{z} \biggl(\text{EL}_{s_1,\dots,s_n}\brc{t}+\sum_{k_1=1}^{s_1}\frac{(-1)^{k_1-1}}{(k_1-1)!\,(s_1-k_1)!}\text{cLi}_{k_1,s_{2},\dots,s_n}\,\log(-z)^{s_1-k_1}\biggr)\frac{\mathrm{d}t}{t}=\\
&\text{EL}_{s_1+1,s_2,\dots,s_n}(z)+\int_{-1}^0\text{EL}_{s_1,\dots,s_n}(t)\frac{\mathrm{d}t}{t}+\sum_{k_1=1}^{s_1}\frac{(-1)^{k_1-1}}{(k_1-1)!\,(s_1+1-k_1)!}\text{cLi}_{k_1,s_2,\dots,s_n}\log(-z)^{s_1+1-k_1}=\\
&\text{EL}_{s_1+1,s_2,\dots,s_n}(z)+\frac{(-1)^{s_1}}{s_1!}\int_0^{1}\mathrm{e}^{-t}\log(t)^{s_1}\text{EL}_{s_2,\dots,s_n}(t)\frac{\mathrm{d}t}{t}+\sum_{k_1=1}^{s_1}\frac{(-1)^{k_1-1}}{(k_1-1)!\,(s_1+1-k_1)!}\text{cLi}_{k_1,s_2,\dots,s_n}\log(-z)^{s_1+1-k_1}.
\end{aligned}
\end{equation}
Putting the two results together, we conclude that the first integral is given by
\begin{equation}
\begin{aligned}
\text{EL}_{s_1+1,s_2,\dots,s_n}(z)+\sum_{k_1=1}^{s_1+1}\frac{(-1)^{k_1-1}}{(k_1-1)!\,(s_1+1-k_1)!}\text{cLi}_{k_1,s_2,\dots,s_n}\log(-z)^{s_1+1-k_1},
\end{aligned}
\end{equation}
where $\text{cLi}_{s_1+1,s_2,\dots,s_n}$ is defined consistently with \eqref{eq:cLi_gen_def}.

The evaluation of the second type of integral proceeds as in the previous Appendix (in analogy with (\ref{eq:ELimn_sub_proof4})), giving the result
\begin{equation}
\sum_{l=k_1}^{s_i}(-1)^{l-k_1}\frac{(s_i-k_1)!}{(s_i-l)!}\,\log\left((-1)^{i}z\right)^{s_i-l}\text{ELi}_{s_1+l-k_2+1,s_2+k_2-k_3,\dots,s_{i-1}+k_{i-1}-k_{i}}\brc{z},
\end{equation}
for every $i\in\{2,\dots,n\}$.
Finally, reintroducing the summations in front of the second type of integrals and using the Hockey-Stick identity \eqref{eq:hockey_id} as in the previous Appendix (see identities \eqref{eq:ELi_HS1} and \eqref{eq:ELi_HS2}), the claim is proved:
\begin{equation}
\begin{aligned}
\text{ELi}_{s_1+1,s_2,\dots,s_n}\brc{z}=&\,\text{EL}_{s_1+1,s_2,\dots,s_n}(z)+\sum_{k_1=1}^{s_1+1}\frac{(-1)^{k_1-1}}{(k_1-1)!\,(s_1+1-k_1)!}\text{cLi}_{k_1,s_2,\dots,s_n}\log(-z)^{s_1+1-k_1}\\
&+\sum_{i=2}^{n-1}\sum_{s_i\geq k_1\geq \dots\geq k_i\geq 1} \frac{\brc{-1}^{k_1}}{\brc{k_i-1}!\brc{s_i-k_1}!} \tbinom{s_1+k_1-k_{2}}{s_1}\,\prod_{j=2}^{i-1} \tbinom{s_j-1+k_j-k_{j+1}}{s_j-1}\, \text{cLi}_{k_i,s_{i+1},\dots,s_n}\\
&\times \log(\brc{-1}^i z)^{s_i-k_1}\text{ELi}_{s_1+1+k_1-k_2,\dots,s_{i-1}+k_{i-1}-k_i}\brc{z}\\
&+\sum_{s_n\geq k_1\geq\dots\geq k_n\geq 0} \frac{\brc{-1}^{k_1+1}}{k_n!\brc{s_n-k_1}!}\tbinom{s_1+k_1-k_{2}}{s_1} \,\prod_{j=2}^{n-1} \tbinom{s_j-1+k_j-k_{j+1}}{s_j-1} \, \Gamma^{\brc{k_n}}\brc{1}  \times\\
&\times \log\brc{\brc{-1}^n z}^{s_n-k_1}\text{ELi}_{s_1+1+k_1-k_2,\dots,s_{n-1}+k_{n-1}-k_n}\brc{z}.
\end{aligned}
\end{equation}

\bibliographystyle{JHEP}
\bibliography{biblio}

\providecommand{\href}[2]{#2}\begingroup\raggedright\begin{thebibliography}{10}

\bibitem{Berti:2009kk}
E.~Berti, V.~Cardoso and A.~O. Starinets, \emph{{Quasinormal modes of black
  holes and black branes}},
  \href{http://dx.doi.org/10.1088/0264-9381/26/16/163001}{\emph{Class. Quant.
  Grav.} {\bf 26} (2009) 163001}, [\href{http://arxiv.org/abs/0905.2975}{{\tt
  0905.2975}}].

\bibitem{PhysRevLett.116.061102}
{\scshape LIGO Scientific Collaboration and Virgo Collaboration} collaboration,
  B.~P. Abbott, R.~Abbott, T.~D. Abbott, M.~R. Abernathy, F.~Acernese,
  K.~Ackley et~al., \emph{Observation of gravitational waves from a binary
  black hole merger},
  \href{http://dx.doi.org/10.1103/PhysRevLett.116.061102}{\emph{Phys. Rev.
  Lett.} {\bf 116} (Feb, 2016) 061102}.

\bibitem{PhysRevD.30.295}
V.~Ferrari and B.~Mashhoon, \emph{New approach to the quasinormal modes of a
  black hole}, \href{http://dx.doi.org/10.1103/PhysRevD.30.295}{\emph{Phys.
  Rev. D} {\bf 30} (Jul, 1984) 295--304}.

\bibitem{PhysRevD.35.3632}
S.~Iyer, \emph{Black-hole normal modes: A wkb approach. ii. schwarzschild black
  holes}, \href{http://dx.doi.org/10.1103/PhysRevD.35.3632}{\emph{Phys. Rev. D}
  {\bf 35} (Jun, 1987) 3632--3636}.

\bibitem{Mano:1996mf}
S.~Mano, H.~Suzuki and E.~Takasugi, \emph{{Analytic solutions of the
  Regge-Wheeler equation and the postMinkowskian expansion}},
  \href{http://dx.doi.org/10.1143/PTP.96.549}{\emph{Prog. Theor. Phys.} {\bf
  96} (1996) 549--566}, [\href{http://arxiv.org/abs/gr-qc/9605057}{{\tt
  gr-qc/9605057}}].

\bibitem{Mano:1996vt}
S.~Mano, H.~Suzuki and E.~Takasugi, \emph{{Analytic solutions of the Teukolsky
  equation and their low frequency expansions}},
  \href{http://dx.doi.org/10.1143/PTP.95.1079}{\emph{Prog. Theor. Phys.} {\bf
  95} (1996) 1079--1096}, [\href{http://arxiv.org/abs/gr-qc/9603020}{{\tt
  gr-qc/9603020}}].

\bibitem{Mano:1996gn}
S.~Mano and E.~Takasugi, \emph{{Analytic solutions of the Teukolsky equation
  and their properties}},
  \href{http://dx.doi.org/10.1143/PTP.97.213}{\emph{Prog. Theor. Phys.} {\bf
  97} (1997) 213--232}, [\href{http://arxiv.org/abs/gr-qc/9611014}{{\tt
  gr-qc/9611014}}].

\bibitem{Cardoso:2003vt}
V.~Cardoso, J.~P.~S. Lemos and S.~Yoshida, \emph{{Quasinormal modes of
  Schwarzschild black holes in four-dimensions and higher dimensions}},
  \href{http://dx.doi.org/10.1103/PhysRevD.69.044004}{\emph{Phys. Rev. D} {\bf
  69} (2004) 044004}, [\href{http://arxiv.org/abs/gr-qc/0309112}{{\tt
  gr-qc/0309112}}].

\bibitem{Konoplya:2004ip}
R.~A. Konoplya, \emph{{Quasinormal modes of the Schwarzschild black hole and
  higher order WKB approach}}, {\emph{J. Phys. Stud.} {\bf 8} (2004) 93--100}.

\bibitem{novaes}
F.~Novaes, C.~Marinho, M.~Lencs{\'e}s and M.~Casals, \emph{{Kerr-de Sitter
  Quasinormal Modes via Accessory Parameter Expansion}},
  \href{http://dx.doi.org/10.1007/JHEP05(2019)033}{\emph{JHEP} {\bf 05} (2019)
  033}, [\href{http://arxiv.org/abs/1811.11912}{{\tt 1811.11912}}].

\bibitem{novaes2014}
F.~Novaes and B.~Carneiro~da Cunha, \emph{{Isomonodromy, Painlev{\'e}
  transcendents and scattering off of black holes}},
  \href{http://dx.doi.org/10.1007/JHEP07(2014)132}{\emph{JHEP} {\bf 07} (2014)
  132}, [\href{http://arxiv.org/abs/1404.5188}{{\tt 1404.5188}}].

\bibitem{CarneirodaCunha:2015hzd}
B.~Carneiro~da Cunha and F.~Novaes, \emph{{Kerr Scattering Coefficients via
  Isomonodromy}}, \href{http://dx.doi.org/10.1007/JHEP11(2015)144}{\emph{JHEP}
  {\bf 11} (2015) 144}, [\href{http://arxiv.org/abs/1506.06588}{{\tt
  1506.06588}}].

\bibitem{Novaes:2018fry}
F.~Novaes, C.~Marinho, M.~Lencs\'es and M.~Casals, \emph{{Kerr-de Sitter
  Quasinormal Modes via Accessory Parameter Expansion}},
  \href{http://dx.doi.org/10.1007/JHEP05(2019)033}{\emph{JHEP} {\bf 05} (2019)
  033}, [\href{http://arxiv.org/abs/1811.11912}{{\tt 1811.11912}}].

\bibitem{CarneirodaCunha:2019tia}
B.~Carneiro~da Cunha and J.~a.~P. Cavalcante, \emph{{Confluent conformal blocks
  and the Teukolsky master equation}},
  \href{http://dx.doi.org/10.1103/PhysRevD.102.105013}{\emph{Phys. Rev. D} {\bf
  102} (2020) 105013}, [\href{http://arxiv.org/abs/1906.10638}{{\tt
  1906.10638}}].

\bibitem{Amado:2020zsr}
J.~B. Amado, B.~Carneiro~da Cunha and E.~Pallante, \emph{{Vector perturbations
  of Kerr-AdS$_{5}$ and the Painlev\'e VI transcendent}},
  \href{http://dx.doi.org/10.1007/JHEP04(2020)155}{\emph{JHEP} {\bf 04} (2020)
  155}, [\href{http://arxiv.org/abs/2002.06108}{{\tt 2002.06108}}].

\bibitem{Aminov:2020yma}
G.~Aminov, A.~Grassi and Y.~Hatsuda, \emph{{Black Hole Quasinormal Modes and
  Seiberg\textendash{}Witten Theory}},
  \href{http://dx.doi.org/10.1007/s00023-021-01137-x}{\emph{Annales Henri
  Poincare} {\bf 23} (2022) 1951--1977},
  [\href{http://arxiv.org/abs/2006.06111}{{\tt 2006.06111}}].

\bibitem{Hatsuda:2020sbn}
Y.~Hatsuda, \emph{{Quasinormal modes of Kerr-de Sitter black holes via the Heun
  function}}, \href{http://dx.doi.org/10.1088/1361-6382/abc82e}{\emph{Class.
  Quant. Grav.} {\bf 38} (2020) 025015},
  [\href{http://arxiv.org/abs/2006.08957}{{\tt 2006.08957}}].

\bibitem{Hatsuda:2020egs}
Y.~Hatsuda and M.~Kimura, \emph{{Semi-analytic expressions for quasinormal
  modes of slowly rotating Kerr black holes}},
  \href{http://dx.doi.org/10.1103/PhysRevD.102.044032}{\emph{Phys. Rev. D} {\bf
  102} (2020) 044032}, [\href{http://arxiv.org/abs/2006.15496}{{\tt
  2006.15496}}].

\bibitem{BarraganAmado:2021uyw}
J.~Barrag\'an~Amado, B.~Carneiro~da Cunha and E.~Pallante, \emph{{Remarks on
  holographic models of the Kerr-AdS$_{5}$ geometry}},
  \href{http://dx.doi.org/10.1007/JHEP05(2021)251}{\emph{JHEP} {\bf 05} (2021)
  251}, [\href{http://arxiv.org/abs/2102.02657}{{\tt 2102.02657}}].

\bibitem{Bonelli:2021uvf}
G.~Bonelli, C.~Iossa, D.~P. Lichtig and A.~Tanzini, \emph{{Exact solution of
  Kerr black hole perturbations via CFT2 and instanton counting: Greybody
  factor, quasinormal modes, and Love numbers}},
  \href{http://dx.doi.org/10.1103/PhysRevD.105.044047}{\emph{Phys. Rev. D} {\bf
  105} (2022) 044047}, [\href{http://arxiv.org/abs/2105.04483}{{\tt
  2105.04483}}].

\bibitem{Bianchi:2021xpr}
M.~Bianchi, D.~Consoli, A.~Grillo and J.~F. Morales, \emph{{QNMs of branes, BHs
  and fuzzballs from quantum SW geometries}},
  \href{http://dx.doi.org/10.1016/j.physletb.2021.136837}{\emph{Phys. Lett. B}
  {\bf 824} (2022) 136837}, [\href{http://arxiv.org/abs/2105.04245}{{\tt
  2105.04245}}].

\bibitem{Bianchi:2021mft}
M.~Bianchi, D.~Consoli, A.~Grillo and J.~F. Morales, \emph{{More on the SW-QNM
  correspondence}},
  \href{http://dx.doi.org/10.1007/JHEP01(2022)024}{\emph{JHEP} {\bf 01} (2022)
  024}, [\href{http://arxiv.org/abs/2109.09804}{{\tt 2109.09804}}].

\bibitem{Amado:2021erf}
J.~B. Amado, B.~C. da~Cunha and E.~Pallante, \emph{{Quasinormal modes of scalar
  fields on small Reissner-Nordstr\"om-AdS5 black holes}},
  \href{http://dx.doi.org/10.1103/PhysRevD.105.044028}{\emph{Phys. Rev. D} {\bf
  105} (2022) 044028}, [\href{http://arxiv.org/abs/2110.08349}{{\tt
  2110.08349}}].

\bibitem{Hatsuda:2021gtn}
Y.~Hatsuda and M.~Kimura, \emph{{Spectral Problems for Quasinormal Modes of
  Black Holes}},
  \href{http://dx.doi.org/10.3390/universe7120476}{\emph{Universe} {\bf 7}
  (2021) 476}, [\href{http://arxiv.org/abs/2111.15197}{{\tt 2111.15197}}].

\bibitem{Fioravanti:2021dce}
D.~Fioravanti and D.~Gregori, \emph{{A new method for exact results on
  Quasinormal Modes of Black Holes}},
  \href{http://arxiv.org/abs/2112.11434}{{\tt 2112.11434}}.

\bibitem{Bonelli:2022ten}
G.~Bonelli, C.~Iossa, D.~P. Lichtig and A.~Tanzini, \emph{{Irregular Liouville
  correlators and connection formulae for Heun functions}},
  \href{http://arxiv.org/abs/2201.04491}{{\tt 2201.04491}}.

\bibitem{Dodelson:2022yvn}
M.~Dodelson, A.~Grassi, C.~Iossa, D.~Panea~Lichtig and A.~Zhiboedov,
  \emph{{Holographic thermal correlators from supersymmetric instantons}},
  \href{http://dx.doi.org/10.21468/SciPostPhys.14.5.116}{\emph{SciPost Phys.}
  {\bf 14} (2023) 116}, [\href{http://arxiv.org/abs/2206.07720}{{\tt
  2206.07720}}].

\bibitem{Consoli:2022eey}
D.~Consoli, F.~Fucito, J.~F. Morales and R.~Poghossian, \emph{{CFT description
  of BH's and ECO's: QNMs, superradiance, echoes and tidal responses}},
  \href{http://arxiv.org/abs/2206.09437}{{\tt 2206.09437}}.

\bibitem{Imaizumi:2022qbi}
K.~Imaizumi, \emph{{Quasi-normal modes for the D3-branes and Exact WKB
  analysis}},
  \href{http://dx.doi.org/10.1016/j.physletb.2022.137450}{\emph{Phys. Lett. B}
  {\bf 834} (2022) 137450}, [\href{http://arxiv.org/abs/2207.09961}{{\tt
  2207.09961}}].

\bibitem{Ivanov:2022qqt}
M.~M. Ivanov and Z.~Zhou, \emph{{Vanishing of Black Hole Tidal Love Numbers
  from Scattering Amplitudes}},
  \href{http://dx.doi.org/10.1103/PhysRevLett.130.091403}{\emph{Phys. Rev.
  Lett.} {\bf 130} (2023) 091403}, [\href{http://arxiv.org/abs/2209.14324}{{\tt
  2209.14324}}].

\bibitem{daCunha:2022ewy}
B.~C. da~Cunha and J.~a.~P. Cavalcante, \emph{{Expansions for semiclassical
  conformal blocks}},  \href{http://arxiv.org/abs/2211.03551}{{\tt
  2211.03551}}.

\bibitem{Imaizumi:2022dgj}
K.~Imaizumi, \emph{{Exact conditions for quasi-normal modes of extremal
  M5-branes and exact WKB analysis}},
  \href{http://dx.doi.org/10.1016/j.nuclphysb.2023.116221}{\emph{Nucl. Phys. B}
  {\bf 992} (2023) 116221}, [\href{http://arxiv.org/abs/2212.04738}{{\tt
  2212.04738}}].

\bibitem{Bianchi:2022qph}
M.~Bianchi and G.~Di~Russo, \emph{{2-charge circular fuzz-balls and their
  perturbations}}, \href{http://dx.doi.org/10.1007/JHEP08(2023)217}{\emph{JHEP}
  {\bf 08} (2023) 217}, [\href{http://arxiv.org/abs/2212.07504}{{\tt
  2212.07504}}].

\bibitem{Gregori:2022xks}
D.~Gregori and D.~Fioravanti, \emph{{Quasinormal modes of black holes from
  supersymmetric gauge theory and integrability}},
  \href{http://dx.doi.org/10.22323/1.414.0422}{\emph{PoS} {\bf ICHEP2022} (11,
  2022) 422}.

\bibitem{Dodelson:2023vrw}
M.~Dodelson, C.~Iossa, R.~Karlsson and A.~Zhiboedov, \emph{{A thermal product
  formula}}, \href{http://dx.doi.org/10.1007/JHEP01(2024)036}{\emph{JHEP} {\bf
  01} (2024) 036}, [\href{http://arxiv.org/abs/2304.12339}{{\tt 2304.12339}}].

\bibitem{Bianchi:2023rlt}
M.~Bianchi, C.~Di~Benedetto, G.~Di~Russo and G.~Sudano, \emph{{Charge
  instability of JMaRT geometries}},
  \href{http://dx.doi.org/10.1007/JHEP09(2023)078}{\emph{JHEP} {\bf 09} (2023)
  078}, [\href{http://arxiv.org/abs/2305.00865}{{\tt 2305.00865}}].

\bibitem{Bianchi:2023sfs}
M.~Bianchi, G.~Di~Russo, A.~Grillo, J.~F. Morales and G.~Sudano, \emph{{On the
  stability and deformability of top stars}},
  \href{http://dx.doi.org/10.1007/JHEP12(2023)121}{\emph{JHEP} {\bf 12} (2023)
  121}, [\href{http://arxiv.org/abs/2305.15105}{{\tt 2305.15105}}].

\bibitem{Giusto:2023awo}
S.~Giusto, C.~Iossa and R.~Russo, \emph{{The black hole behind the cut}},
  \href{http://dx.doi.org/10.1007/JHEP10(2023)050}{\emph{JHEP} {\bf 10} (2023)
  050}, [\href{http://arxiv.org/abs/2306.15305}{{\tt 2306.15305}}].

\bibitem{Hatsuda:2023geo}
Y.~Hatsuda and M.~Kimura, \emph{{Perturbative quasinormal mode frequencies}},
  \href{http://dx.doi.org/10.1103/PhysRevD.109.044026}{\emph{Phys. Rev. D} {\bf
  109} (2024) 044026}, [\href{http://arxiv.org/abs/2307.16626}{{\tt
  2307.16626}}].

\bibitem{Fioravanti:2023zgi}
D.~Fioravanti and D.~Gregori, \emph{{New Developments in $\mathcal {N}=2$
  Supersymmetric Gauge Theories: from Integrability to Black Holes}},
  \href{http://dx.doi.org/10.5506/APhysPolBSupp.16.5-A31}{\emph{Acta Phys.
  Polon. Supp.} {\bf 16} (2023) 31}.

\bibitem{Saketh:2023bul}
M.~V.~S. Saketh, Z.~Zhou and M.~M. Ivanov, \emph{{Dynamical tidal response of
  Kerr black holes from scattering amplitudes}},
  \href{http://dx.doi.org/10.1103/PhysRevD.109.064058}{\emph{Phys. Rev. D} {\bf
  109} (2024) 064058}, [\href{http://arxiv.org/abs/2307.10391}{{\tt
  2307.10391}}].

\bibitem{Lei:2023mqx}
Y.~Lei, H.~Shu, K.~Zhang and R.-D. Zhu, \emph{{Quasinormal modes of C-metric
  from SCFTs}}, \href{http://dx.doi.org/10.1007/JHEP02(2024)140}{\emph{JHEP}
  {\bf 02} (2024) 140}, [\href{http://arxiv.org/abs/2308.16677}{{\tt
  2308.16677}}].

\bibitem{Bautista:2023sdf}
Y.~F. Bautista, G.~Bonelli, C.~Iossa, A.~Tanzini and Z.~Zhou, \emph{{Black hole
  perturbation theory meets CFT2: Kerr-Compton amplitudes from
  Nekrasov-Shatashvili functions}},
  \href{http://dx.doi.org/10.1103/PhysRevD.109.084071}{\emph{Phys. Rev. D} {\bf
  109} (2024) 084071}, [\href{http://arxiv.org/abs/2312.05965}{{\tt
  2312.05965}}].

\bibitem{Ivanov:2024sds}
M.~M. Ivanov, Y.-Z. Li, J.~Parra-Martinez and Z.~Zhou, \emph{{Gravitational
  Raman Scattering in Effective Field Theory: A Scalar Tidal Matching at
  O(G3)}}, \href{http://dx.doi.org/10.1103/PhysRevLett.132.131401}{\emph{Phys.
  Rev. Lett.} {\bf 132} (2024) 131401},
  [\href{http://arxiv.org/abs/2401.08752}{{\tt 2401.08752}}].

\bibitem{Arnaudo:2024rhv}
P.~Arnaudo, G.~Bonelli and A.~Tanzini, \emph{{One loop effective actions in
  Kerr-(A)dS Black Holes}},  \href{http://arxiv.org/abs/2405.13830}{{\tt
  2405.13830}}.

\bibitem{Aminov:2024mul}
G.~Aminov and P.~Arnaudo, \emph{{Black hole scattering amplitudes via analytic
  small-frequency expansion and monodromy}}, {\emph{ArXiv} (2024) },
  [\href{http://arxiv.org/abs/2409.06681}{{\tt 2409.06681}}].

\bibitem{Aminov:2023jve}
G.~Aminov, P.~Arnaudo, G.~Bonelli, A.~Grassi and A.~Tanzini, \emph{{Black hole
  perturbation theory and multiple polylogarithms}},
  \href{http://dx.doi.org/10.1007/JHEP11(2023)059}{\emph{JHEP} {\bf 11} (2023)
  059}, [\href{http://arxiv.org/abs/2307.10141}{{\tt 2307.10141}}].

\bibitem{goncharov2001multiple}
A.~B. Goncharov, \emph{Multiple polylogarithms and mixed tate motives},
  {\emph{arXiv preprint math/0103059} (2001) }.

\bibitem{Wald}
M.~Waldschmidt, \emph{Multiple polylogarithms: An introduction}, .

\bibitem{boyadzhiev2007polyexponentials}
K.~N. Boyadzhiev, \emph{Polyexponentials}, {\emph{arXiv preprint
  arXiv:0710.1332} (2007) }.

\bibitem{kim2020degenerate}
T.~Kim, D.~Kim, H.~Kim and L.~Jang, \emph{Degenerate poly-bernoulli numbers and
  polynomials}, {\emph{Informatica} {\bf 31} (2020) 2--8}.

\bibitem{Kim2019ANO}
D.~S. Kim and T.~G. Kim, \emph{A note on polyexponential and unipoly
  functions}, {\emph{Russian Journal of Mathematical Physics} {\bf 26} (2019)
  40--49}.

\bibitem{KIM2020124017}
T.~Kim and D.~S. Kim, \emph{Degenerate polyexponential functions and degenerate
  bell polynomials},
  \href{http://dx.doi.org/https://doi.org/10.1016/j.jmaa.2020.124017}{\emph{Journal
  of Mathematical Analysis and Applications} {\bf 487} (2020) 124017}.

\bibitem{komatsu1}
T.~Komatsu, \emph{Poly-cauchy numbers},
  \href{http://dx.doi.org/10.2206/kyushujm.67.143}{\emph{Kyushu Journal of
  Mathematics} {\bf 67} (07, 2012) }.

\bibitem{komatsu2}
T.~Komatsu, V.~Laohakosol and K.~Liptai, \emph{A generalization of poly-cauchy
  numbers and their properties},
  \href{http://dx.doi.org/10.1155/2013/179841}{\emph{Abstract and Applied
  Analysis} {\bf 2013} (01, 2013) }.

\bibitem{lacpao2019hurwitz}
N.~Lacpao, R.~Corcino and M.~A.~R. Vega, \emph{Hurwitz-lerch type
  multi-poly-cauchy numbers}, {\emph{Mathematics} {\bf 7} (2019) 335}.

\bibitem{Borwein1996}
J.~M. Borwein and R.~Girgensohn, \emph{Evaluation of triple euler sums.},
  {\emph{The Electronic Journal of Combinatorics [electronic only]} {\bf 3}
  (1996) Research paper R23, 27 p.--Research paper R23, 27 p.}

\bibitem{Blumlein:2009cf}
J.~Blumlein, D.~J. Broadhurst and J.~A.~M. Vermaseren, \emph{{The Multiple Zeta
  Value Data Mine}},
  \href{http://dx.doi.org/10.1016/j.cpc.2009.11.007}{\emph{Comput. Phys.
  Commun.} {\bf 181} (2010) 582--625},
  [\href{http://arxiv.org/abs/0907.2557}{{\tt 0907.2557}}].

\bibitem{em/1062621000}
D.~H. Bailey, J.~M. Borwein and R.~Girgensohn, \emph{{Experimental evaluation
  of Euler sums}}, {\emph{Experimental Mathematics} {\bf 3} (1994) 17 -- 30}.

\bibitem{Flajolet1998EulerSA}
P.~Flajolet and B.~Salvy, \emph{Euler sums and contour integral
  representations}, {\emph{Exp. Math.} {\bf 7} (1998) 15--35}.

\bibitem{xu2020explicit}
C.~Xu and W.~Wang, \emph{Explicit formulas of euler sums via multiple zeta
  values}, {\emph{Journal of Symbolic Computation} {\bf 101} (2020) 109--127}.

\bibitem{XU2017443}
C.~Xu, \emph{Multiple zeta values and euler sums},
  \href{http://dx.doi.org/https://doi.org/10.1016/j.jnt.2017.01.018}{\emph{Journal
  of Number Theory} {\bf 177} (2017) 443--478}.

\bibitem{CHEN2015107}
K.-W. Chen, \emph{Applications of stuffle product of multiple zeta values},
  \href{http://dx.doi.org/https://doi.org/10.1016/j.jnt.2015.01.003}{\emph{Journal
  of Number Theory} {\bf 153} (2015) 107--116}.

\bibitem{kuba2019multisets}
M.~Kuba, \emph{On multisets, interpolated multiple zeta values and limit laws},
  {\emph{arXiv preprint arXiv:1903.07346} (2019) }.

\bibitem{hoffman2021logarithmic}
M.~E. Hoffman and M.~Kuba, \emph{Logarithmic integrals, zeta values, and tiered
  binomial coefficients}, {\emph{Monatshefte f{\"u}r Mathematik} {\bf 195}
  (2021) 119--154}.

\bibitem{hoffman2015quasi}
M.~E. Hoffman, \emph{Quasi-symmetric functions and mod p multiple harmonic
  sums}, {\emph{Kyushu Journal of Mathematics} {\bf 69} (2015) 345--366}.

\bibitem{seki2020ohno}
S.-i. Seki and S.~Yamamoto, \emph{Ohno-type identities for multiple harmonic
  sums}, {\emph{Journal of the Mathematical Society of Japan} {\bf 72} (2020)
  673--686}.

\bibitem{Si+2021+1612+1619}
X.~Si, \emph{Euler-type sums involving multiple harmonic sums and binomial
  coefficients}, \href{http://dx.doi.org/doi:10.1515/math-2021-0124}{\emph{Open
  Mathematics} {\bf 19} (2021) 1612--1619}.

\bibitem{kuba2019note}
M.~Kuba and A.~Panholzer, \emph{A note on harmonic number identities, stirling
  series and multiple zeta values}, {\emph{International journal of number
  theory} {\bf 15} (2019) 1323--1348}.

\bibitem{batir2017some}
N.~Batir, \emph{On some combinatorial identities and harmonic sums},
  {\emph{arXiv preprint arXiv:1703.06401} (2017) }.

\bibitem{muneta2007some}
S.~Muneta, \emph{On some explicit evaluations of multiple zeta-star values},
  {\emph{arXiv preprint arXiv:0710.3219} (2007) }.

\bibitem{Lewin'81}
L.~Lewin, \emph{Polylogarithms and associated functions}.
\newblock Elsevier North Holland, Inc., 1981.

\bibitem{Goncharov:1998kja}
A.~B. Goncharov, \emph{{Multiple polylogarithms, cyclotomy and modular
  complexes}}, \href{http://dx.doi.org/10.4310/MRL.1998.v5.n4.a7}{\emph{Math.
  Res. Lett.} {\bf 5} (1998) 497--516},
  [\href{http://arxiv.org/abs/1105.2076}{{\tt 1105.2076}}].

\bibitem{MINH2000217}
H.~N. Minh, M.~Petitot and J.~V.~D. Hoeven, \emph{Shuffle algebra and
  polylogarithms},
  \href{http://dx.doi.org/https://doi.org/10.1016/S0012-365X(00)00155-2}{\emph{Discrete
  Mathematics} {\bf 225} (2000) 217--230}.

\bibitem{MINH2000273}
H.~N. Minh and M.~Petitot, \emph{Lyndon words, polylogarithms and the riemann
  zeta function},
  \href{http://dx.doi.org/https://doi.org/10.1016/S0012-365X(99)00267-8}{\emph{Discrete
  Mathematics} {\bf 217} (2000) 273--292}.

\bibitem{HOFFMAN1997477}
M.~E. Hoffman, \emph{The algebra of multiple harmonic series},
  \href{http://dx.doi.org/https://doi.org/10.1006/jabr.1997.7127}{\emph{Journal
  of Algebra} {\bf 194} (1997) 477--495}.

\end{thebibliography}\endgroup

\end{document}